\documentclass[11pt,reqno]{amsart}

\usepackage{amsmath,epsfig}
\usepackage{amssymb}
\usepackage{booktabs}
\usepackage{graphicx}
\graphicspath{{.}{./figures/}}

\usepackage{natbib}
\usepackage{hypernat}
\usepackage[dvips,ps2pdf,pagebackref=false,colorlinks=false,citecolor=red,linkcolor=blue,plainpages=false,breaklinks=true,pdfpagelabels=true]{hyperref}
\usepackage[all]{hypcap}

\numberwithin{equation}{section}

\def\1{\mbox{1\hspace{-.35em}1}} 
\def\R{\mathbb{R}}

\def\N{\mathbb{N}}
\def\P{\mathbb{P}}
\def\E{\mathbb{E}}
\def\L{\mathbb{L}}

\def\R{\mathbb{R}}

\def\Z{\mathbb{Z}}

\def\Lip{\mbox{Lip\,}}
\def\cov{\mbox{Cov}}%

\newtheorem{rem}{Remark}
\newtheorem{theo}{Theorem}[section]
\newtheorem{lem}{Lemma}[section]
\newtheorem{prop}{Proposition}[section]

\begin{document}

\title[Adaptive density estimation under weak dependence]
{Adaptive density estimation under weak dependence}

\author[I. Gannaz]{Gannaz Ir\`{e}ne$^{(1)}$}
\address[1]{Laboratoire Jean Kuntzmann, Grenoble-INP,
BP 53,
38041 Grenoble Cedex 9, France}

\email{irene.gannaz@imag.fr}

\author[O. Wintenberger]{Wintenberger Olivier$^{(2)}$}
\address[2]{SAMOS-MATISSE 
(Statistique Appliqu\'ee et MOd\'elisation Stochastique)
Centre d'\'Economie de la Sorbonne
Universit\'e Paris 1 - Panth\'eon-Sorbonne, CNRS
90, rue de Tolbiac
75634 Paris Cedex 13, France}

\email{olivier.wintenberger@univ-paris1.fr}


\medskip

\begin{abstract}
Assume that $(X_t)_{t\in\Z}$ is a real valued time series
admitting a common marginal density $f$ with respect to Lebesgue's measure. Donoho {\it et al.}~(1996) propose  near-minimax estimators $\widehat f_n$ based on thresholding wavelets to estimate $f$ on a compact set in an independent and identically distributed setting. The aim of the present work is to extend these results to general weak dependent contexts. Weak dependence assumptions are expressed as decreasing bounds of covariance terms and are detailed for different examples. The threshold levels in estimators $\widehat f_n$ depend on weak dependence properties of the sequence $(X_t)_{t\in\Z}$ through the constant. If these properties are unknown, we propose cross-validation procedures to get new estimators. These procedures are illustrated via simulations of dynamical systems and non causal infinite moving averages. We also discuss the efficiency of our estimators with respect to the decrease of covariances bounds. 
\end{abstract}

\maketitle

\section{Introduction}
\label{intro}
Let $(X_t)_{t\in\Z}$ be a real valued time series admitting a common marginal density
$f$ that is compactly supported. The general purpose of this paper is to estimate $f$ by wavelet estimators $\widehat f_n$ constructed from $n$ observations $(X_1,\dots,X_n)$. In their seminal paper Donoho {\it et al.}~(1996)~\cite{DonJohnKerkPic} showed that projection-like linear estimators are not optimal: introduction of nonlinearity via thresholds of wavelet coefficients is investigated. Wavelets thresholding provides estimators which adapt themselves to the unknown smoothness of $f$, we refer to Vannucci (1998)~\cite{Vannucci} for a survey of the use of wavelet bases in density estimation. The present work extends near minimax results of soft and hard-threshold estimators from the independent and identically distributed (iid for short) framework, see Theorem 5 in Donoho {\it et al.}~(1996)~\cite{DonJohnKerkPic}, to cases where weak dependence between variables occurs.

Our main assumptions give bounds for covariance terms as decreasing sequences which tend to zero when the gap between the past and the future of the time series goes to infinity. In order to give examples satisfying these conditions, we introduce coefficients that give bounds of covariance terms and that are computable for a large class of models. Weak dependent coefficients, introduced by Doukhan and Louhichi (1999) \cite{DoukhanLouhichi}, as mixing ones, are well-adapted to that purpose. Using $\beta$-mixing coefficients Tribouley and Viennet (1998)~\cite{TribouleyViennet} proposed minimax estimators with respect to the Mean Integrated Square Error (MISE for short). Comte and Merlev\`ede (2002)~\cite{ComteMerlevede} obtained near-minimax results using $\alpha$-mixing coefficients. The loss of a logarithmic factor in the convergence rate in this last paper is balanced by the generality of the context, as the class of $\alpha$-mixing models is larger than the one of $\beta$-mixing models.

However, $\alpha$ and $\beta$-mixing coefficients are not easy to compute for some models and are useless for others; Andrews~(1984)~\cite{Andrews} proved that the mixing coefficients of the stationary solution of the AR($1$) model
\begin{equation}\label{andrews}X_t=\frac{1}{2}\left(X_{t-1}+\xi_t\right),\text{ where }(\xi_t)_{t\in\Z}\text{ iid with a Bernouilli law of parameter } 1/2,\end{equation}
do not tend to zero as the gap from the past to the future of the time series goes to infinity. Mixing coefficients do not behave nicely in this case as, through a reversion of time, the Markov chain solution of \eqref{andrews} is a dynamical system, i.e. $X_{t-1}=T(X_t)$ for some transformation $T$, namely $T(x)=2x\1_{0\le x<1/2}+ (2x-1)\1_{1/2\le x\le 1}$. So called weak dependence coefficients have been recently developed to deal with such processes, see Dedecker {\it et al.} \cite{Dedecker2007}, Maume-Deschamps \cite{Maume-Deschamps2006} and references therein.
Introduced by Dedecker and
Prieur (2005) in~\cite{DedeckerPrieur}, $\widetilde \phi$-weak dependence coefficients give sharp bounds on the covariance terms of dynamical systems, such as the stationary solution to \eqref{andrews}. Using these coefficients, we prove near-minimax results of thresholded wavelet estimators for dynamical systems called expanding maps. To our knowledge, only non adaptive density estimation has been studied in this non-mixing context, see for instance Bosq et Guegan \cite{Bosq1995}, Prieur \cite{Prieur} and Maume-Deschamps \cite{Maume-Deschamps2006}.\\ 

The advantage of our approach is also to treat in one draw many other contexts of dependence. We prove that near-minimaxity still holds for a very large class of models using $\lambda$-weak dependence coefficients, defined by Doukhan and
Wintenberger,~2007,~\cite{Doukhan2007}. We pay for generality by adding up conditions on the joint densities of the couples $(X_0,X_r)$ for all $r>0$. These conditions are not restrictive as it is satisfied for many econometric models such as ARMA, GARCH, ARCH, LARCH, MA models.\\

The estimation scheme is based on Donoho {\it et al.}'s
procedures developed in~\cite{DonJohnKerkPic} for the iid case, and it is adaptive with respect to the regularity of $f$. Soft and hard-threshold
levels $(\lambda_j)_{j_0\le j\le j_1}$ are chosen equal to $K\sqrt{j/n}$ for some $K>0$. Note that the constant $K$ and the highest resolution level $j_1$ depend on the weak dependence properties of the observations. If weak dependence properties are known, estimators are near minimax: same rates as in the iid setting are achieved for  mean-$L^p$ errors with $1\le p<\infty$. If weak dependence properties are unknown, we develop cross-validation procedures to approximate threshold levels $\hat\lambda_j$ and the highest resolution level $\widehat j_1$.
We check on simulations that corresponding estimators are adaptive with respect to the regularity of $f$ and to  weak dependence properties of $(X_t)_{t\in\Z}$. The order of errors of approximations in the dependence cases are very close to the one in the iid cases.\\

We believe that we obtain such good results as we work on simulations of processes satisfying our main Assumption {\bf (D)}. This assumption consists of the exponential decay of the covariance terms. We give in this paper some simulations study of dynamical systems that do not satisfy this Assumption {\bf (D)}. Comparing the behavior of our estimators with the kernel ones show that ours are less efficient when covariance terms do not decrease exponentially fast. We also prove that the error terms of our estimators are unbounded in some cases, due to terms of covariances that decrease too slowly. This is a restriction of our procedure, based on thresholding wavelet coefficients.\\

The paper is structured as follows. In Section \ref{Preliminaries}, we give notation, we introduce estimation procedures and we formulate weak dependence assumptions. Main results are given in Section \ref{main} and examples of models in Section~\ref{model}. Cross-validation procedures and their accuracy on simulations are developed in Section~\ref{simus}. Proofs are relegated in the last Section.

\section{Preliminaries}\label{Preliminaries}

\subsection{Notation}
We restrict ourselves to the estimation of a density $f$ which is compactly supported. We suppose that $f$ is supported by $[-B,B]$ for some $B>0$. For all $p\ge1$, $L^p$ denotes the space of all functions $f$ such that $\|f\|^p_{p}=\int|f(x)|^pdx<\infty$.

\subsection{Density estimators}
Throughout the paper, we work within an $r$-regular orthonormal multiresolution analysis of $L^2$ (endowed with the usual inner product), associated with a compactly supported scaling function $\phi$ and a compactly supported mother wavelet $\psi$. Without loss of generality, we suppose that the support of functions $\phi$ and $\psi$ is included in an interval $[-A,A]$ for some $A>0$. Let us recall that $\phi$ and $\psi$ generate orthonormal basis by dilatations and translations: for a given primary resolution level $j_0$, the functions $\{\phi_{j,k}:x\mapsto 2^{j/2}\phi(2^{j}x-k) \}_{k\in\Z}$ and $\{\psi_{j,k}:x\mapsto 2^{j/2}\psi(2^{j}x-k)\}_{k\in\Z}$ are such that the family 
$$
\{\phi_{j_0,k},\,k\in\Z,\, \psi_{j,k}\,j\geq j_0,\,k\in\Z\}
$$
is an orthonormal base of $L^2$. Any function $f\in L^2$ can thus be decomposed as 
$$
f=\sum_{k\in\Z}\alpha_{j_0,k}\phi_{j_0,k}+\sum_{j=j_0}^{\infty}\sum_{k\in\Z}\beta_{j,k}\psi_{j,k},
$$ 
where $\alpha_{j,k}=\int f(x)\phi_{j,k}(x) dx$, $\beta_{j,k}=\int f(x)\psi_{j,k}(x) dx$.

The nonlinear estimator developed in \cite{DonJohnKerkPic} is defined by the equation
$$\widehat f_n=\sum_{k\in\Z}\widehat\alpha_{j_0,k}\phi_{j_0,k}+
\sum_{j=j_0}^{j_1}\sum_{k\in\Z}\gamma_{\lambda_j}(\widehat\beta_{j,k})\psi_{j,k},$$ where $\hat\alpha_{j,k}= n^{-1}\sum_{i=1}^n \phi_{j,k}(X_i)$ and
$\hat\beta_{j,k}=n^{-1}\sum_{i=1}^n \psi_{j,k}(X_i)$, and where $\gamma_\lambda$ is a threshold function of level $\lambda$. The authors consider both hard and soft thresholding functions, corresponding respectively to $\gamma_\lambda(\beta)=\beta\1_{|\beta|>\lambda}$ and $\gamma_\lambda(\beta)=(|\beta|-\lambda)_+\mbox{sign}(\beta)$. If no distinction is done in the sequel, both hard and soft thresholding estimators are concerned. 

Let $s,\pi,r$ be three positive real numbers satisfying $s+1/2-1/\pi>0$. We assume that $f$ belongs to the Besov ball $\mathcal B^s_{\pi,r}(M_1)$ on the real line, {\it i.e.} $\|f\|_{s,\pi,r}\le M_1$ where $$\|f\|_{s,\pi,r}=|\alpha_{0,0}|+\left(\sum_{j\in\N}\left(2^{j(s\pi+\pi/2-1)}\sum_{k\in\Z}|\beta_{j,k}|^\pi\right)^{r/\pi}\right)^{1/r}.$$
Approximation errors of an estimator $f_n$ are expressed as
$\E\|f_n-f\|_{p}^p$ for $p\ge 1$. Associated minimax rates are the best convergence decrease $\alpha$ of the worst approximation error we may achieve over all estimators $f_n$:
$$
\inf_{f_n} \sup_{f\in \mathcal B^s_{\pi,r}(C)}\E\|f_n-f\|_{p}^p\,=\,\mathcal{O}\left(n^{-p\alpha}\right).
$$ 
The minimax rate $\alpha$ is determined in \cite{DonJohnKerkPic} as:
\begin{equation}\label{minimaxrate}
\alpha=\begin{cases}\alpha_+=(s/(1+2s)
&\mbox{if }\epsilon\ge 0,\\
\alpha_-=(s-1/\pi+1/p)/(1+2s-2/\pi)&\mbox{if }\epsilon\le 0,\end{cases}\qquad\mbox{ where }\epsilon=s\pi-(p-\pi)/2.\end{equation}

\subsection{Weakly dependent assumption}
Throughout the paper, the symbol $\delta$ denotes with no distinction $\phi$ or $\psi$ and $\widetilde\delta_{j,k}(x)=\delta_{j,k}(x)-\E\,\delta_{j,k}(X_0)$ for all integers $j\ge0$ and $k$. Define for all positive integers $u,v$ the quantities
\begin{equation}\label{ckr}
C_{u,v}^{j,k}(r)=\sup_{\max
s_{i+1}-s_i=s_{u+1}-s_u=r}\{|\cov(\widetilde\delta_{j,k}(X_{s_1})\cdots
\widetilde\delta_{j,k}(X_{s_u}),
\widetilde\delta_{j,k}(X_{s_{u+1}})\cdots
\widetilde\delta_{j,k}(X_{s_{u+v}}))|\}.
\end{equation}
Functions $\phi$ and $\psi$ play a symmetric role through $\delta$ in this setting. As stressed in \cite{DoukhanLouhichi}, bounds on covariance terms $C_{u,v}^{j,k}(r)$ are useful to extend asymptotic results from the iid case. Now we can state the main assumption of this paper:
\begin{description}
\item[(D)] There exists a sequence $\rho(r)$ such that for all $r\ge 0$, all indexes $j,k,u,v$, we have
\begin{equation}\label{cpr}\tag{D1}
C_{u,v}^{j,k}(r)\le (u+v+uv)/2\, (2^{j/2}M_2)^{u+v-2} \rho(r)\mbox{ where }M_2 \mbox{ is a constant satisfying }\|\delta\|_\infty\le M_2. \end{equation}
Moreover, there exist real numbers $a,b,C_0>0$ depending only on $\delta$, $f$ and on the dependence properties of $(X_t)_{t\in\Z}$ such that
\begin{equation}\label{condrho}\tag{D2}
\rho(r)\le C_0 \exp(-ar^b)\mbox{ for
all }r>0.
\end{equation}
\end{description}
In Section \ref{model}, we give explicit conditions on the stationary process $(X_t)_{t\in\Z}$ in order that it satisfies Assumption~{\bf (D)}. When it is possible, the values of the constants $a,b,C_0$ and $M_2$ are given.

\section{Main results}\label{main}

Let $(X_t)_{t\in\Z}$ be a stationary real valued time series. \vspace{-0.5cm}
\subsection{A useful lemma}\label{theoric}
Under similar conditions than Assumption~{\bf (D)}, moments inequalities of even  orders and Bernstein's type inequalities for the sums $\sum_{i=1}^n \widetilde \delta_{j,k}(X_i)$ are respectively given in \cite{DoukhanLouhichi} and \cite{DoukhanNeumann}. The following Lemma recall these inequalities applied on the quantities $\beta_{j,k}$. They remain valid for $\alpha_{j,k}$ as $\phi$ and $\psi$ play the same role in {\bf (D)}.
\begin{lem}\label{useful}
If {\bf (D)} holds then for all even integer $q\ge 2$, for all $\lambda\ge0$ and for $j,k$ such that $0\le j\le \log n$ and $0\le k\le 2^{j}-1$:
\begin{eqnarray}\label{Rosenthal}
&&\E |\widehat\beta_{j,k}-\beta_{j,k}|^q\le C_1
n^{-q/2},\\ 
\label{bernstdep}&& P\left(
|\widehat\beta_{j,k}-\beta_{j,k}|\ge\lambda\right)\leq 2\exp\left(
- \frac{n\lambda^{2}}{C_2+ C_3 (2^j/n)^{b/(2(1+2b))}
(\sqrt n\lambda)^{(2+3b)/(1+2b)}} \right),
\end{eqnarray}
where $C_1$, $C_2$ and $C_3$ are constants depending on $q$ and on the constants of Assumption~{\bf (D)}: $a$, $b$, $C_0$ and $M_2$.
\end{lem}
The proof of this Lemma is given in Section \ref{prooflem}. Notice that $C_3$ depends deeply on the dependence context through $C_0$, see Propositions \ref{phiD} and \ref{etaD} for more details.

\subsection{Near-minimax results of thresholded wavelet estimators}\label{res}
Following results are extensions to weak dependence settings of Theorem 5 of \cite{DonJohnKerkPic}.

\begin{theo}\label{th}
Suppose that $f\in\mathcal B^s_{\pi,r}(M_1)$
with $1/\pi<s< N/2$ where $N$ is the regularity of the function $\psi$. If {\bf (D)} holds, then for each $1\le p<\infty$ there
exists a constant $C(N,p,a,b,C_0,M_1,M_2,A,B)$ such that
$$\E[\|\widehat f_n -f\|_p^p]\le C \begin{cases}
\displaystyle \left(\frac{\log n}n\right)^{p\alpha} & \text{if}~\epsilon\neq0\\
\displaystyle \left(\frac{\log n}n\right)^{p\alpha} (\log n)^{(p/2-(1\wedge \pi)/r)_+} & \text{if}~\epsilon=0\\
\end{cases}$$
where the minimax rate $\alpha$ and the parameter $\epsilon$ are given in (\ref{minimaxrate}). Here
$j_0$ is chosen as the smallest integer larger than $\log(n)(1+N)^{-1}$, $j_1$ is the largest integer smaller than $\log( n \log^{-2/b -3} (n))$ and $\lambda_j=K\sqrt{j/n}$ for a sufficiently large
constant $K>0$.
\end{theo}
A sketch of the proof of this Theorem is given in Section \ref{proofth}.

We refer to \cite{Meyer} for the definition of the parameter $N$, the regularity of the wavelet function $\psi$. The condition $s<N/2$ ensures the sparsity of the wavelet coefficients of the density~$f$. This condition is not restrictive as $N$ can be chosen sufficiently large in practice.

The estimators $\widehat f_n$ are the same than in the iid case given in \cite{DonJohnKerkPic}, except for the highest resolution level $j_1$ which is smaller here. This restriction is needed in the weak dependence context due to the Bernstein's type inequality \eqref{bernstdep} which is not as sharp as the one in the iid case. But this restriction does not perturb the rate of convergence which is the same as the one obtained in the iid case by \cite{DonJohnKerkPic}.

The constant $K>0$ plays a key role in the asymptotic behavior of $\widehat f_n$. From \eqref{bernstdep}, we infer that the constant $K$ depends on the parameters $a,b,C_0,M_2$ through $C_1,C_2,C_3$ in an intricate way. Then how sufficiently large the constant $K$ must be deeply depends on the dependence structure of observations $(X_1,\ldots,X_n)$.
Then contrarily to the iid case, we are not able to develop direct procedures based on the observations $(X_1,\ldots,X_n)$ that chose a convenient parameter $K$ like in \cite{Juditsky2004}. In Section \ref{simus}, we propose cross validation procedures to determine the threshold levels $\lambda_j$ when the weak dependence properties of the process $(X_t)_{t\in\Z}$ are unknown.

\section{Examples}\label{model}

In this Section, we give examples of models that satisfy Assumption~{\bf (D)} using $\widetilde\phi$ and $\lambda$-weak dependence coefficients introduced respectively in \cite{DedeckerPrieur} and \cite{Doukhan2007}. For each of them, we proceed in two steps: Firstly we give sufficient conditions on these coefficients for ensuring Assumption~{\bf (D)} in Propositions \ref{phiD} and \ref{etaD} and secondly we give in Subsections \ref{exphi} and \ref{exeta} examples of models satisfying such conditions. 

A Lipschitz function $h:\R^{u}\to \R$ for some
$u\in\N^\ast$ is a function such that $\mbox{Lip
\!}(h)<\infty$ with
\begin{eqnarray*}
\mbox{Lip \!}(h)&=&\sup_{(a_1,\ldots,a_u)\ne
(b_1,\ldots,b_u)}\frac{\left|h(a_1,\ldots,a_u)- h
(b_1,\ldots,b_u)\right|}{|a_1-b_1|+\cdots+|a_u-b_u|}.
\end{eqnarray*}
As Lipschitz functions play an important role in weak dependence contexts we restrict ourselves to the cases where $N>4$. This assumption on the regularity of the wavelet functions implies that $\phi$ and $\psi$ can be chosen as Lipschitz functions, as established in \cite{Daubechies}. Note also that $\|f\|_\infty<\infty$ as $f\in\mathcal B^s_{\pi,r}(M_1)$, see equation (15) in \cite{DonJohnKerkPic}. For convenience, we denote with no distinction $\psi_{j,k}(x)-\beta_{j,k}$ and $\phi_{j,k}(x)-\alpha_{j,k}$ as $\widetilde\delta_{j,k}(x)$ for any integers $j\ge 0$ and $k$.\\

Weak dependence coefficients is to generalize mixing ones. Let us recall that $\alpha$-mixing coefficients can be defined in a similar way by two equations: 
\begin{eqnarray*}
\alpha(r)&=&\sup_{\tiny\begin{matrix}
	1\le \ell \\
i+r\le j_1\le\cdots\le j_\ell\\
\|g\|_\infty\leq 1\end{matrix}} \E|\E\left(g(X_{i_1},\dots,X_{i_u})|\sigma(\{X_j,j\le i\})\right) -\E\left(g(X_{i_1},\dots,X_{i_u})\right)|,\\
\alpha(r)&=&\sup_{\tiny\begin{matrix}
(u,v)\in \N^\ast\times \N^\ast \\
i_1\le\cdots\le i_u\le
i_u+r\le i_{u+1}\le\cdots\le i_{u+v}\\
\|f\|_\infty,\|g\|_\infty\leq 1\end{matrix}}|\cov(f(X_{i_1},\dots,X_{i_u}),g(X_{i_{u+1}},\dots,X_{i_{u+v}}))|.
\end{eqnarray*}
These coefficients measure the dependence between the past and the future values of the process $(X_t)_{t\in\Z}$ as the gap $r$ between past and future goes to infinity. As these coefficients are often too restrictive, the authors of \cite{Dedecker2007a} release them by considering a supremum taken on functions with bounded variations or on Lipschitz bounded functions rather than uniformly bounded functions. These different choices of functions sets lead to different coefficients of weak dependence, namely $\widetilde \phi$ and $\lambda$ respectively. We give hereafter the precise definition of $\widetilde\phi$ and $\lambda$-weak dependence coefficients.

\subsection{\texorpdfstring{$\widetilde\phi$-weak dependence}{(phi tilde)-weak dependence}}\label{secphi} 

Bounded variations functions are defined as follows: let $BV$ and $BV_1$ denote the sets of functions $g$ supported on
$[-A,A]$ satisfying respectively $\|g\|_{BV}<+\infty$ and $\|g\|_{BV}\le 1$
where
$$\|g\|_{BV}=|g(-A)|+\sup_{n\in\N}\sup_{a_0=-A<a_1<\cdots<a_n=A}\sum_{i=1}^n|g(a_{i})-g(a_{i-1})|.$$
Let $(\Omega,\mathcal A,\mathbb{P})$ be a probability space,
$\mathcal M$ a $\sigma$-algebra of $\mathcal A$ and  let $X=(X_1,\ldots,X_v)$, for $v\geq 1$, be a collection of real valued random variables $X_i$ defined on $\mathcal A$. We define the coefficient $\widetilde \phi$ as it was introduced in \cite{Dedecker2007a} by the equation:
$$
\widetilde \phi({\mathcal M},X_1,\ldots,X_v)=\sup_{g_1,\ldots,g_v \in
BV_1}\left\|\int\prod_{i=1}^v g_i(x_i)\P_{X|{\mathcal
M}}(dx)-\int\prod_{i=1}^v g_i(x_i)\P_{X}(dx)\right\|_\infty,
$$
where $dx=(dx_1,\ldots,dx_v)$.
The coefficients $\widetilde \phi(r)$ are now defined by the equation
$$
\widetilde \phi (r)=\max_{1\le \ell}\,\frac{1}{\ell}\sup_{i+r\le
j_1\le\cdots \le j_\ell}\widetilde \phi(\sigma(\{X_j; j\le i\}),X_{j_1},\ldots,X_{j_\ell})\;.
$$
These coefficients are multivariate extensions of $\widetilde\phi_1(r)$ defined in \cite{DedeckerPrieur}, see also \cite{Maume-Deschamps2006}. Instead of mixing coefficients, they efficiently treat the dependence structure of dynamical systems and associated Markov chains.
A process $(X_t)_{t\in\Z}$ is said to be $\widetilde\phi$-weakly dependent if the series $\widetilde\phi(r)$ goes to $0$ when $r$ goes to infinity, {\it i.e.} when the gap between observations from the future and observations from the past goes to infinity. Introducing $\widetilde\phi$-weakly dependent processes is useful in the present framework due to the following links with the Assumption {\bf (D)}:

\begin{prop}\label{phiD}
Assume that $(X_t)_{t\in\Z}$ is a process such that there exist $a,b,c>0$ satisfying
\begin{equation}\label{condphi}
\widetilde\phi(r)\le c^v \exp(-ar^b), 
\end{equation}
then Assumption~{\bf (D)} holds with $M_2=2(\|\delta\|_\infty+A\Lip\delta)$ and $C_0=4c(\|\delta\|_\infty+A\Lip\delta)\|f\|_\infty\|\delta\|_1$.
\end{prop}
The proof of this Proposition is given in Subsection \ref{proofsimus}.

\subsection{\texorpdfstring{Examples of $\widetilde\phi$-weakly dependent processes}{Examples of (phi tilde)-weakly dependent processes}}\label{exphi}
Following the work of \cite{Dedecker2007a}, we give a general class of models where \eqref{condphi} is satisfied.
\begin{lem}\label{tildephi}
Assume that $(X_t)_{t\in\Z}$ is a process satisfying the Markov property, taking values in $[0,1]$ and such that there exist constants $a,b,c>0$ satisfying, for any functions $g,k$ with $g\in BV_1$ and $\E |k(X_0)|<\infty$, and for all $r\ge0$, 
\begin{eqnarray}
\label{condan}|\cov(k(X_0),g(X_r))|&\le& \E |k(X_0)|\exp(ar^{-b}),\\
\label{condC'}\left\|\E(g(X_r|X_0=\cdot)\right\|_{BV}&\le& c,
\end{eqnarray}
then $\widetilde \phi_{v}(r)\le  c^v \exp(-ar^b)$ and the conclusions of Proposition \ref{phiD} follow.
\end{lem}
The proof of this Lemma is given in Subsection \ref{proofsimus}.
Various examples of processes satisfying conditions of Lemma \ref{tildephi} are given in \cite{Dedecker2007a}. We recall here the case of Markov chains obtained by time reversing expanding maps, as they are extensions of Andrew's example \eqref{andrews}. 

Let us define stationary Markov chains $(X_t)_{t\in \Z}$ associated with dynamical systems through a reversion of the time as non degenerate stationary solutions of the recurrent equation 
\begin{equation}\label{sysdyn}X_t=T^i(X_{t-i}),\qquad\forall
t\in\Z,i\in\N\end{equation} 
where $T:[0,1]\to [0,1]$ is a deterministic function.
\begin{rem}
For such Markov chains, mixing coefficients are useless. Future values write simply as functions of past values via \eqref{sysdyn} and then it is easy to check that $\alpha(r)$ are constant from their definitions. Thus $\alpha(r)$ does not tend to $0$ and $\alpha$-mixing coefficients do not evaluate the dependence of such processes.
\end{rem}
A Markov chain $(X_t)_{t\in \Z}$ is associated with an expanding map through a reversion of the time if $T$ 
satisfies
\begin{itemize}
\item\textit{(Regularity)} The function $T$ is differentiable, with a continuous derivate $T'$ and there exists a grid $0=a_0\le a_1\cdots\le a_k=1$
such that $|T'(x)|>0$ on $]a_{i-1},a_i[$ for each $i=1,\dots,k$.
\item\textit{(Expansivity)} For any integer $i$, let $I_i$ be the set on which the first derivate of $T^i$, $(T^i)'$, is defined.
There exists $a>0$ and $s>1$ such that $\displaystyle{\inf}_{x\in
I_i}\{|(T^i)'(x)|\}>as^i$.
\item\textit{(Topological mixing)} For any nonempty open sets $U$, $V$, there
exists $i_0\ge 1$ such that $T^{-i}(U)\cap V\neq\varnothing$ for all
$i\ge i_0$.
\end{itemize}
Under these three conditions a non degenerate stationary solution $(X_t)_{t\in\N}$ to \eqref{sysdyn} exists and has remarkable properties, see \cite{Viana} for a nice survey. For instance, the process $(X_t)_{t\in\Z}$ satisfies the conditions of Proposition~\ref{tildephi} for some $a,c>0$ and $b=1$, see \cite{DedeckerPrieur}. Moreover, the marginal distribution is absolutely regular and its distribution belongs to $BV$. Noticing that $\mathcal B^1_{1,1}\subset
BV\subset \mathcal B^1_{1,\infty}$, see e.g. \cite{DonJohnKerkPic}, Theorem~\ref{th} provides adaptive estimators of the marginal density of $(X_t)_{t\in\Z}$ in that context and extends Theorem~2.2 of \cite{Vanharen2006} where rates of the MISE for non adaptive estimators are given in such context.

\subsection{\texorpdfstring{$\lambda$-weak dependence}{lambda-weak dependence}}\label{seceta}
The stationary process $(X_t)_{t\in\Z}$ is
$\lambda$-weakly dependent, as defined in \cite{Doukhan2007}, if there exists a sequence of non-negative
real numbers $\lambda(r)$ satisfying $\lambda(r)\to0$ as
$r\to\infty$ and such that:
\begin{multline*}
\left|\mbox{Cov}\left(\,h\left(X_{i_1}, \ldots X_{i_u}\right)\, , \,
k\left(X_{i_{u + 1}},\ldots, X_{i_{u+v}}\,\right)\right)\right|
\leq \\
(u\,\|k\|_\infty{\rm Lip}(h)\,+\,v\|h\|_\infty\,{\rm Lip}(k)+uv\Lip(h)\Lip(k))\; \lambda(r)
\end{multline*} 
for all $p$-tuples, $(i_1, \ldots
,i_p)$ with $i_1 \leq \cdots \leq i_u \le i_u+r \leq i_{u+1}\le
\cdots \le i_p$, and for all $h\in\Lambda_u$ and $h\in\Lambda_v$ where
\begin{eqnarray*}
{\Lambda}_u &= &\left\{ 
h:\R^{u}\to\ \R, \; \Lip (h)<\infty, \;\|h\|_\infty<\infty
\right\},\mbox{ for any }u\le1.
\end{eqnarray*}
The $\lambda$-weak dependence provides simple bounds of the covariance
terms:
$$
\left| \mbox{Cov} \left(\widetilde\delta_{j,k}(X_{i_1})\cdots \widetilde\delta_{j,k}(X_{i_u}), \widetilde\delta_{j,k}(X_{i_{u +
1}})\cdots \widetilde\delta_{j,k}( X_{i_{u+v}}) \right) \right| \leq
(u+v+uv)\,\|\widetilde\delta_{j,k}\|_\infty^{u+v-2}(\Lip \widetilde\delta_{j,k})^2 \lambda (r).
$$
The right hand side term is bounded by $(2^{j/2}2\|\delta\|_\infty)^{u+v-2}\Lip \delta^2 2^{3j}\lambda (r)$. At this point, we do not achieve Assumption \eqref{cpr} as this bound depends on $j$ via an extra term $2^{3j}$. An additional assumption is needed to ensure Assumption~{\bf (D)}:
\begin{description}
\item[(J)] The joint densities $f_{r}$ of $(X_0,X_r)$ exist and are
bounded for $r>0$.
\end{description}
\begin{prop}\label{etaD}
Assume that $(X_t)_{t\in\Z}$ is a $\lambda$-weakly dependent process satisfying {\bf (J)} and such that there exist $a',b',c'>0$ and $a'',b'',c''>0$ with:
$$
\lambda(r)\le c' \exp(-a'r^{b'})\mbox{ and }\|f_r\|_\infty \le c'' \exp(a''r^{b''})\mbox{ for all }r>0.
$$
If $b''<b'$ then {\bf (D)} holds for some $C_0>0$ with $M_2=2\|\delta\|_\infty$, $a=a'/4$ and $b=b'$.
\end{prop}
The proof of this Proposition is given in Subsection \ref{proofsimus}.

\subsection{\texorpdfstring{Examples of $\lambda$-weakly dependent processes}{Examples of lambda-weakly dependent processes}}\label{exeta}
Firstly we give in this Section two generic $\lambda$-weakly dependent models, Bernoulli shifts and random processes with infinite connections. Secondly we detail the conditions of Proposition \ref{etaD} in three more specific models.

Let $H:\R^\Z\to [0,1]$ be a measurable function. A Bernoulli shift with innovations $\xi_t$ is defined as
$$
X_t=H\left((\xi_{t-i})_{i\in \Z}\right), \qquad t\in \Z.
$$
According to \cite{Doukhan2007}, such Bernoulli shifts are $\lambda-$weakly dependent with $\lambda(r)\le 2v([r/2])$. If $(\xi_t)_{t\in \Z}$ is iid, $(v(r))_{r>0}$ is a non-increasing sequence satisfying
$$\E\left|H\left(\xi_{j},j\in\Z\right)-H\left(\xi_{j}',j\in\Z\right)\right|\le v(r)\mbox{ for all }r>0,$$
where the iid sequence $(\xi_{j}')_{j\in\Z}$ is such that $\xi_{j}'=\xi_j$ for $|j|\le r$ and $\xi_{j}'$ independent of $\xi_{j}$ otherwise. If the weak dependence coefficients $\lambda_{\xi}(r)$ refer to $(\xi_t)_{t\in \Z}$, we can compute these of $(X_t)_{t\in \Z}$. More precisely, the result in \cite{Doukhan2007} states that if there exists $\ell\ge 0$ such that
$$ |H(x)-H(y)|\le b_s(\|z\|^\ell\vee 1)|x_s-y_s|,$$
for some sequence $b_j\ge0$ satisfying $\sum_j|j|b_j<\infty$ and where $x,y\in\R^{\Z}$ coincide except for the index $s\in\Z$ and $\|x\|=\sup_{i\in\Z}|x_i|$, if $\E|\xi_0|^{m'}<\infty$ for some $m'\ge \ell+1$
then $(X_t)_{t\in \Z}$ is $\lambda$-weakly dependent with
$$
\lambda(k)\le c\inf_{r\le [k/2]} \left[\sum_{|j|\ge r}|j|b_j+ (2r+1)^{2 }\lambda_{\xi}(k-2r)^{\frac{{m'}-1-\ell}{{m'}-1+\ell}}\right],\mbox{ for some }c>0.
$$
Different values of $b$ in Assumption \eqref{condrho} may arise naturally
when realizing the minimum of the equation above, see below for some classical examples.

Another approach is the one of random processes with infinite connections considered in \cite{Doukhanc}. Let $F:[0,1]^{\Z/\{0\}}\times \R\to [0,1]$ be measurable. Under suitable conditions on $F$, the stationary solution of the equation
$$
X_t=F((X_j,j\neq t),\xi_t),~~a.s.,
$$
exists and is $\lambda$-weakly dependent. We refer the reader to \cite{Doukhanc} for more details and we end the Section with some specific $\lambda$-weakly dependent models.
\subsubsection{Infinite moving average}
A Bernoulli shift is an
infinite moving average process if
\begin{equation}
\label{linear} X_t=\sum_{i\in\Z}\alpha_i\xi_{t-i}.
\end{equation}
If $\xi_t$ are iid random variables satisfying $\E|\xi_0|\le 1$ then $(X_t)_{t\in\Z}$ is
$\lambda$-weakly dependent with $\lambda(r)\le 4\sum_{|j|>[r/2]}|a_j|$.
If $a_j\le K\alpha^{|j|}$ for $j\neq0$,
$K>0$ and $0<\alpha<1$ then the condition on $\lambda$ in Proposition \ref{etaD} holds with
$b'=1$. Then Assumption \eqref{condrho} is ensured with $b=1$.
\subsubsection{{\em LARCH($\infty$)} model} Let $( \xi_t )_{t \in \Z}$ be an iid
centered real valued sequence and $a, a_j , {j \in \N^{\ast}}$ be
real numbers. LARCH($\infty$) models are solutions of the recurrence
equation
\begin{equation}\label{archv1}
X_t = \xi_t \left( a + {\sum_{j\neq 0}} a_j X_{t -j} \right).
\end{equation}
The stationary solution of (\ref{archv1}) 
satisfies the condition on $\lambda$ in Proposition \ref{etaD} with
$b'=1/2$ if there exists
$K,\alpha>0$ and $\alpha<1$ such that $a_j\le K\alpha^{|j|}$ for all
$j\neq0$. Then Assumption \eqref{condrho} is ensured with $b=1/2$. See \cite{DoukhanTeyssiere} for applications of this model in econometrics.

\subsubsection{Affine model}
Let us consider the stationary solution $(X_t)_{t\in\Z}$ of the equation $$X_t=M(X_{t-1},X_{t-2},\ldots) \xi_t+f(X_{t-1},X_{t-2},\ldots),$$ where $M$ and $f$ are both Lipschitz functions. 
This model contains various time series processes such as ARCH, GARCH, ARMA, ARMA-GARCH, etc. If the $\xi_t$ are iid random variables with a bounded marginal density, then {\bf (J)} holds and the joint densities are uniformly bounded, as stated in the Appendix of \cite{Doukhand}. Moreover if the functions $M$ and $f$ have exponentially decreasing Lipschitz coefficients, then conditions of Proposition \ref{etaD} hold with $b'=1/2$, $b''=0$ and Assumption \eqref{condrho} follows with $b=1/2$.

\section{Cross-validation procedures and simulations}\label{simus}

The aim of this Section is to evaluate the applicability and the quality of the procedure on simulated data. Even if the estimators are adaptive with respect to the regularity of the density function, a constant appears in the threshold levels that we cannot calibrate if the dependence properties of the observations are unknown. Then we develop a cross-validation scheme in order to apply concretely the estimator on simulated data. We investigate several examples of dependence for the simulated observations that satisfy the convergence result of this paper. We also give  counter-examples where the estimators fail to converge near-minimaxly.

\subsection{Cross-validation procedures}


According to Theorem \ref{th}, let us fix $j_0$ as the smallest integer larger than $\log(n)(1+N)^{-1}$ and define $\widehat f_n^{HTCV}$ and $\widehat f_n^{STCV}$ respectively as the hard and soft-thresholding estimators associated with threshold levels $((\widehat\lambda_j)_{j_0\le j\le j^\ast})$. Here $j^\ast=\log_2 n$ and $((\widehat\lambda_j)_{j_0\le j\le j^\ast})$ are determined by cross-validation procedures: 
$\widehat\lambda_j=\mathop{\mbox{Argmin}}_{\lambda}CV_j(\lambda)$
where cross validation criterion are respectively defined by the equations
\begin{eqnarray*}
\mbox{HTCV: }&&CV_j(\lambda)=\sum_{k\in\Z}\1_{\{|\widehat\beta_{j,k}|\ge \lambda\}}\left[\widehat\beta_{j,k}^2-\frac{2}{n(n-1)}\sum_{1\le i\neq h\le n}\psi_{j,k}(X_i)\psi_{j,k}(X_h)\right],\\
\mbox{STCV: }&&CV_j(\lambda)=\sum_{k\in\Z}\1_{\{|\widehat\beta_{j,k}|\ge \lambda\}}\left[\widehat\beta_{j,k}^2-\frac{2}{n(n-1)}\sum_{1\le i\neq h\le n}\psi_{j,k}(X_i)\psi_{j,k}(X_h)+\lambda^2\right].
\end{eqnarray*}
These criterion are obtained by approximating the coefficients $\beta_{j,k}$ by $\hat\beta_{j,k}$ and the products
$\beta_{j,k}\beta_{j,k'}$ by $\sum_{h\neq i}
\psi_{j,k}(X_i)\psi_{j,k}(X_h)/(n(n-1))$ in the Integrated Square Error.\\

The estimator $\widehat j_1$ is defined as the smallest integer such that $CV_j (\widehat \lambda_j)=0$ for all $\widehat j_1\le j\le j^\ast$. 
Notice that cross-validation procedures may consider larger resolution levels than the estimators $\widehat f_n$ as $j^\ast$ is larger than $j_1$ given in Theorem \ref{th}. 

\subsection{Different dependent samplings satisfying {\bf (D)}}

To illustrate the behavior of this cross-validation scheme, we simulate three different weak-dependence cases with the same marginal absolutely continuous distribution $F$. The simulations were carried out as follows:
\begin{description}
\item[Case 1] Independent observations are given by $X_i=F^{-1}(U_i)$ where the $U_i$ are simulations of independent variables, uniform on [0,1].
\vspace{16pt}
\item[Case 2] A $\widetilde \phi$-weakly dependent process is obtained by the equation $X_i=F^{-1}(G(Y_i))$ for $i=1,\dots,n$ with $G(x)=2\sqrt{x(1-x)}/\pi$ and $(Y_i)_{i=1,\dots,n}$ given by \begin{center}$Y_1=G^{-1}(U_1)$\\ and, recursively, $Y_i=T^{i-1}(Y_1)$ for $2\le i\le n$ with $T(x)=4x(1-x)$.\end{center} Note that for all $1\le i\le n$ the $Y_i$ admits the repartition function $G$ the invariant distribution of $T$. Moreover the sequence $(Y_1,\ldots,Y_n)$ satisfies the assumptions of Proposition \ref{tildephi}, see~\cite{Prieur} for details. 
\vspace{16pt}

\item[Case 3] A $\lambda$-weakly dependent process resulting from the transform $X_i=F^{-1}(G(Y_i))$ of variables $(Y_t)_{t\in\Z}$ which are solution of 
$$Y_t=2(Y_{t-1}+Y_{t+1})/5 + 5\xi_t/21
$$ where $(\xi_t)_{t\in\Z}$ is an iid sequence of Bernoulli variables with parameter $1/2$.
The stationary solution of this equation admits the representation
$$Y_t=\sum_{j\in\Z}a_j\xi_{t-j},$$
where $a_j=1/3(1/2)^{|j|}$. This solution belongs to $[0, 1]$ and its marginal distribution $G$ is the one of $(U+U'+\xi_0)/3$ where $U$ and $U'$ are independent variables following $\mathcal{U}([0,1])$. For $1\le i\le n$, the solution $Y_i$ is
approximated by $Y_i^{(N)}$  for $1\le j\le N$ and $j-N\le i\le n+N-j$ where $Y_i^{(j)}$ is generated according to the convergent algorithm given in \cite{Doukhanc}: the initial values $Y_i^{(0)}$ are fixed equal to $0$ and, given simulated variables $(\xi_i)_{-N\le i\le n+N}$, we define recursively $Y_i^{(j)}=2(Y_{i-1}^{(j-1)}+Y_{i+1}^{(j-1)})/5 + 5\xi_i/21$ for $1\le i\le N$ and $i-N\le t\le n+N-i$. Error of approximation are negligible as decreasing exponentially fast with the parameter $N$ that we fix $N=n$, see Lemma 6 of \cite{Doukhanc} for more details. Moreover, Proposition 2.1 of \cite{Dedecker2007} ensures that there exists $a,C>0$ such that $\lambda(r)\le C\exp(-ar)$ for the process $(Y_t)_{t\in\Z}$ and consequently for $(X_t)_{t\in\Z}$.
\end{description}
Two different density functions are considered. The first one is a mixture of a sinus function and a uniform distribution, presenting a discontinuitie, and the second one is a mixture of two gaussian distributions.

\subsection{\texorpdfstring{Comparison of $\widehat f_n^{HTCV}$ and $\widehat f_n^{STCV}$}{Comparison of hard and soft thresholding procedures}}
The first density considered is a mixture between a sinus function and a uniform distribution.

The calculations were carried out on MATLAB on a Unix environment. We considered $n=2^{10}$ observations repeated $M=500$ times for each of the three weak dependence cases. Once the data simulated we applied cross-validation procedures. The usual DWT algorithm proposed by \cite{Mallat} and implemented in the Wavelab toolbox by Donoho and his collaborators (available on \cite{wavelab}) only gives values of wavelet functions on an equidistant grid. As one needs to compute these values at given data points, we consider an equidistant grid of $I$ points with the number of points $I$ huge with respect to the number of observations $n$. Then we approximate the values $\psi_{j,k}(X_i)$ by $\psi_{j,k}([X_i I]/I)$ where $[x]$ denotes the closest integer from any real number $x$. The wavelets used for the decomposition are Daubechies Symmlets with $N=8$ zero-moments. Notice that another possible scheme is the algorithm of \cite{Vidakovic} that gives directly the values $\psi_{j,k}(X_i)$. But as it needs much more calculus time it has not been used here for convenience.

\begin{figure}[!ht]
\begin{center}
\hspace{-0.5cm}\includegraphics[height=3cm,width=5cm]{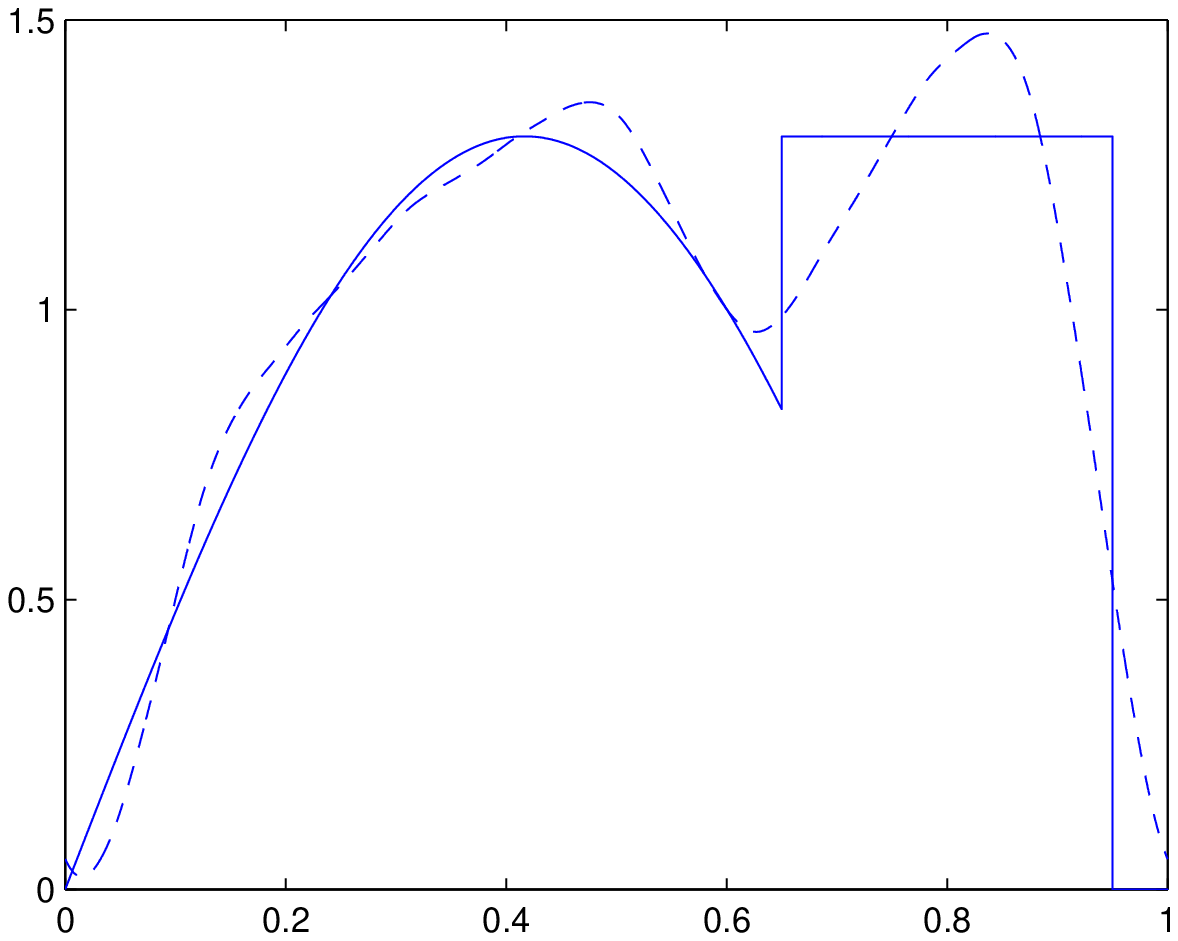}
\hspace{-0.5cm}\includegraphics[height=3cm,width=5cm]{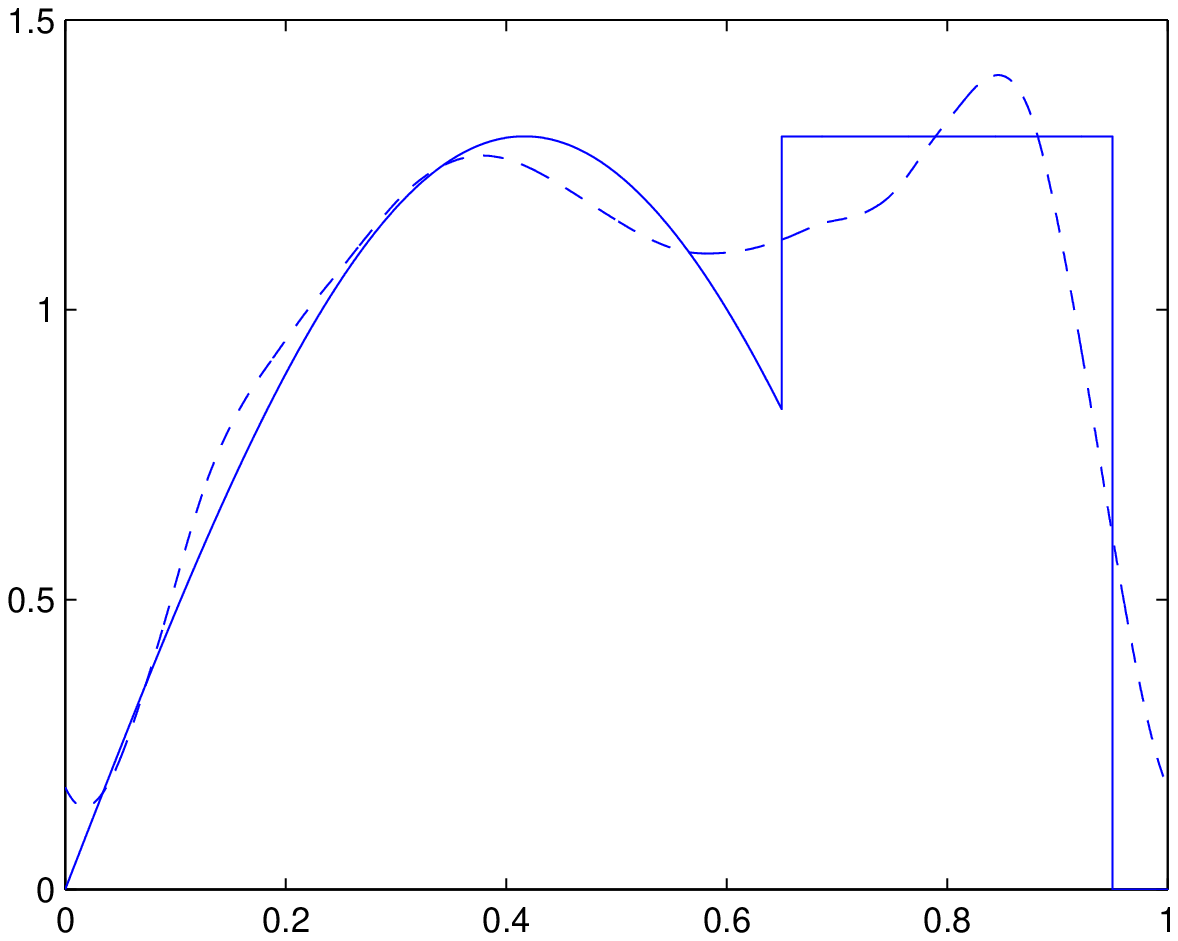}
\hspace{-0.5cm}\includegraphics[height=3cm,width=5cm]{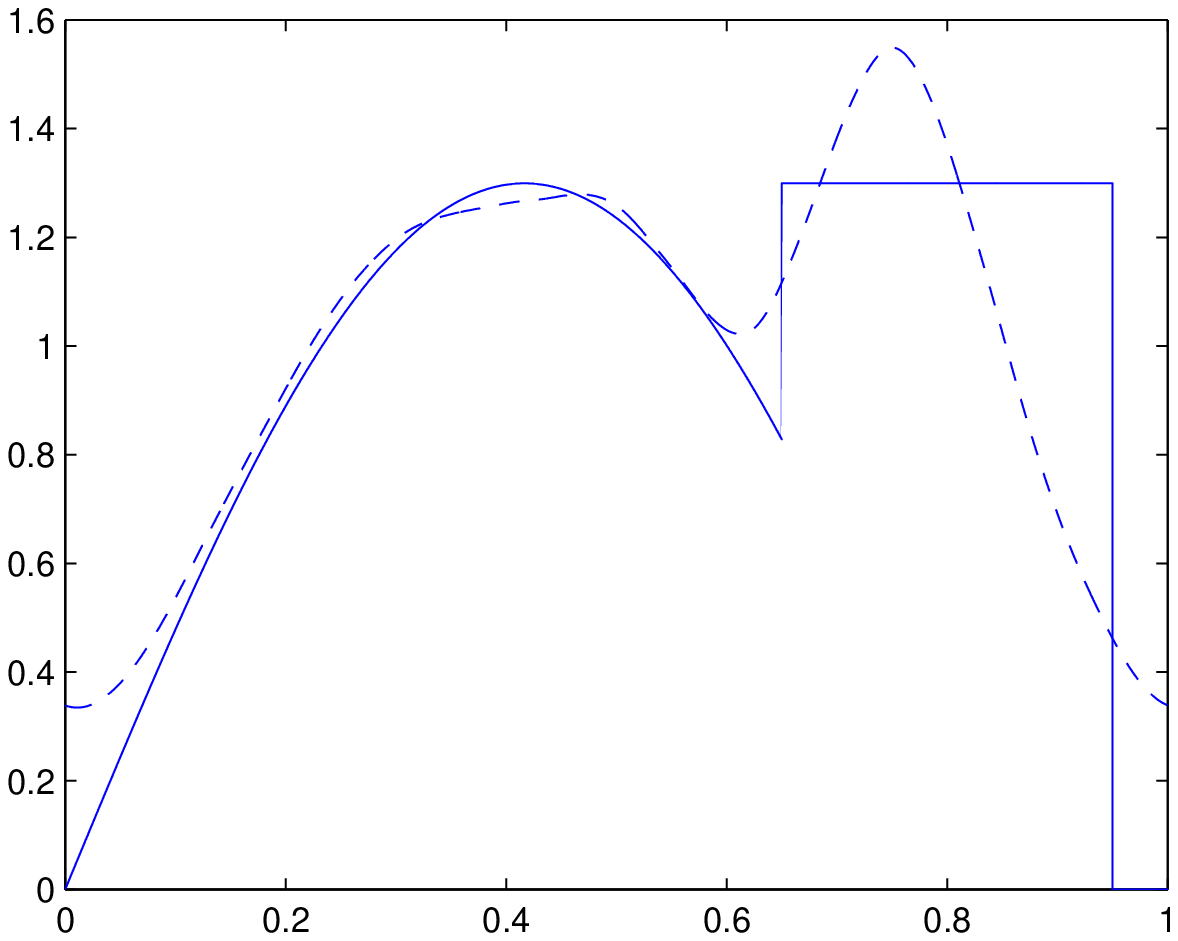}
\end{center}
\vspace{-0.5cm}\caption{{\small Examples of estimators $\widehat f_n^{HTCV}$ obtained on $2^{10}$ observations. The true distribution is represented in dashed lines. Figures from left to right correspond respectively to {\bf Case~1}, {\bf Case~2} and {\bf Case~3}.}}\label{fig:densiteH}
\end{figure}

\begin{figure}[!ht]
\begin{center}
\hspace{-0.5cm}\includegraphics[height=3cm,width=5cm]{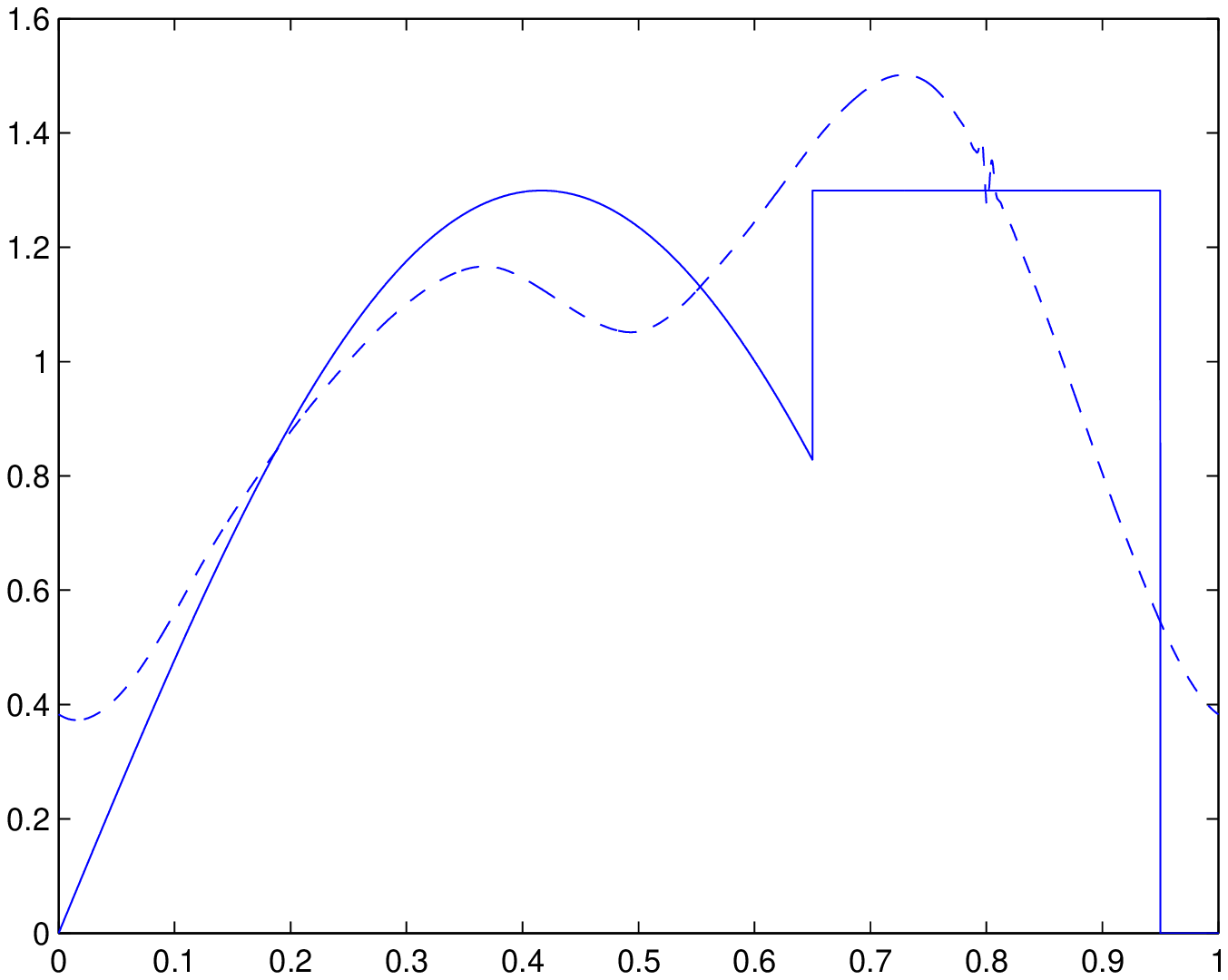}
\hspace{-0.5cm}\includegraphics[height=3cm,width=5cm]{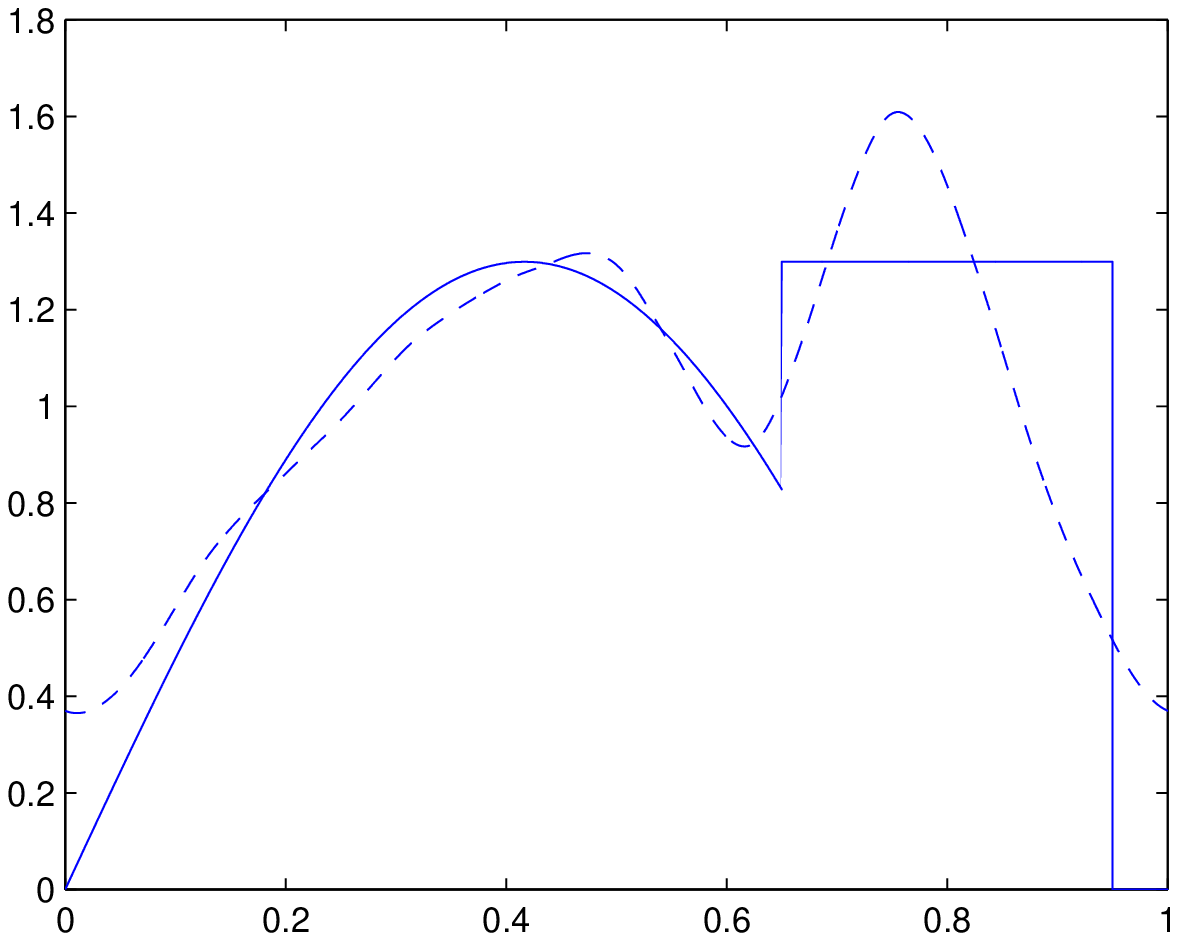}
\hspace{-0.5cm}\includegraphics[height=3cm,width=5cm]{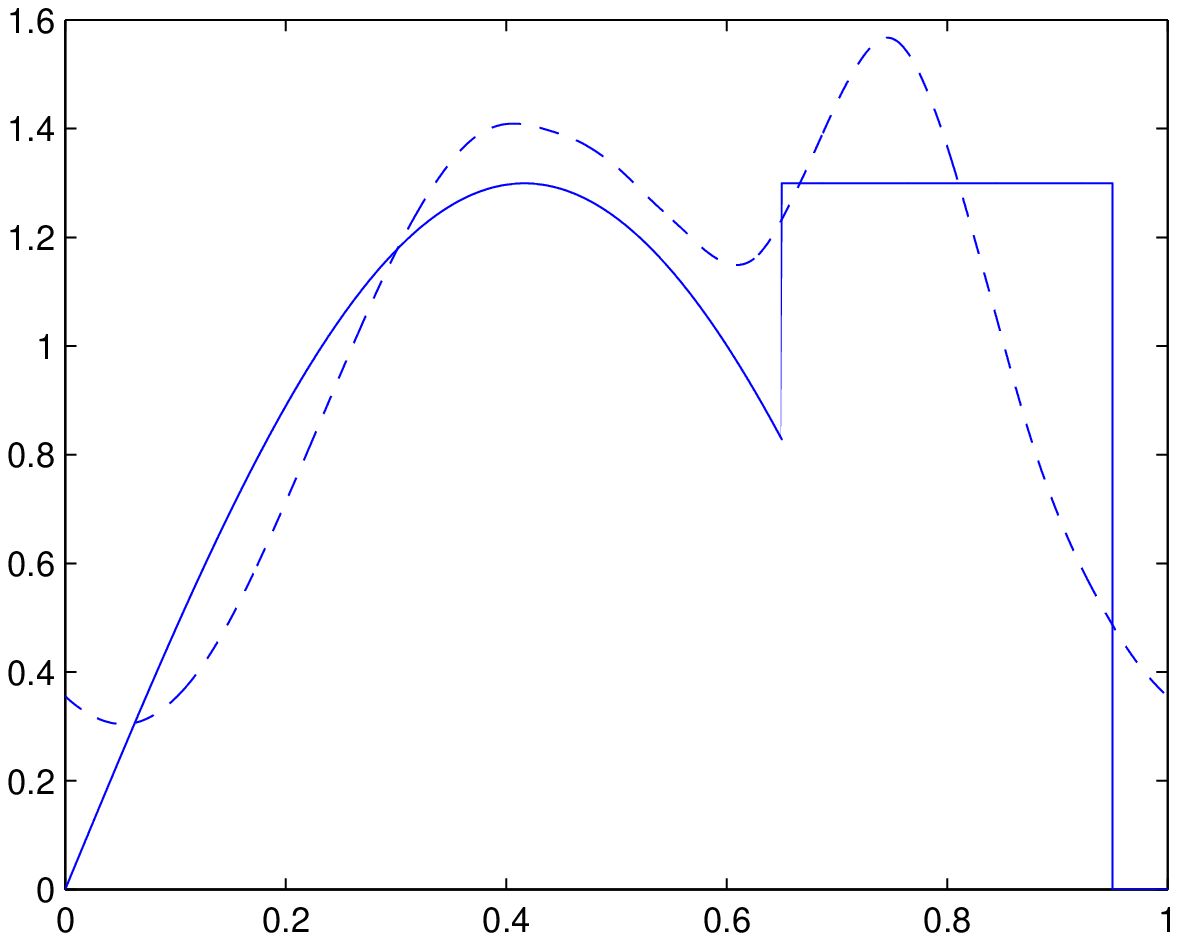}
\end{center}
\vspace{-0.5cm}\caption{{\small Examples of estimators $\widehat f_n^{STCV}$ obtained on $2^{10}$ observations. The true distribution is represented in dashed lines. Figures from left to right correspond respectively to {\bf Case~1}, {\bf Case~2} and {\bf Case~3}.}}\label{fig:densiteS}
\end{figure}

In Figures~\ref{fig:densiteH} and~\ref{fig:densiteS} are represented the estimators $\widehat f_n^{HTCV}$ and $\widehat f_n^{STCV}$ and the true density function $f$ in different weakly dependent cases. The quality of the estimators is  visually good. According to Figures~\ref{fig:densiteH} and~\ref{fig:densiteS} the weak dependence properties of the simulated data do not seem to affect both procedures of estimation. Density estimators presented in Figures~\ref{fig:densiteH} and~\ref{fig:densiteS} do not detect the discontinuity in the density. Actually, for any finite data set, not enough simulated values are concentrated around the discontinuity to allow estimators to detect it. 

Table~\ref{tab:mise} gives approximations by Monte-Carlo method of the MISE. MISE values of the estimators have the same order whereas the weak dependent cases and $\widehat f_n^{STCV}$ is preferable in all the cases. 

\begin{table}[!ht]
\centering
\begin{tabular}{@{ } l@{\hspace{20pt}} c @{\hspace{20pt}} c @{\hspace{20pt}}r @{ }}
\toprule
\multicolumn{4}{c}{MISE of the estimation} \\
\cmidrule(r){1-4} 
& Case 1 & Case 2 & Case 3\\
\midrule
HTCV & 0.096696 & 0.077064 & 0.097193 \\
STCV & 
0.082934 & 0.06586 & 0.097184 \\
\bottomrule
\end{tabular}
\vspace{0.5cm}
\caption{{\small MISE approximated by MC on 500 simulations of samples of size $n=2^{10}$.}}\label{tab:mise}
\end{table}

\begin{figure}[!ht]
\begin{center}
\hspace{-0.5cm}\includegraphics[height=4cm,width=8cm]{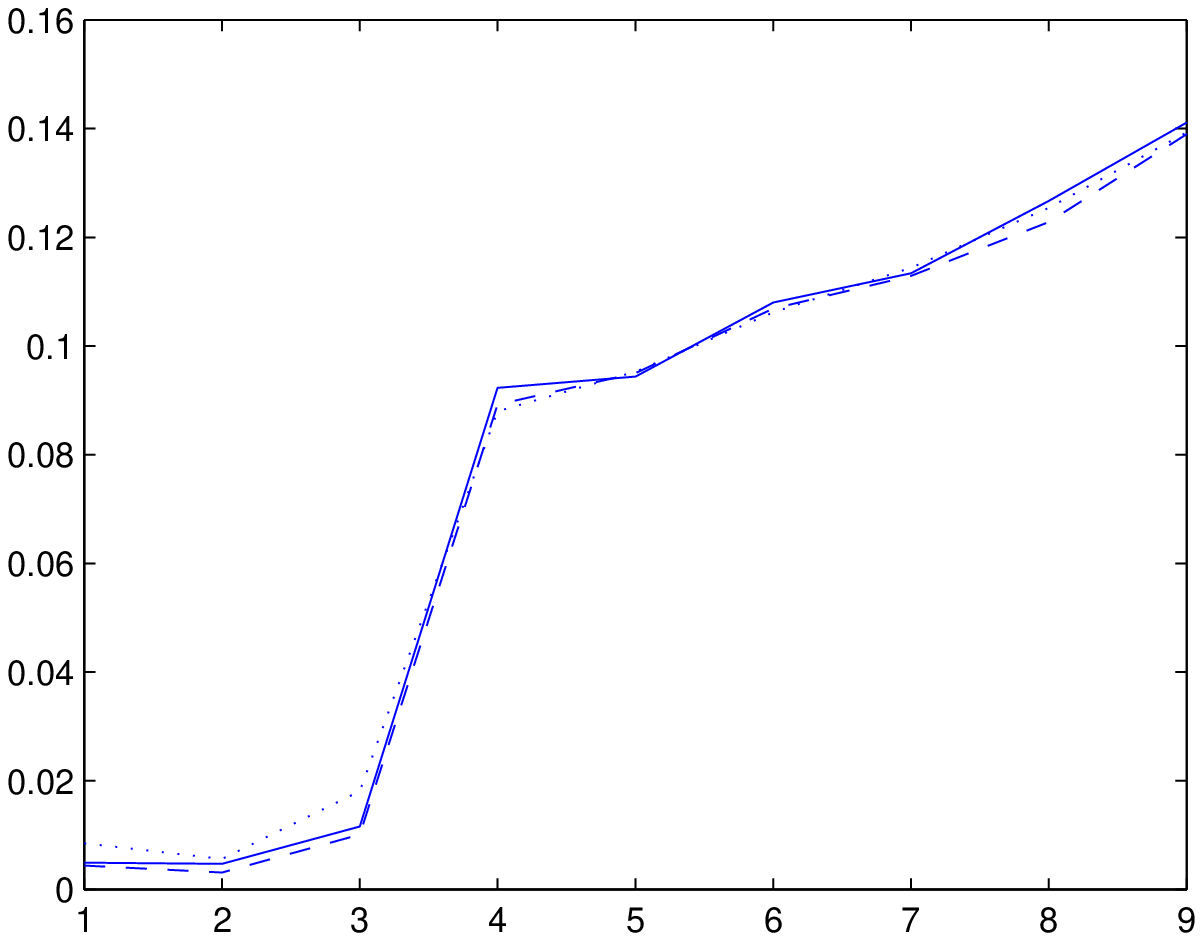}
\includegraphics[height=4cm,width=8cm]{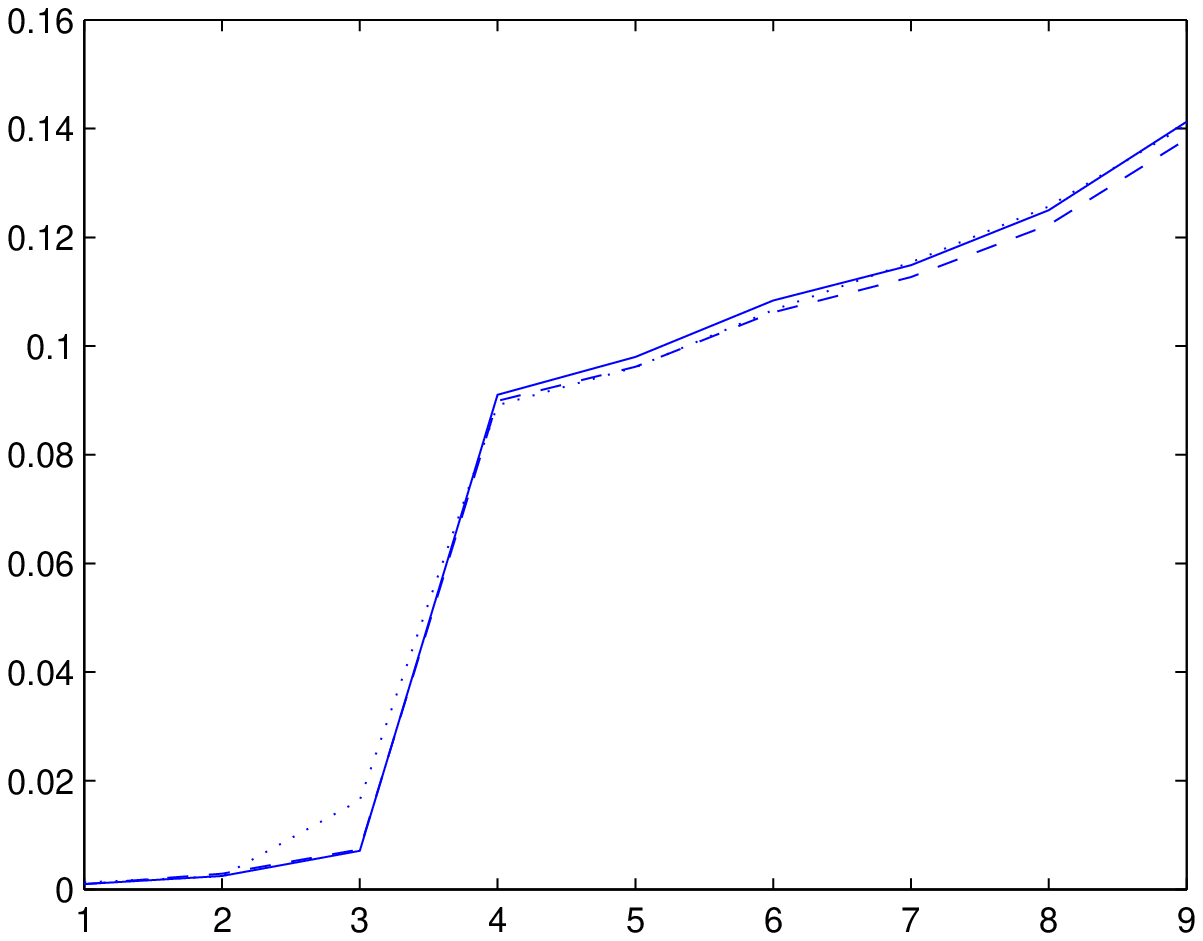}\\
\end{center}
\vspace{-0.5cm}\caption{{\small Means of the proportions of the threshold levels obtained by cross-validation with respect to the resolution levels for hard-thresholding (left) and soft-thresholding (right). 
{\bf Case 1} corresponds to the solid line, {\bf Case 2} to the dashed line and {\bf Case 3} to the dotted line.}}\label{fig:seuils}
\end{figure}

In Figure \ref{fig:seuils} threshold levels are represented with respect to resolution levels. Their behaviors are similar in all cases: the threshold levels increase with respect to resolution levels. For small resolutions, both HTCV and STCV procedures are close as $\widehat\lambda_j^2$ is negligible in $CV_j(\lambda)$. For high resolution it is also the case as $\widehat\lambda_j^2$ is big enough to kill almost all $\widehat\beta_{j,k}$. Moreover these figures tend to confirm that threshold levels do not depend on weak dependence. Finally remark that the curves do not behave in square root of $j$ as theoretical threshold levels given in Theorem~\ref{th}.

\vspace{0.5cm}

\begin{figure}[!ht]
\begin{center}
\hspace{-0.5cm}\includegraphics[height=4cm,width=8cm]{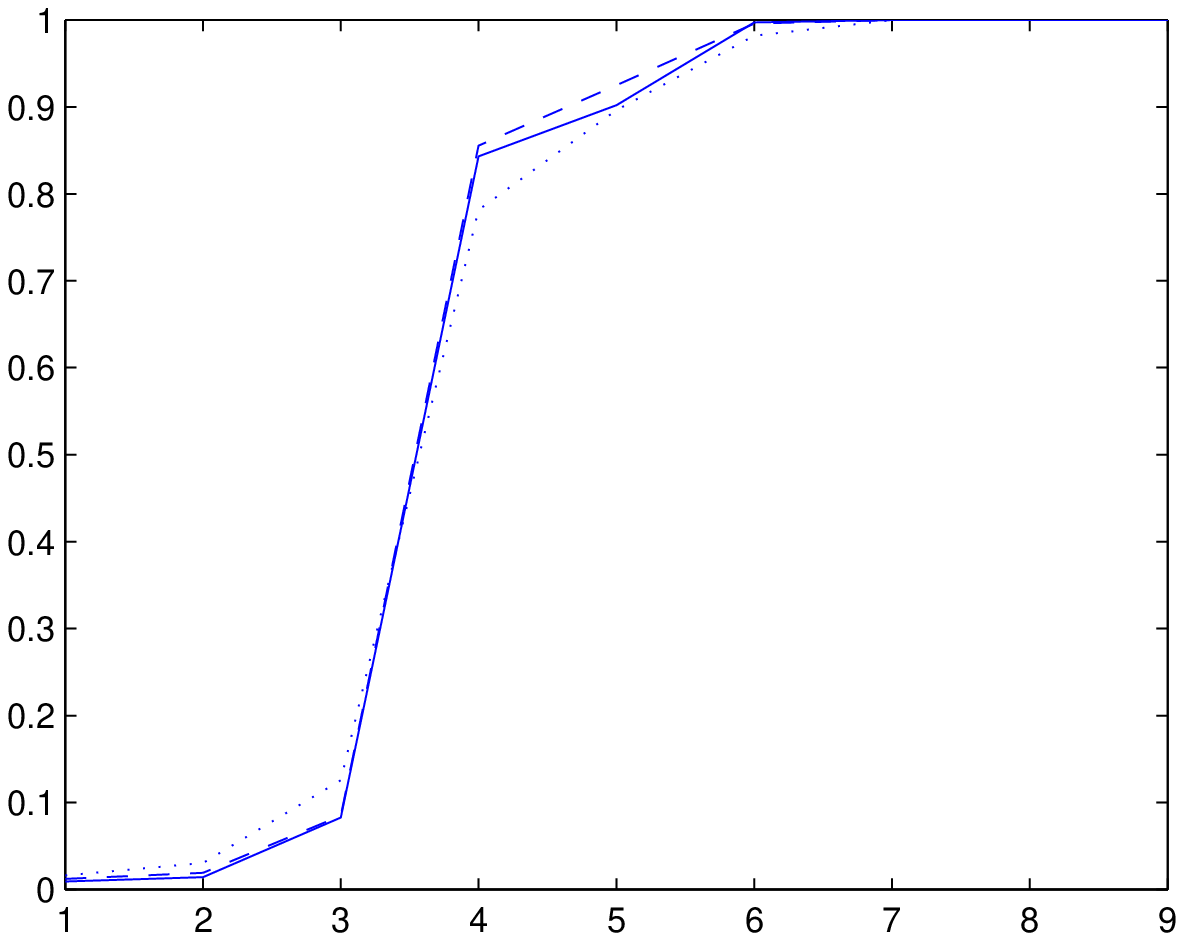}
\includegraphics[height=4cm,width=8cm]{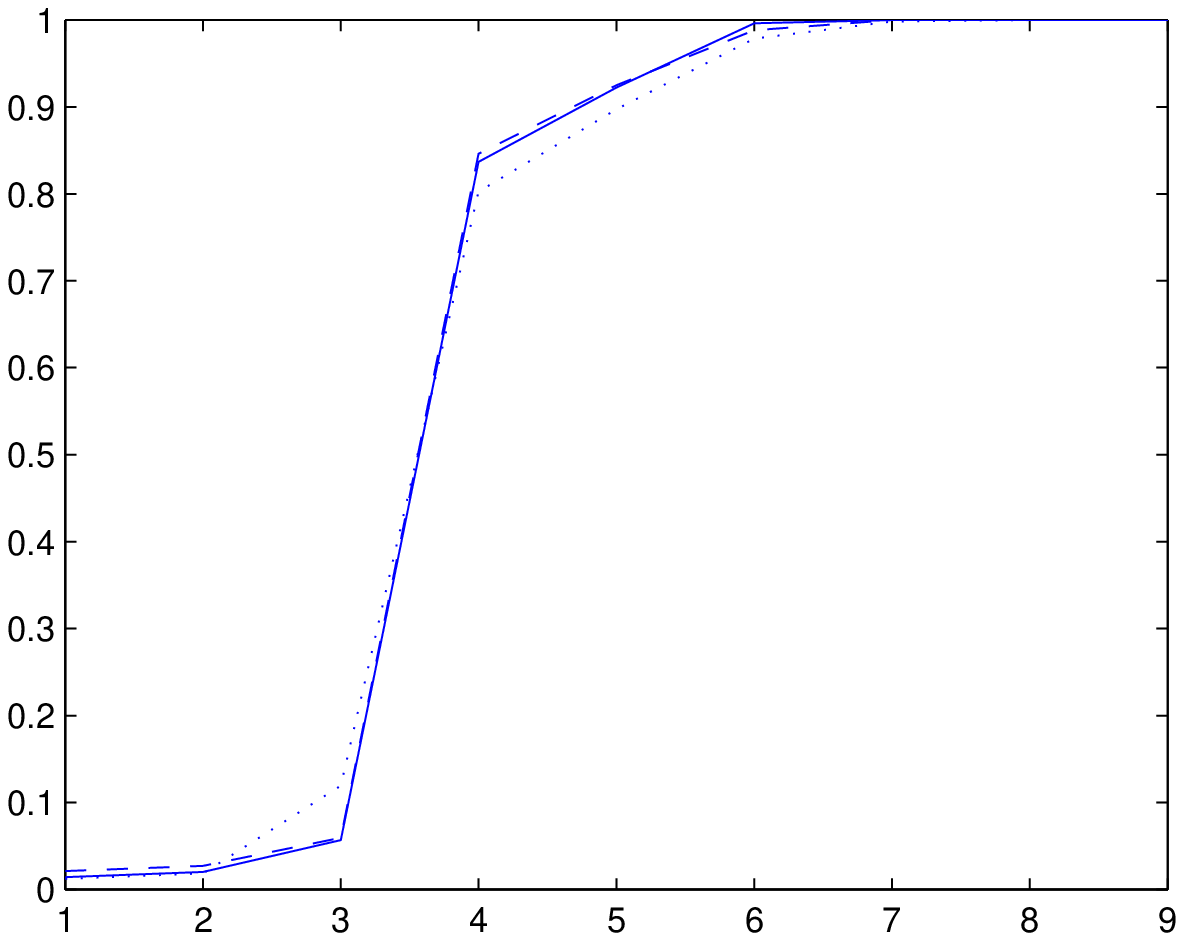}\\
\end{center}
\vspace{-0.5cm}\caption{{\small Means of the proportions of thresholded coefficients with respect to the resolution levels for hard-thresholding (left) and soft-thresholding (right). 
{\bf Case 1} corresponds to the solid line, {\bf Case 2} to the dashed line and {\bf Case 3} to the dotted line.}}\label{fig:nb_seuil}
\end{figure}

After looking at the threshold levels values, we give in Figure \ref{fig:nb_seuil} the frequencies of the $\widehat\beta_{j,k}$ that are less than $\widehat \lambda_j$ with respect to $j$. As these frequencies are not discretized in two values $0$ and $1$, we can infer that both $\widehat f_n^{HTCV}$ and $\widehat f_n^{STCV}$ are not equivalent to linear estimators. It is encouraging for both methods as the results in \cite{DonJohnKerkPic} state that linear estimators are not near-minimax. One can also see in these figures that frequencies of effective thresholds are the same among the weak dependence cases. 


\vspace{0.5cm}

Finally, means of higher resolution levels are given in last Table \ref{tab:hrl}. According to Theorem \ref{th}, the values of this parameter do depend on the cases of weak dependence. But no significant differences appear on simulations. 

\begin{table}[!ht]
\centering
\begin{tabular}{@{ } l@{\hspace{20pt}} c @{\hspace{20pt}} c @{\hspace{20pt}}r @{ }}
\toprule
\multicolumn{4}{c}{Mean of $\widehat j_1$}\\
\cmidrule(r){1-4} 
& Case 1 & Case 2 & Case 3\\
HTCV & 5.168 & 5.14 & 5.13 \\
\midrule
STCV & 
5.14 & 5.04 & 5.13 \\
\bottomrule
\end{tabular}
\vspace{0.5cm}
\caption{{\small Means of $\widehat j_1$ on 500 simulations of $n=2^{10}$ observations.}}\label{tab:hrl}
\end{table}

\subsection{Comparison with kernel estimators}
As results
computed in the last Subsection are systematically better for the STCV estimators than for the HTCV ones, we only present in the
sequel results for the STCV estimators.  Let us now consider the case of a density function that is a mixture of normal distributions. Like in \cite{HardleKerkPicTsyb}, we compare the quality of the wavelet estimator $\hat f^{STCV}$ with linear kernel estimators. The kernel used is Epanechnikov's one and we computed two choices for its width parameter: firstly we consider the width given by the rule of thumb of Matlab, (more precisely, the width is equal to $(q_3-q_1)/(2*0.6745)*(4/(3*n))^{1/5}$, where $q_1$ and $q_3$ denote respectively the first and the third quartile of the empirical distribution) and secondly we consider the width obtained by cross validation on the mean integrated squared error risk.

\begin{figure}[!ht]
\begin{center}
\hspace{-0.5cm}\includegraphics[height=3cm,width=5cm]{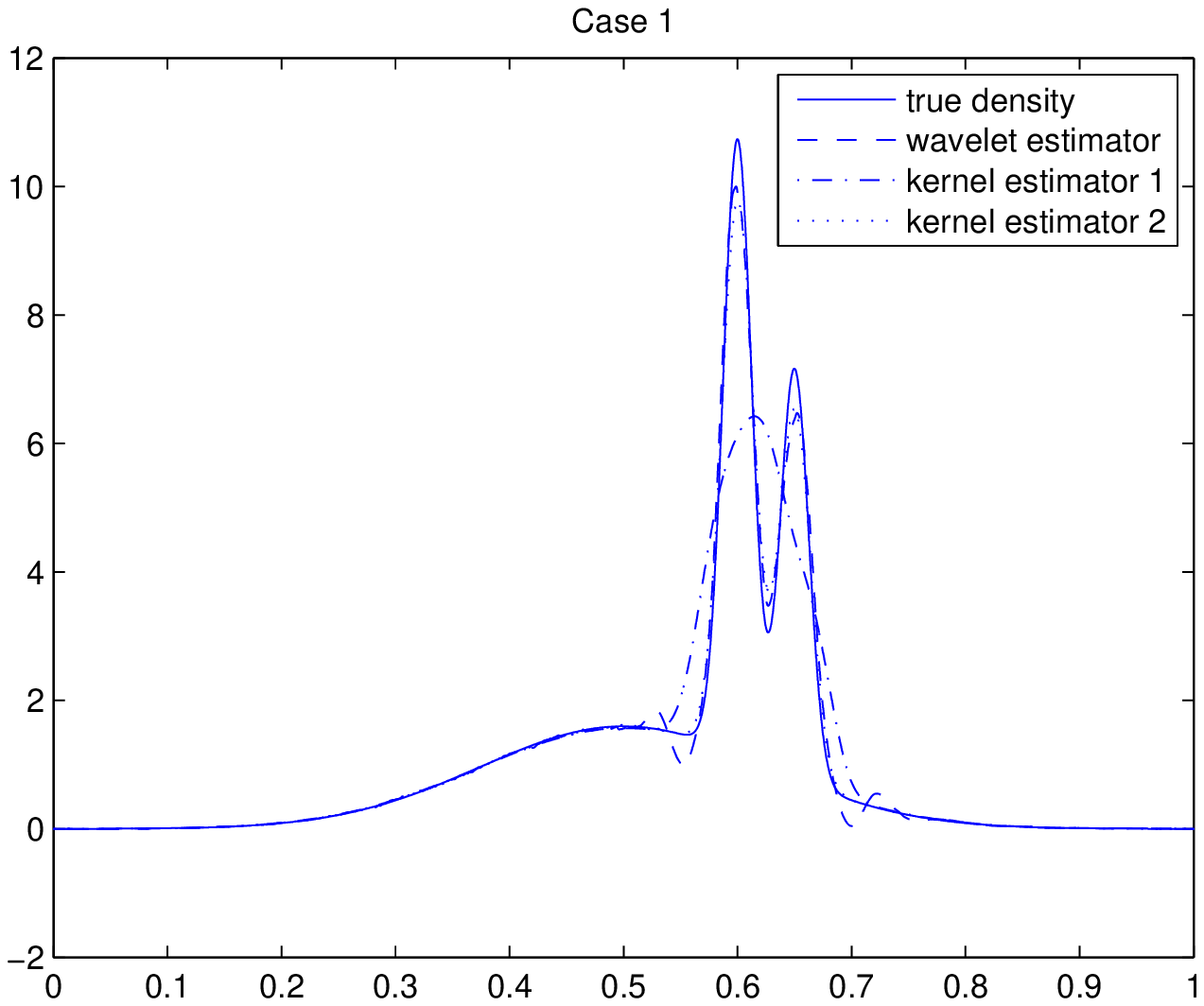} \hspace{-0.5cm}\includegraphics[height=3cm,width=5cm]{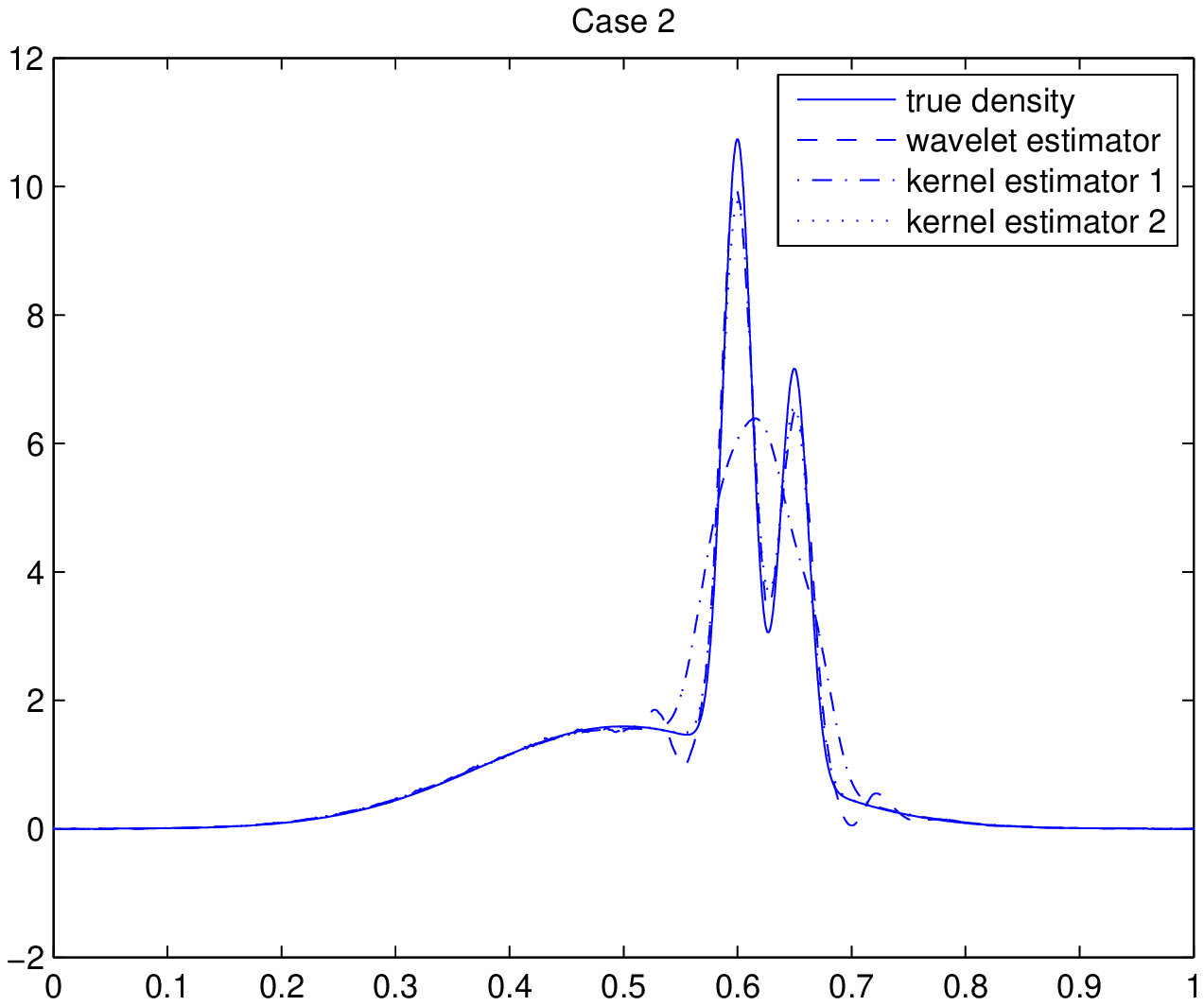} \hspace{-0.5cm}\includegraphics[height=3cm,width=5cm]{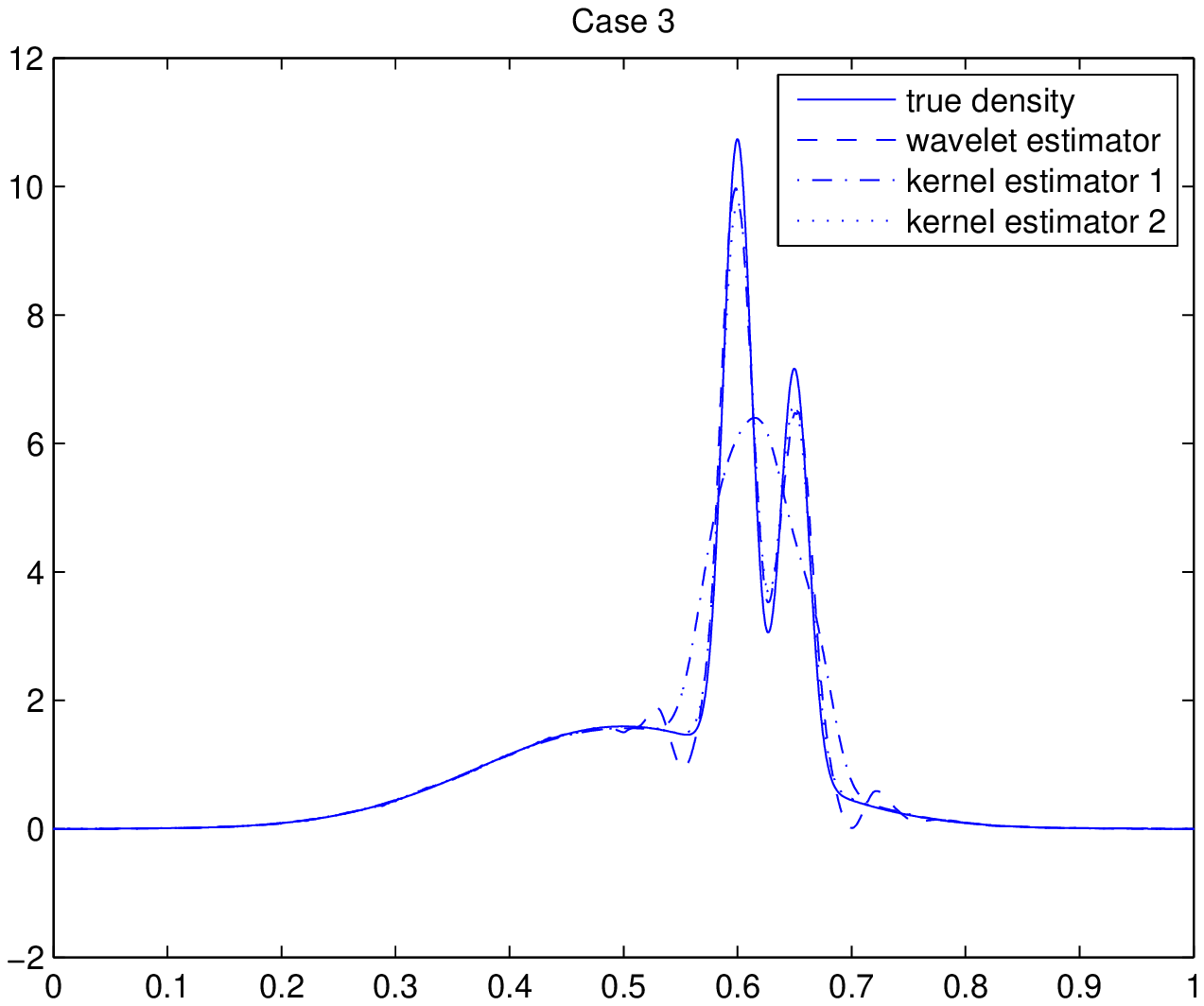}
\end{center}
\vspace{-0.5cm}\caption{{\small Means of estimators $\widehat f_n^{STCV}$ obtained from $2^{10}$ observations on $500$ simulations. The mean of the wavelet estimators is represented in dashed lines while the means of the kernel estimators are represented respectively in line with dots for the rule of thumb's width (kernel estimator 1) and in dots for the cross-validation width (kernel estimator 2). Figures from left to right correspond respectively to {\bf Case~1}, {\bf Case~2} and {\bf Case~3}.}}\label{fig:densiteS2}
\end{figure}

In Figure \ref{fig:densiteS2} are represented the means of wavelet and   kernel estimators in the different cases. Like in the last Subsection, there is no visual difference between the different cases of dependence. In this Figure, it appears that the mean of the kernel estimators fails to detect the two modes of the density when the width is choosing according to the rule of thumb. The bandwidth is overestimated in this case. The quality of the wavelet estimators $\widehat f_n^{STCV}$ and the kernel estimators with the width from a cross validation procedure are visually equivalent.

\begin{figure}[!ht]
\begin{center}
\hspace{-0.5cm}\includegraphics[height=3cm,width=5cm]{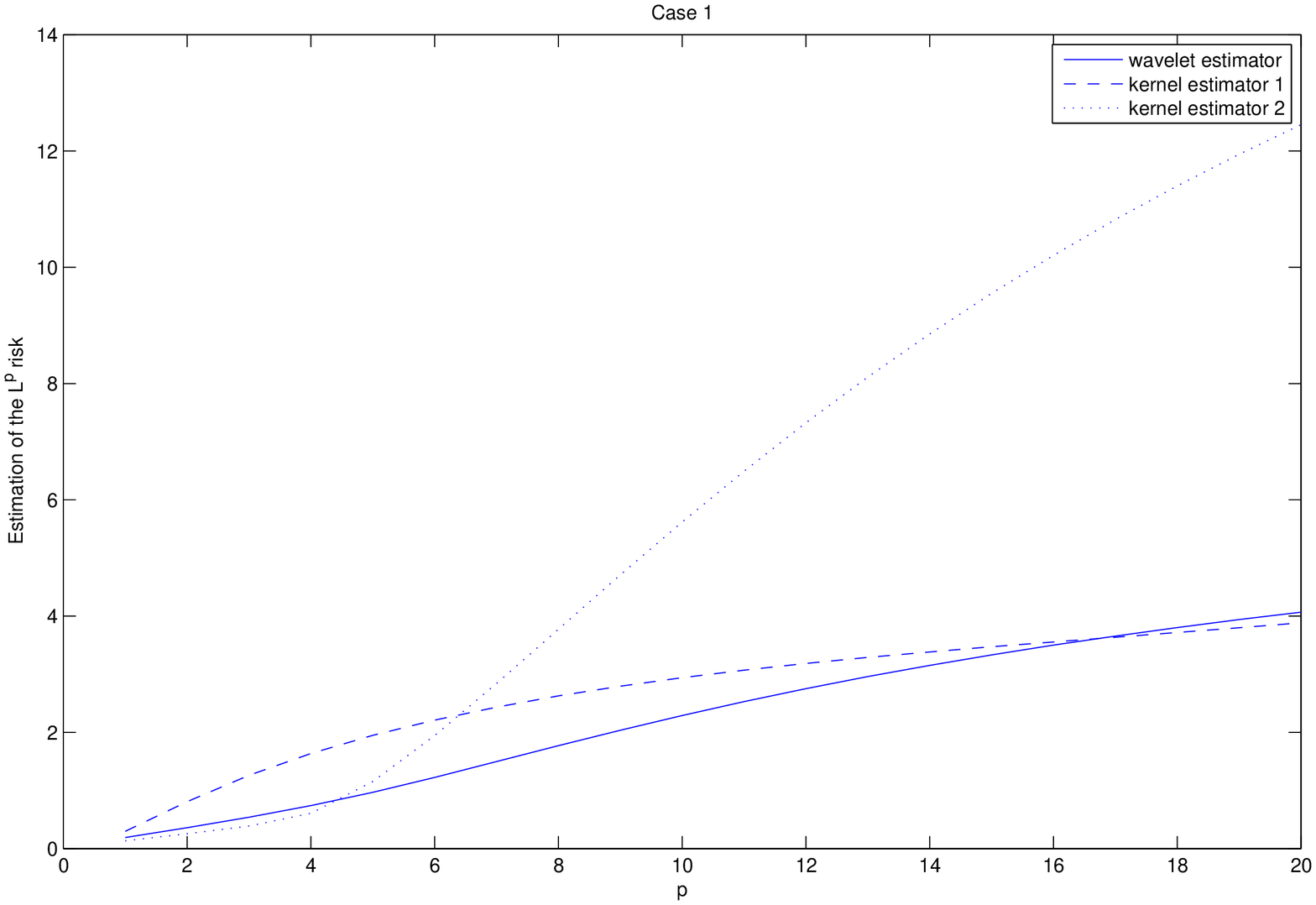} \hspace{-0.5cm}\includegraphics[height=3cm,width=5cm]{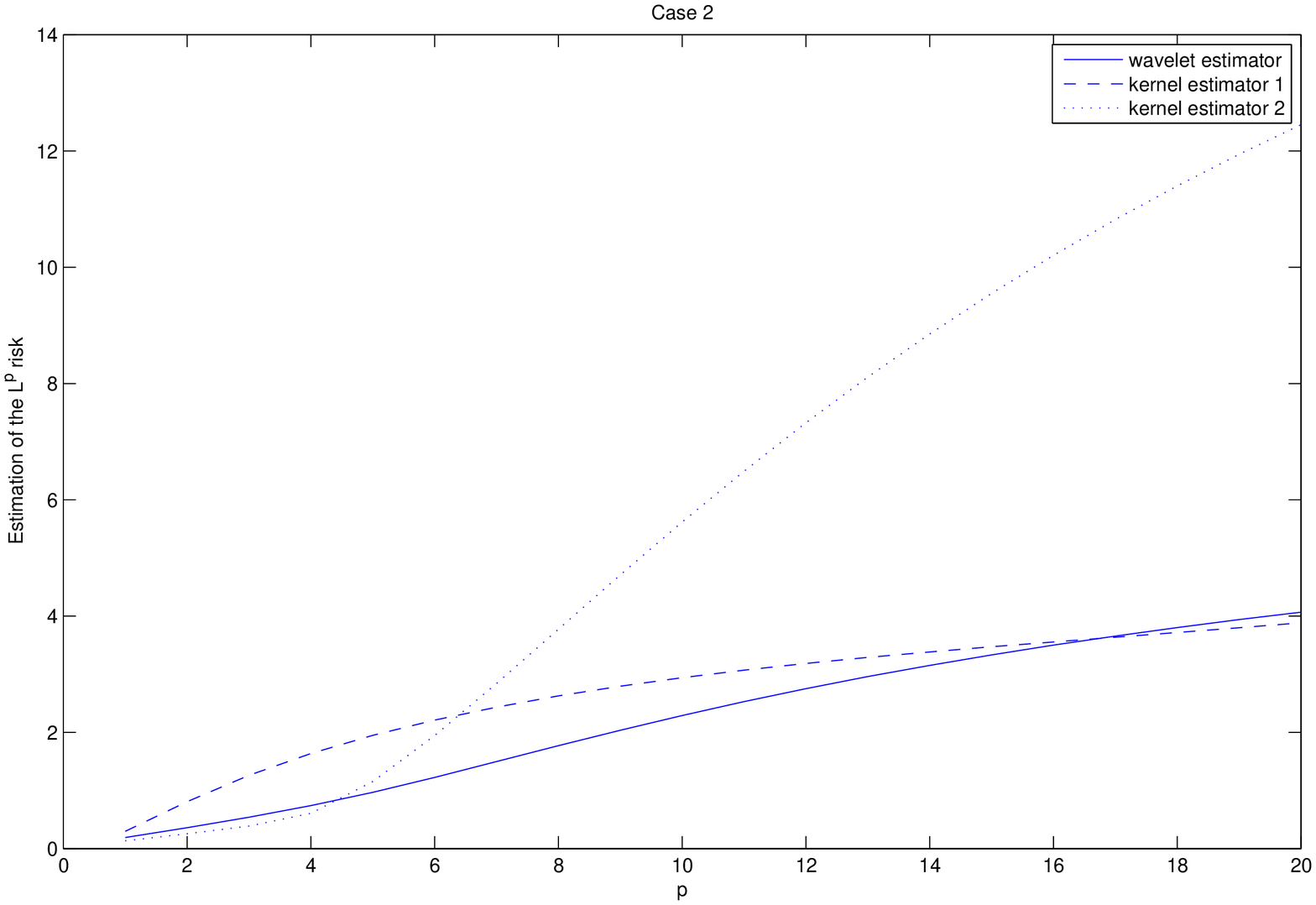} \hspace{-0.5cm}\includegraphics[height=3cm,width=5cm]{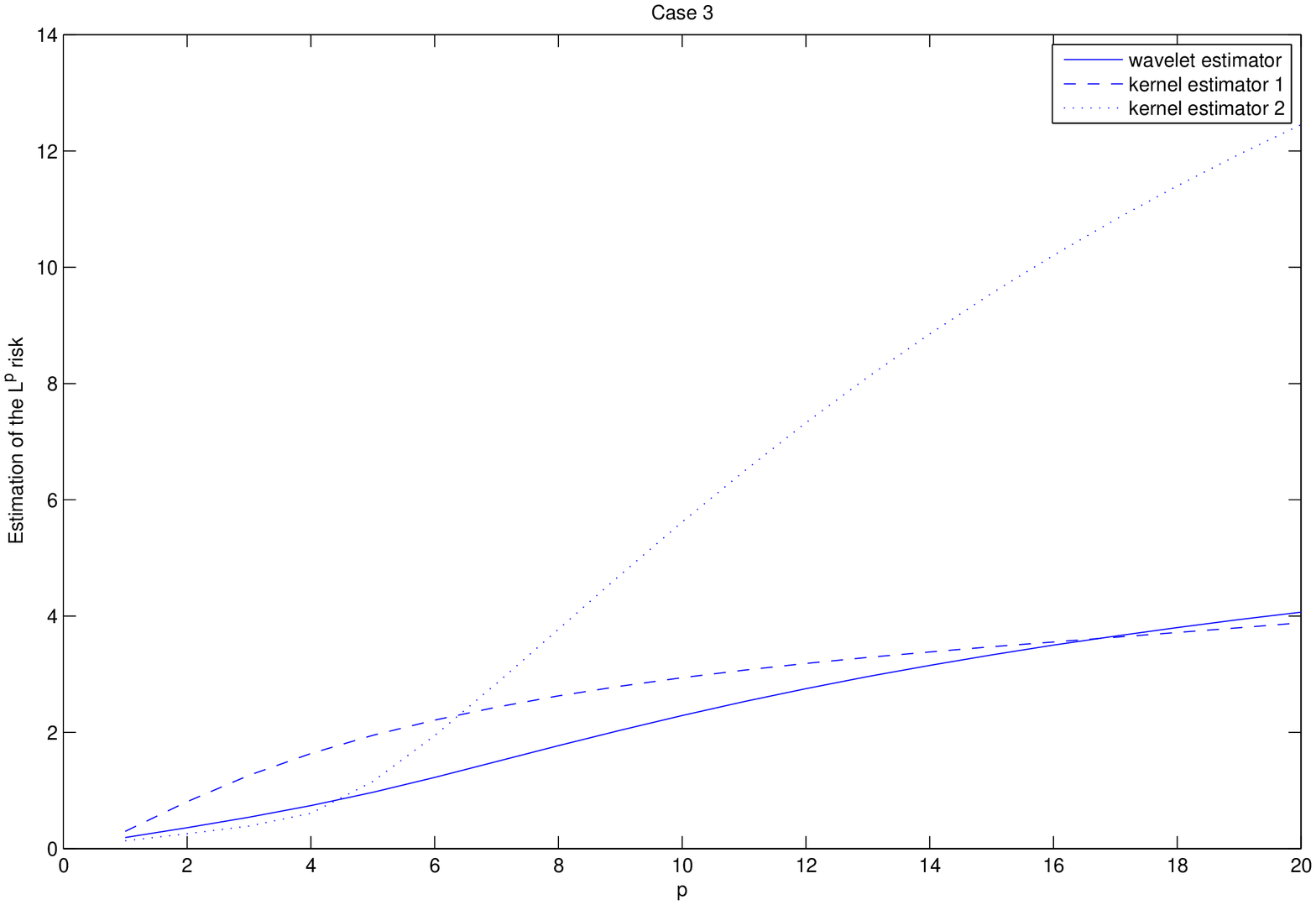}
\end{center}
\vspace{-0.5cm}\caption{{\small Evolution of the $L ^p$ risk of estimators $\widehat f_n^{STCV}$ and kernel estimators obtained on $2^{10}$ observations. The wavelet estimator is represented in dashed lines while the kernel estimators are represented respectively in dots line for the rule of thumb width (kernel estimator 1) and in dots for the cross-validation width (kernel estimator 2). Figures from left to right correspond respectively to {\bf Case~1}, {\bf Case~2} and {\bf Case~3}.}}\label{fig:comparaison}
\end{figure}

To analyse more precisely approximations realized by $\widehat f_n^{STCV}$ and by kernel estimators, we represent in Figure \ref{fig:comparaison} the evolution of the mean $L^p$ risk with respect to $p$ for the three estimators in each case of dependence, i.e. $\E(\|g-f\|_p^p)^{1/p}$ with $g$ equals to one of the three estimators. Even if these risks are  close to each others, kernel estimator with cross validation bandwidth has the smallest risk for small values $p\le 4$. Yet, approximations of this kernel estimator clearly get worse with higher values of $p$, while risks of $\widehat f_n^{STCV}$ seem relatively stable for different value of $p$. Concerning kernel estimator with the width parameter taken according to the rule of thumb, the $L^p$ risk is worse for small values of $p$ but comparable with the one of wavelet estimator for higher values, even if the modes of the density are not detected.

These graphs show that an advantage of $\widehat f_n^{STCV}$ is that its mean $\L^p$ risk seems stable for high values of $p$. Nevertheless, the mean of the $L^2$ risk is larger than the one of kernel estimators with cross-validation width and the computation time higher. One possible way to improve the quality of approximation for the cross validation procedure may be to consider different levels of thresholding at each resolution level. We do not investigate this axis of research as then the time of computation exploded.

\subsection{Different dependent samplings that do not satisfy {\bf (D)}}
In this Subsection, we discuss the necessity of Assumption {\bf (D)}. For this, we study the convergence of the density estimators on some dynamical systems that do not satisfy this assumption. More precisely, we focus on
Liverani-Saussol-Vaienti maps, see \cite{Liverani1999}, defined as the solution of $X_t=T^i(X_{t-i})$ with 
$$
T(x)=\begin{cases} x(1+2^{\alpha'} x^{\alpha'}),&0\le x\le 1/2\mbox{ for some }0<\alpha'<1\\
2x-1,&1/2<x\le 1.
\end{cases}
$$
The process $(X_t)_{t\in\Z}$ is stationary and such that the covariance terms $\cov(f(X_0),g(X_r))$ are of order $r^{1-1/\alpha'}$, see \cite{Young1999} and refinement in \cite{Gouezel2004}. Thus the Assumption {\bf (D)} is not satisfied in this case and we have the non-minimaxity of any thresholded wavelets estimators:
\begin{prop}\label{contrex}
Suppose that the father wavelet $\phi$ is such that $\int\phi>0$ and that the assumptions of Theorem \ref{th} are satisfied. If $1>\alpha'\geq 1/(2\alpha+1)$ with $\alpha$ defined by \eqref{minimaxrate}, then for the thresholded estimators of the marginal density of the Liverani-Saussol-Vaienti map of index $\alpha$ there exists some $C>0$ such that: 
$$
n^{2\alpha}\E[\|\widehat f_n -f\|_2^2]\ge C, \mbox{ for $n$ sufficiently large.}
$$
The same result also holds for the cross validation thresholded estimator $\widehat f_n^{STCV}$.
\end{prop}
The proof of this Proposition is given in Section \ref{proof}.


To simulate these dynamical systems, we simulate $Z_0$ according to the Lebesgue measure on $[0,1]$, then we apply recursively $T$ to determine $Z_i$ and finally  we set $(X_1,\ldots,X_n)=(Z_{n+1},\ldots,Z_{2n})$. This approximation of the stationary solution does not affect the study of the convergence rates as the dynamical system $Z$ is ergodic in mean with rate $O(n^{1-1/\alpha'})$, see Theorem 5 of \cite{Young1999}.\\

\begin{figure}[!ht]
\begin{center}
\hspace{-0.5cm}\includegraphics[height=3cm,width=5cm]{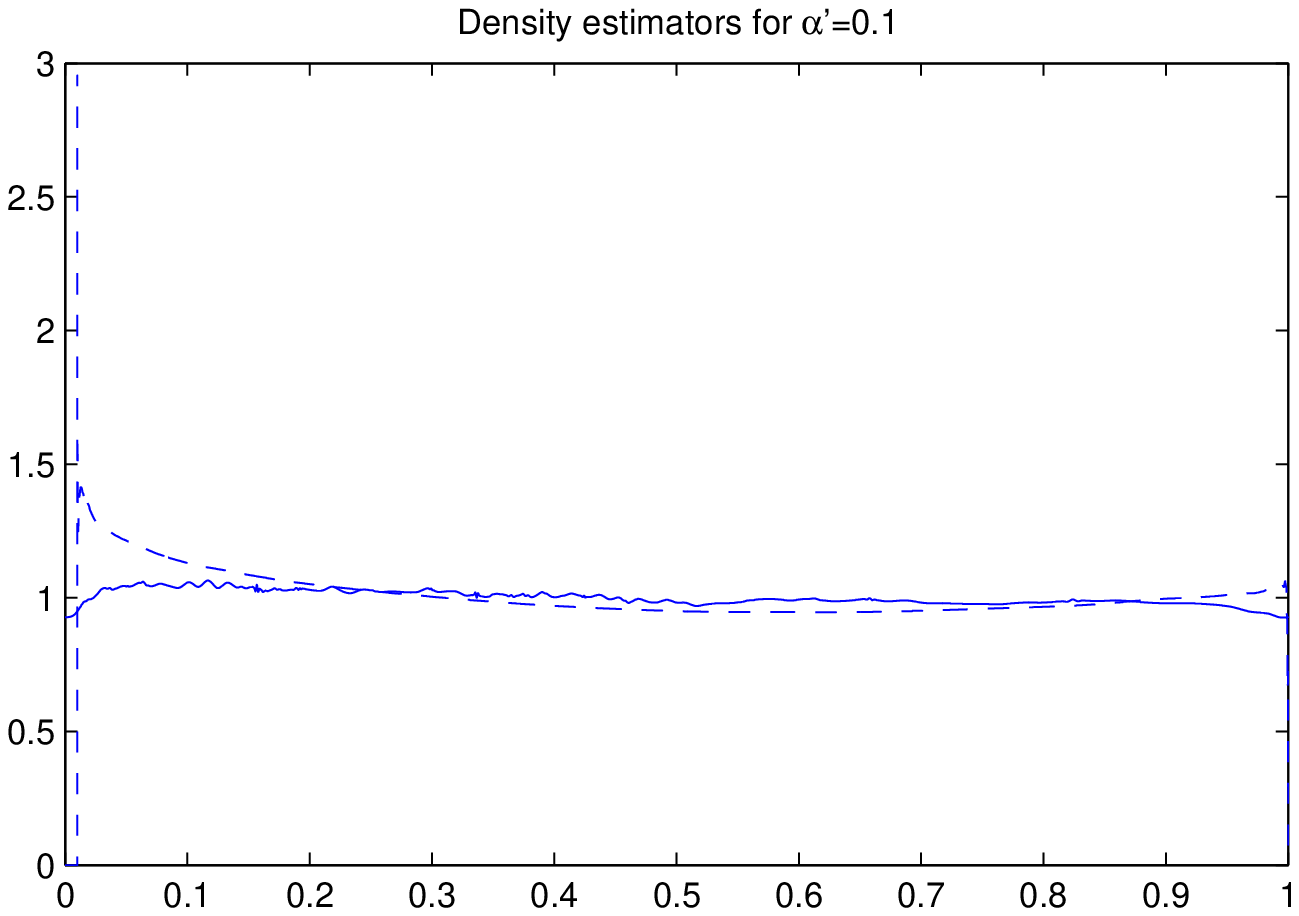}
\hspace{-0.5cm}\includegraphics[height=3cm,width=5cm]{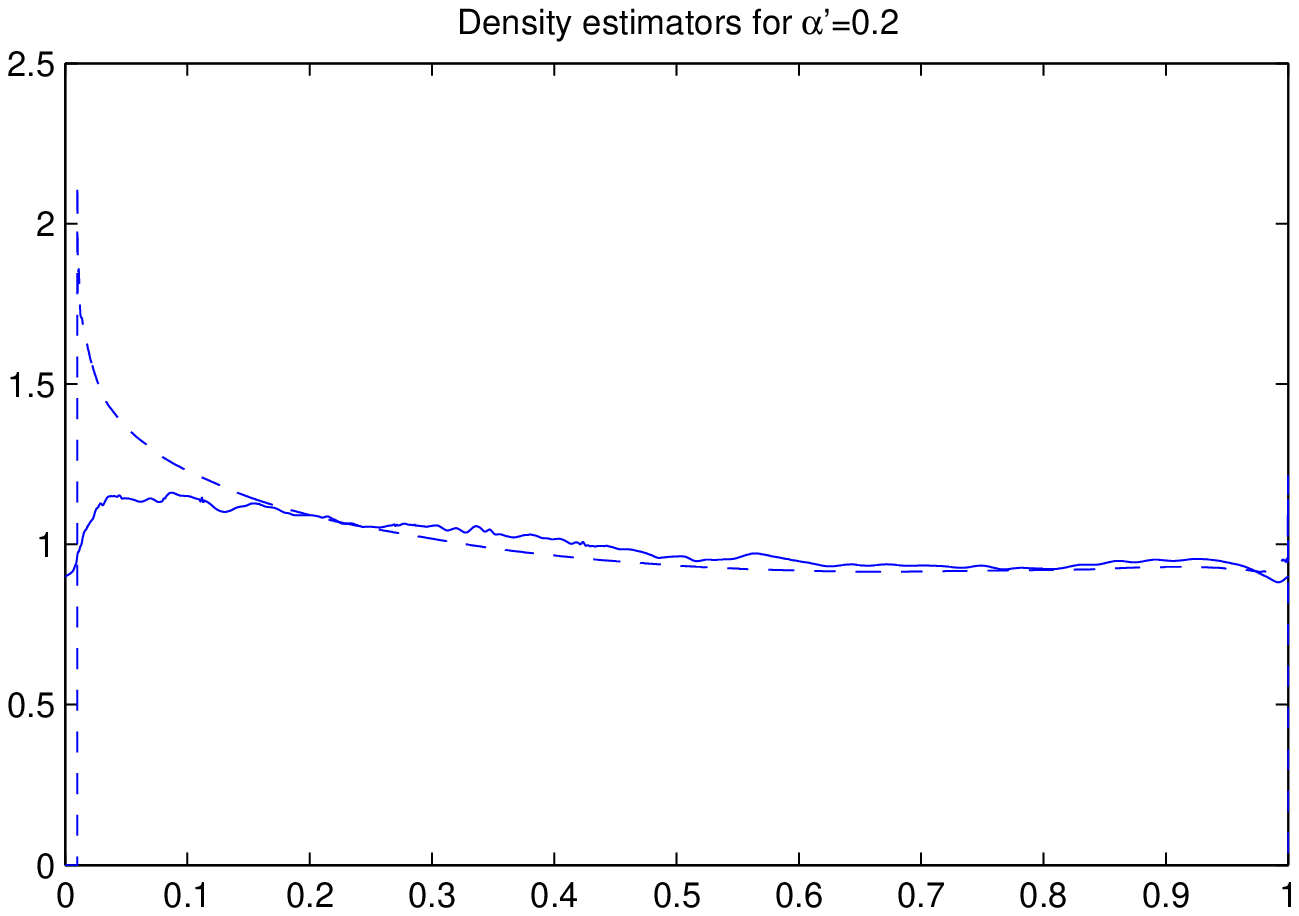} \hspace{-0.5cm}\includegraphics[height=3cm,width=5cm]{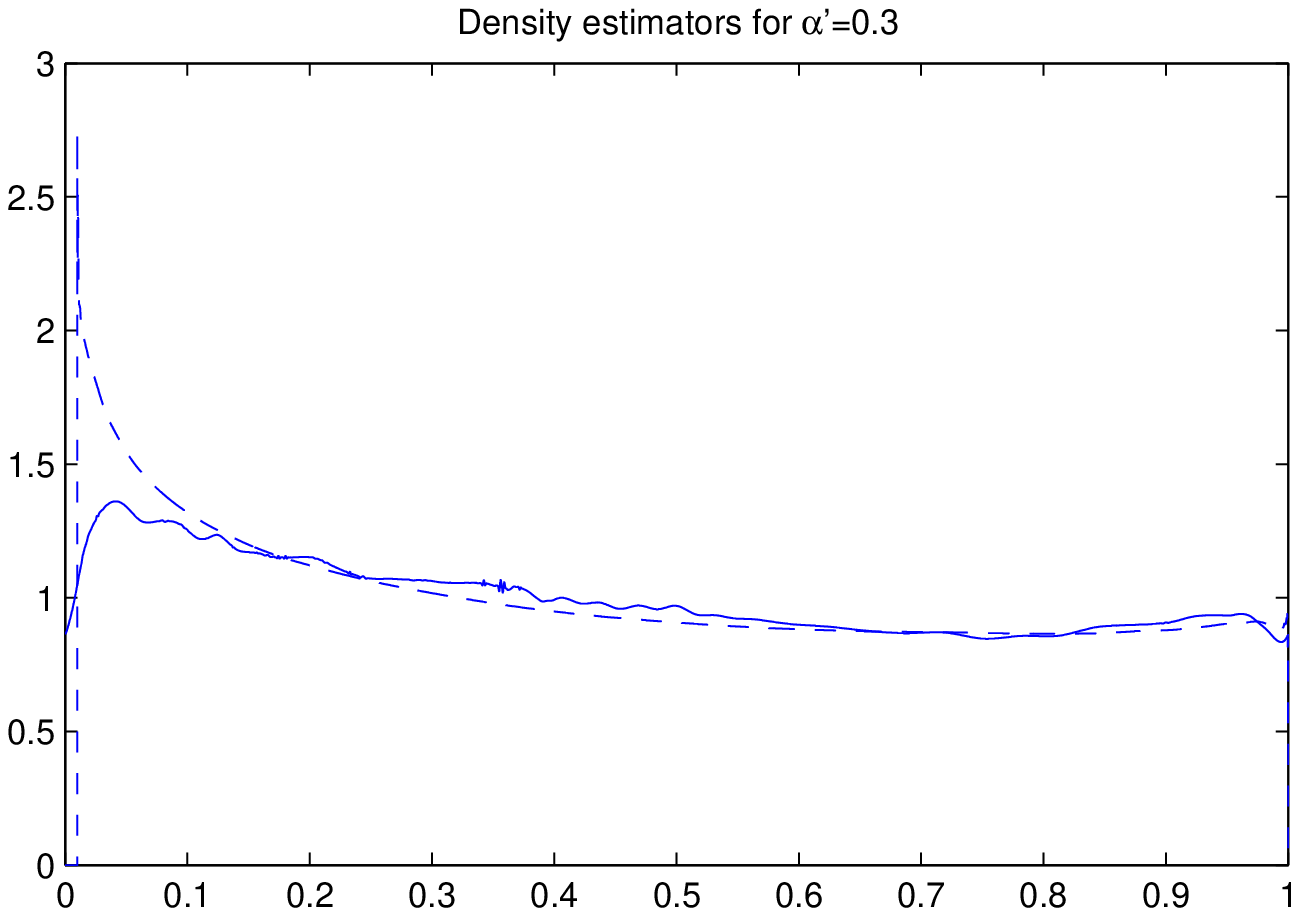}\\
\hspace{-0.5cm}\includegraphics[height=3cm,width=5cm]{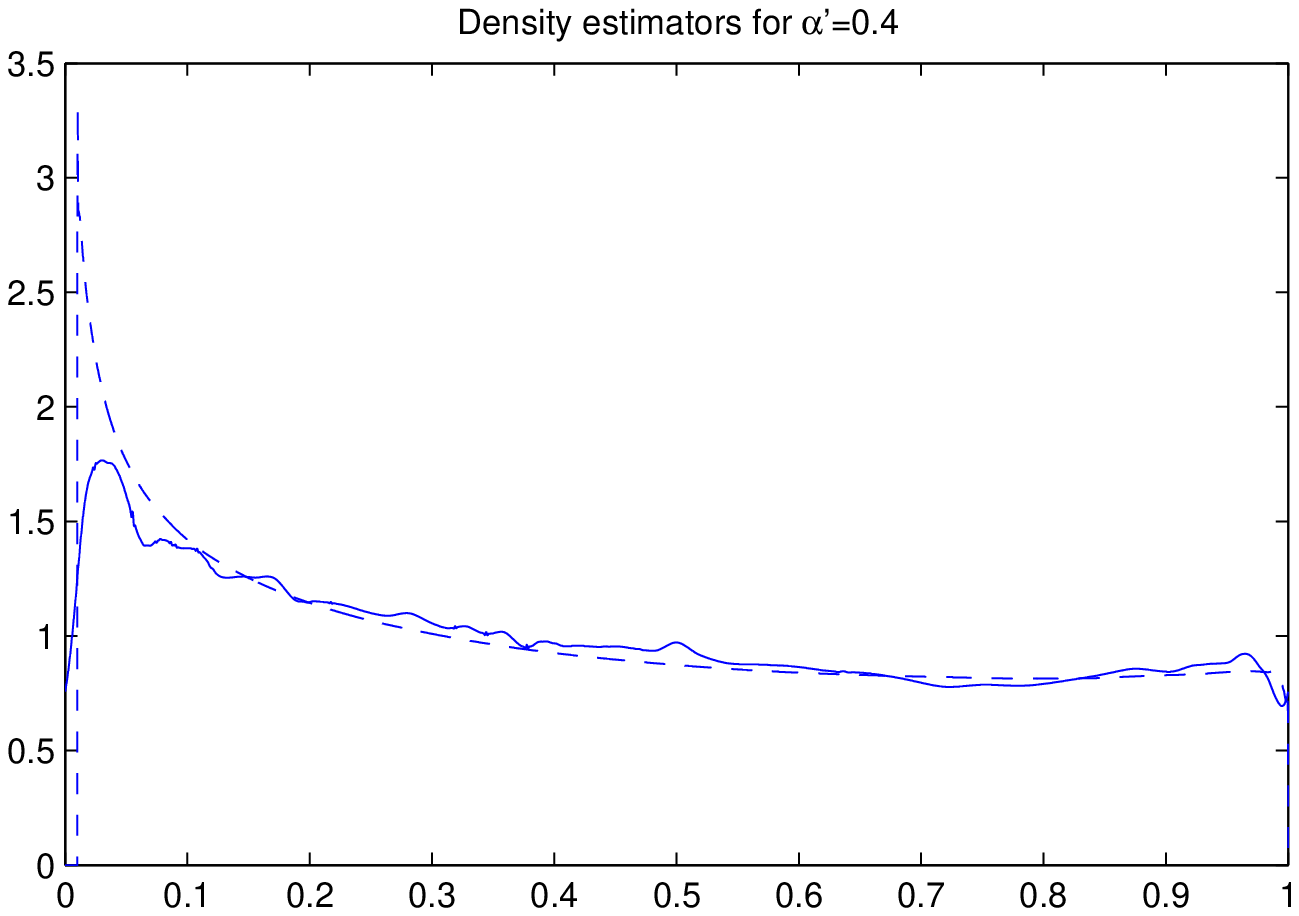} \hspace{-0.5cm}\includegraphics[height=3cm,width=5cm]{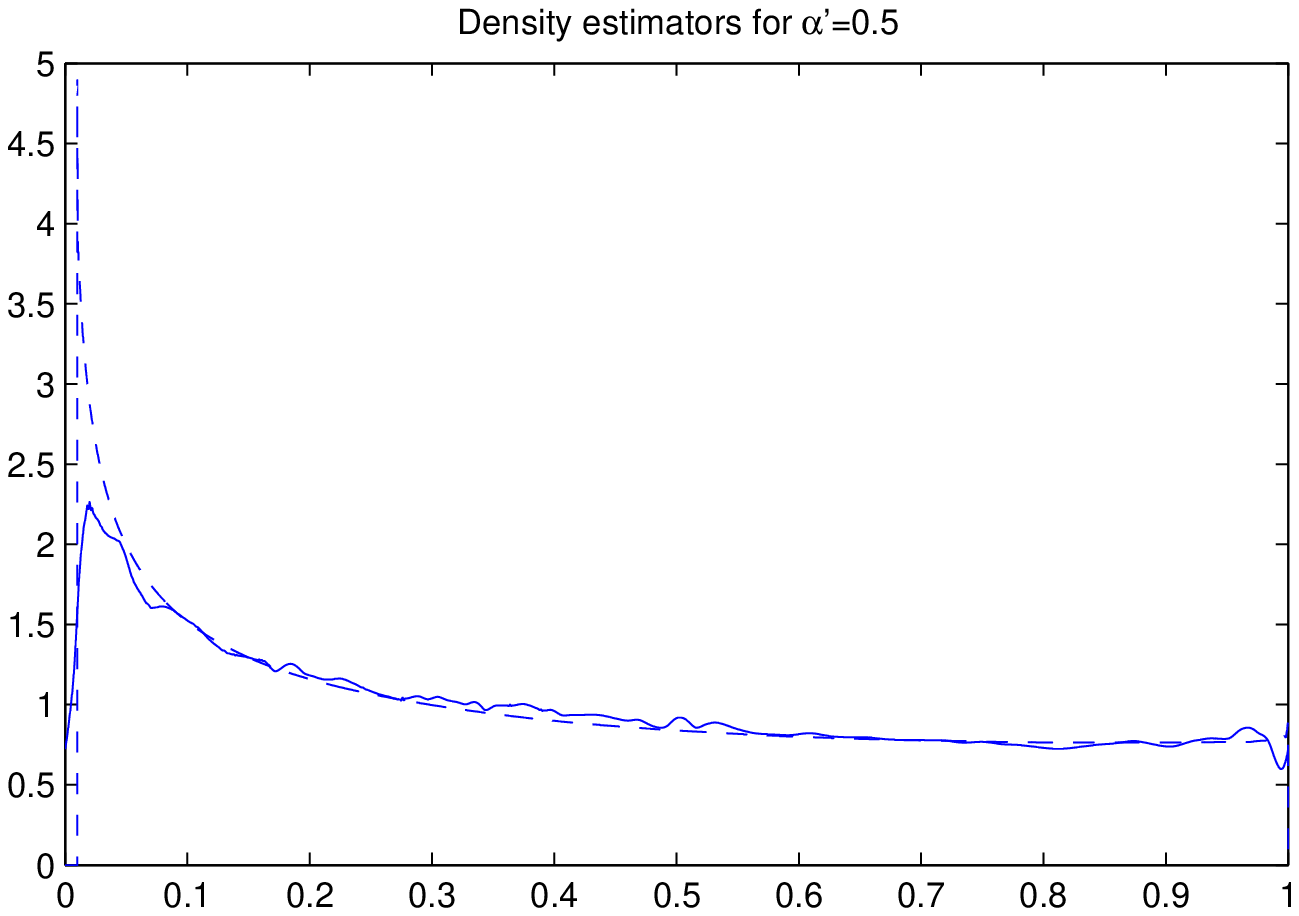} \hspace{-0.5cm}\includegraphics[height=3cm,width=5cm]{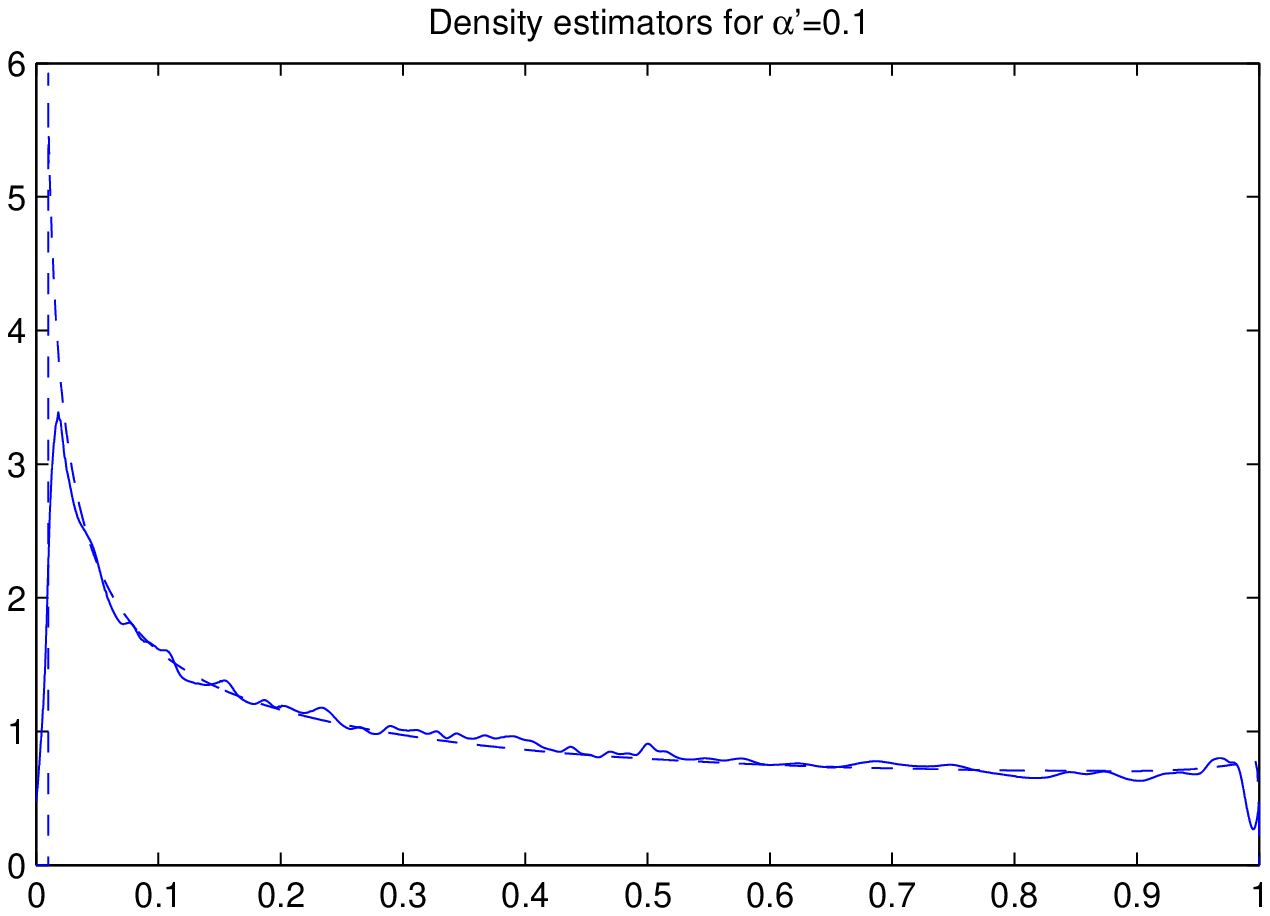}\\
\hspace{-0.5cm}\includegraphics[height=3cm,width=5cm]{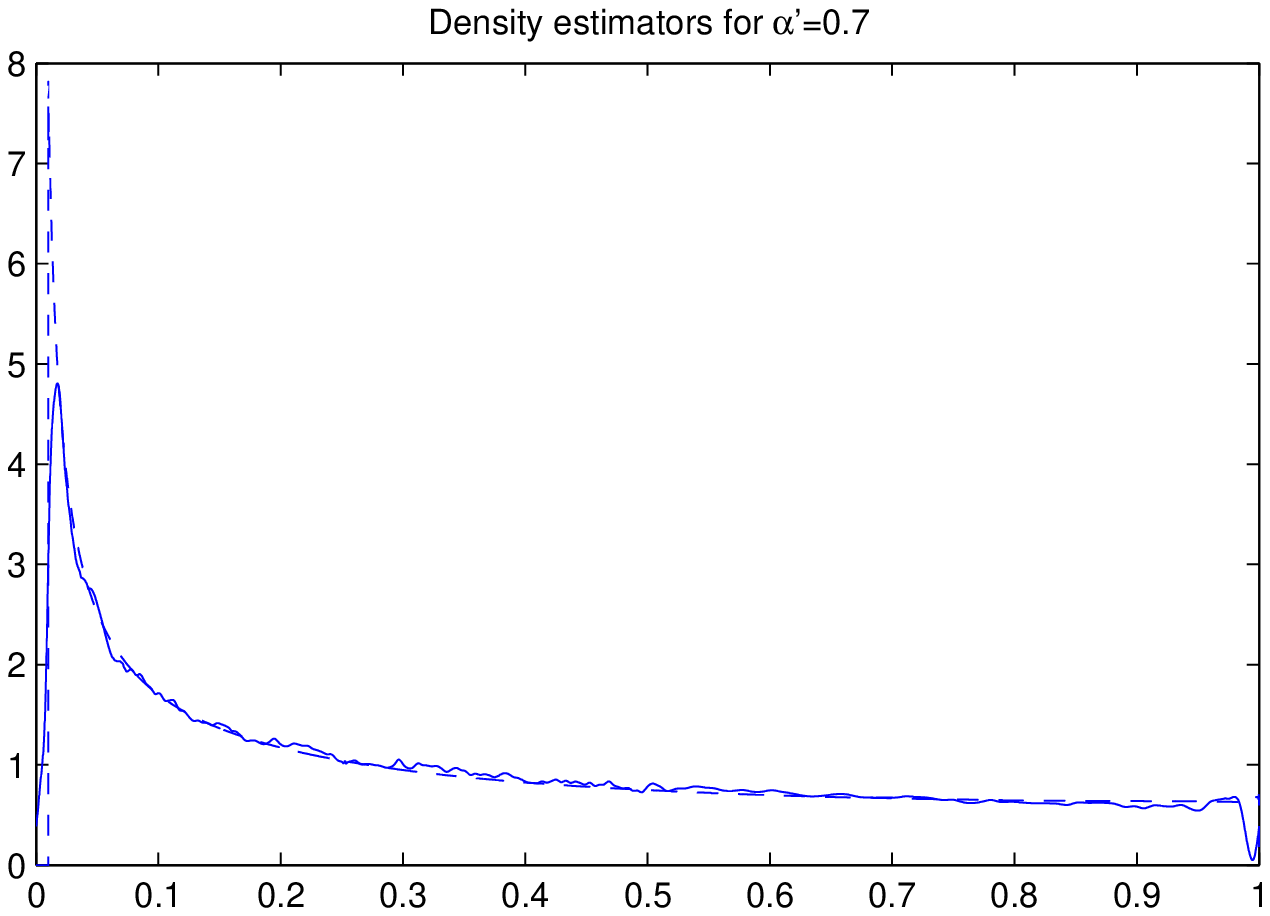} \hspace{-0.5cm}\includegraphics[height=3cm,width=5cm]{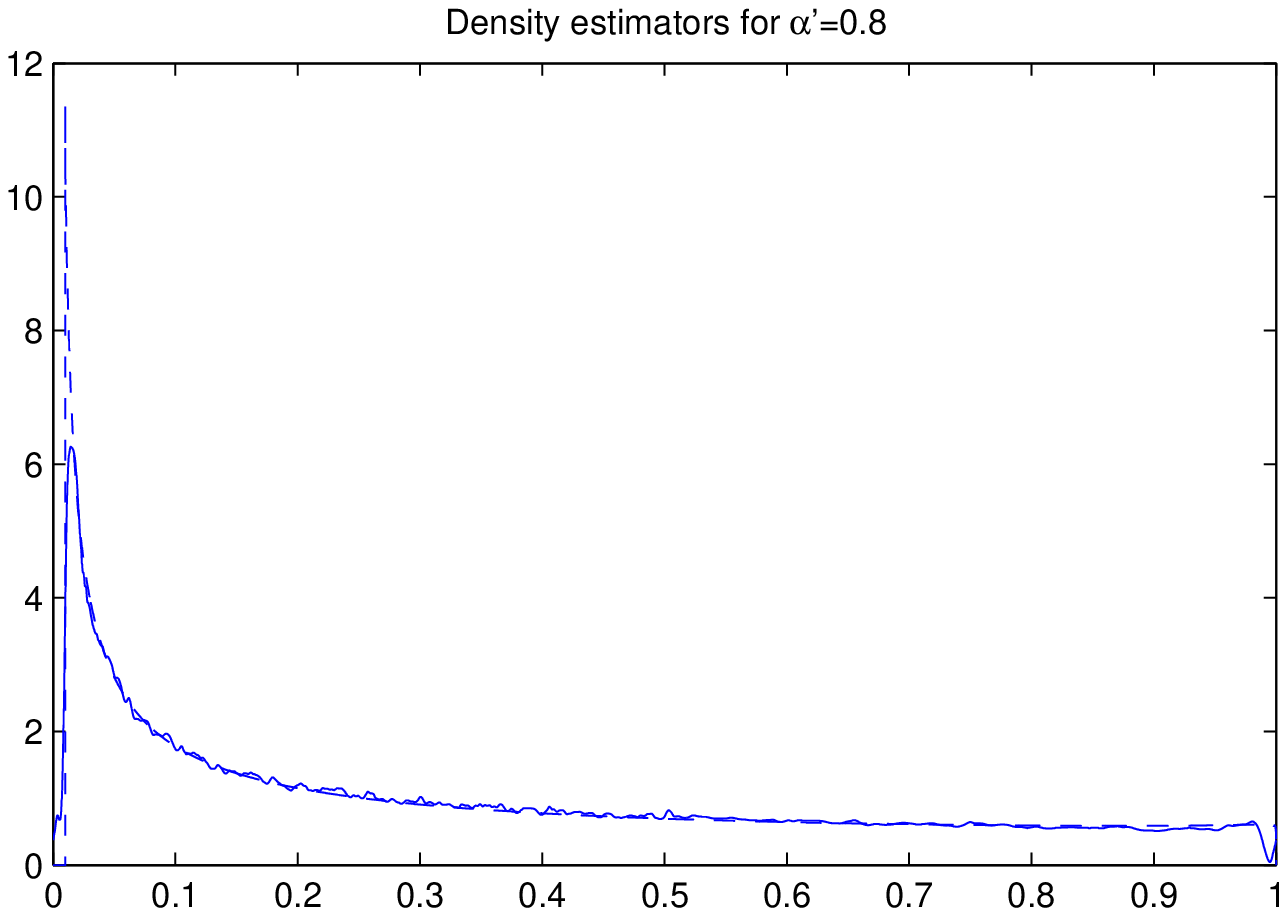} \hspace{-0.5cm}\includegraphics[height=3cm,width=5cm]{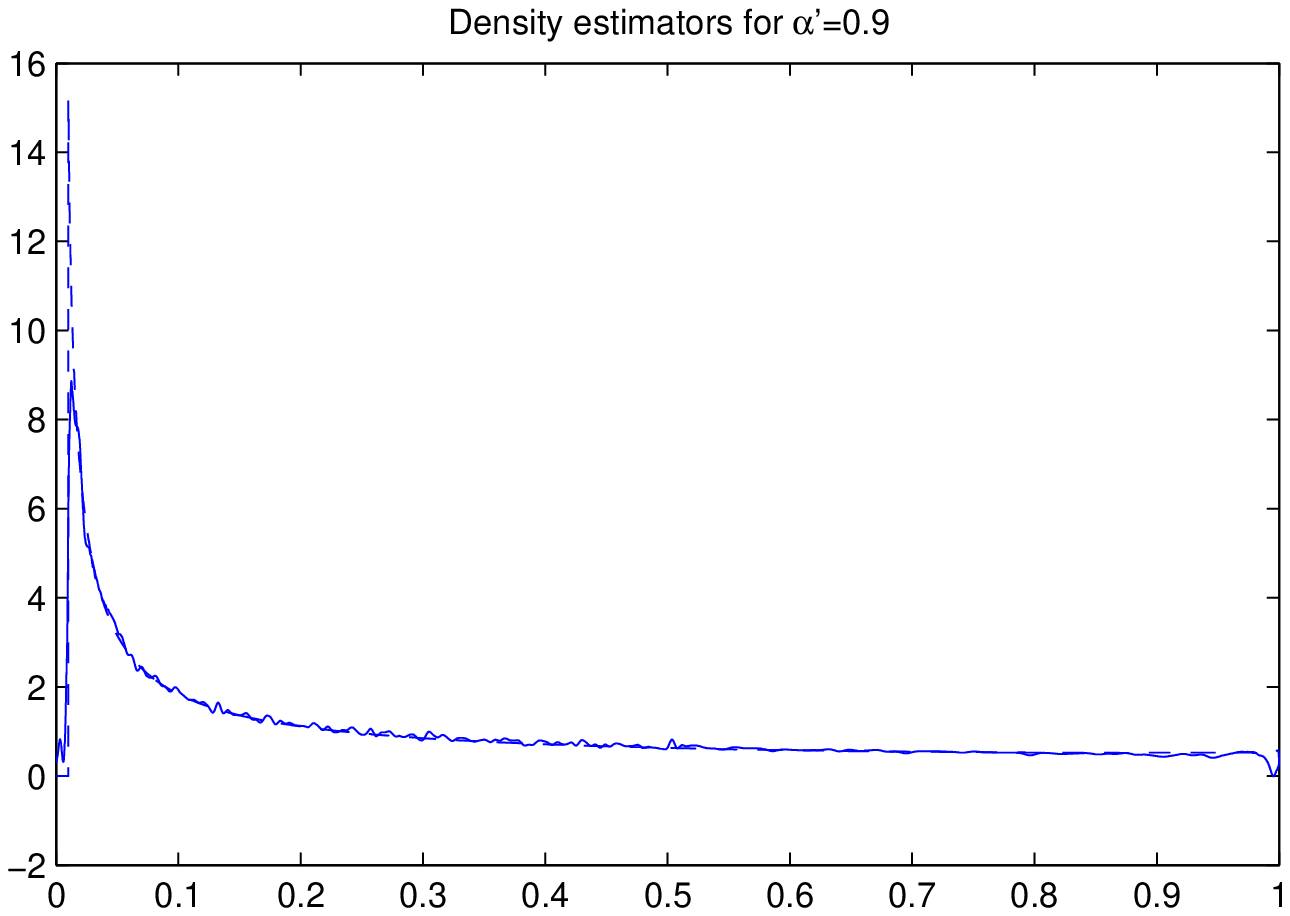}
\end{center}
\vspace{-0.5cm}\caption{{\small Means of the estimators $\widehat f_n^{STCV}$ obtained on $2^{10}$ observations and $500$ simulations. Means of the kernel estimators are represented in dashed lines.}}\label{fig:densitedep}
\end{figure}

The analytic expression of the density $f$ is unknown but it is proved to be continuous, locally Lipschitz and to behave like $x^{-\alpha'}$ as $x\to 0$, see \cite{Liverani1999} and \cite{Young1999}. In particular, $f$ is unbounded on $[0,1]$. Then we restrict our study on $[0.01,1]$ where $f$ is bounded. As the true density $f$ is unknown, we compare here the estimators $\widehat f_n^{STCV}$ with other estimators. In Figure \ref{fig:densitedep}, we plot the mean of $M=100$ estimators $\widehat f_n^{STCV}$  and Epanechnikov kernel estimators given by Matlab's rule of thumb for $9$ different values of $\alpha'$.

Visually, means of both estimators are closed to each other. To detect some difference between the estimators behavior, we decide to compute the moments of order $k=1,\ldots,20$ of the estimators integrated on $[0,1]$:
$$
\int_{0.01}^1\left(\E\left[g^k(t)\right]\right)^{1/k}dt
$$
where the random function $g$ is alternatively $\widehat f_n^{STCV}$ or kernel estimators. 

\begin{figure}[!ht]
\begin{center}
\hspace{-0.5cm}\includegraphics[height=3cm,width=5cm]{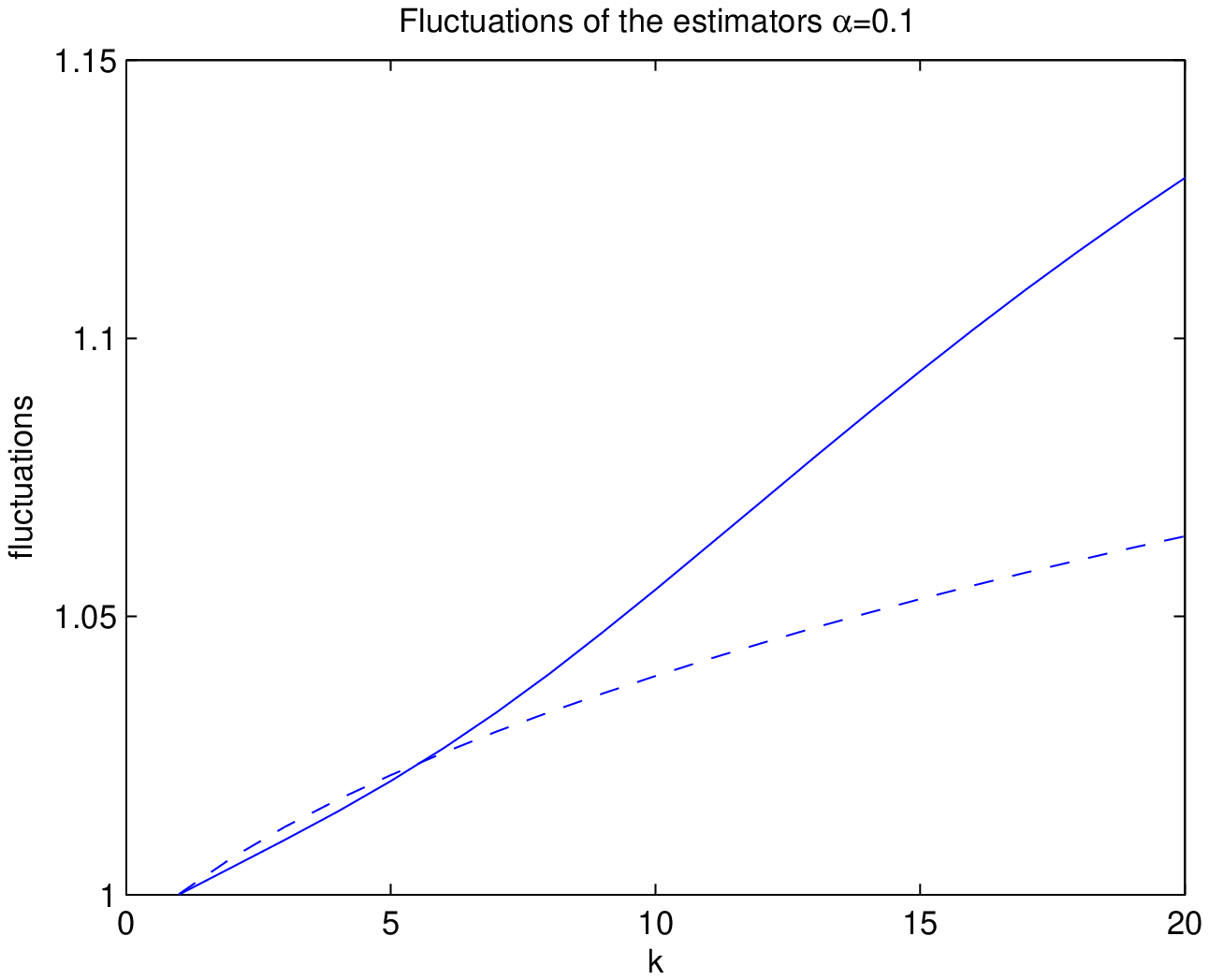}
\hspace{-0.5cm}\includegraphics[height=3cm,width=5cm]{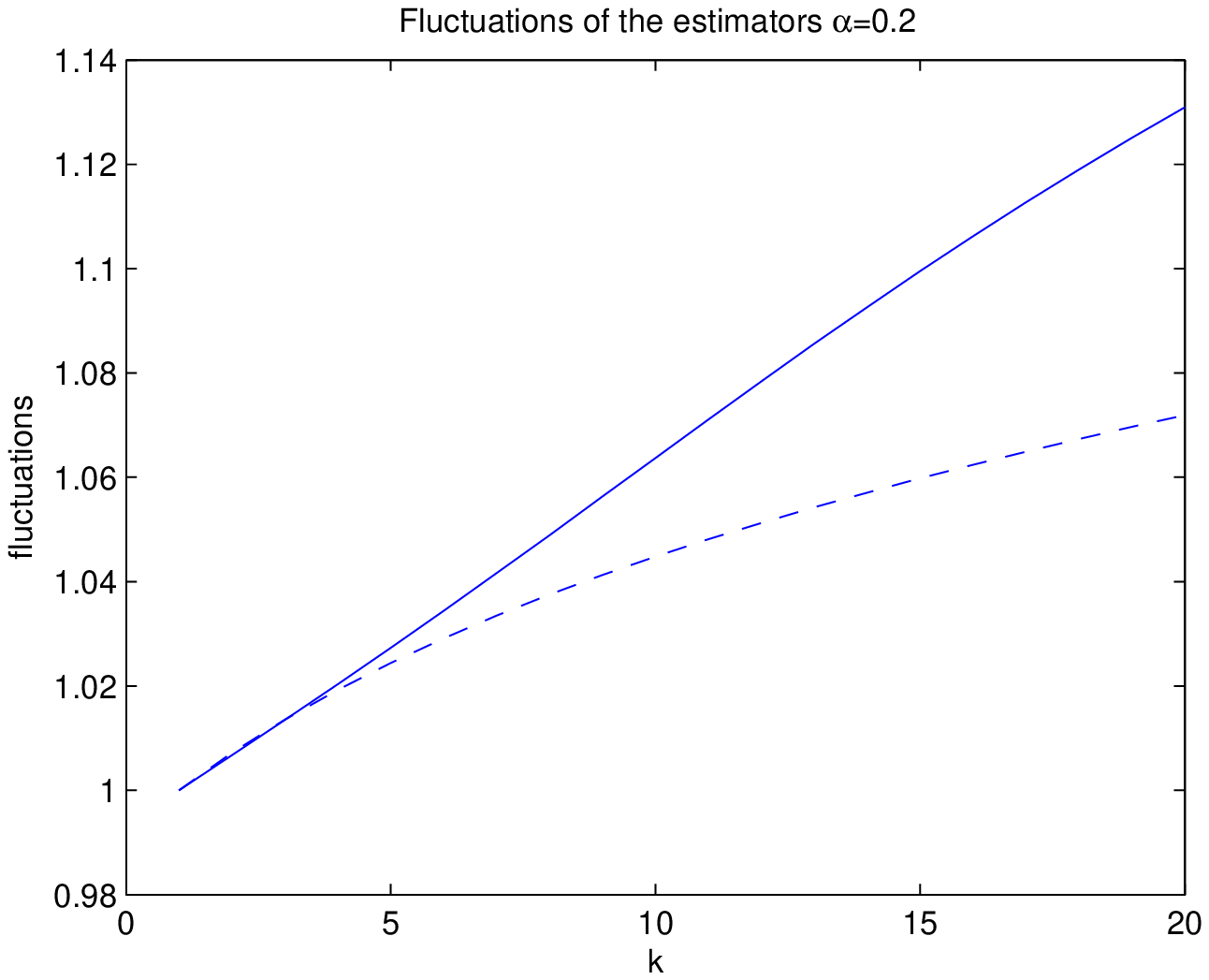} \hspace{-0.5cm}\includegraphics[height=3cm,width=5cm]{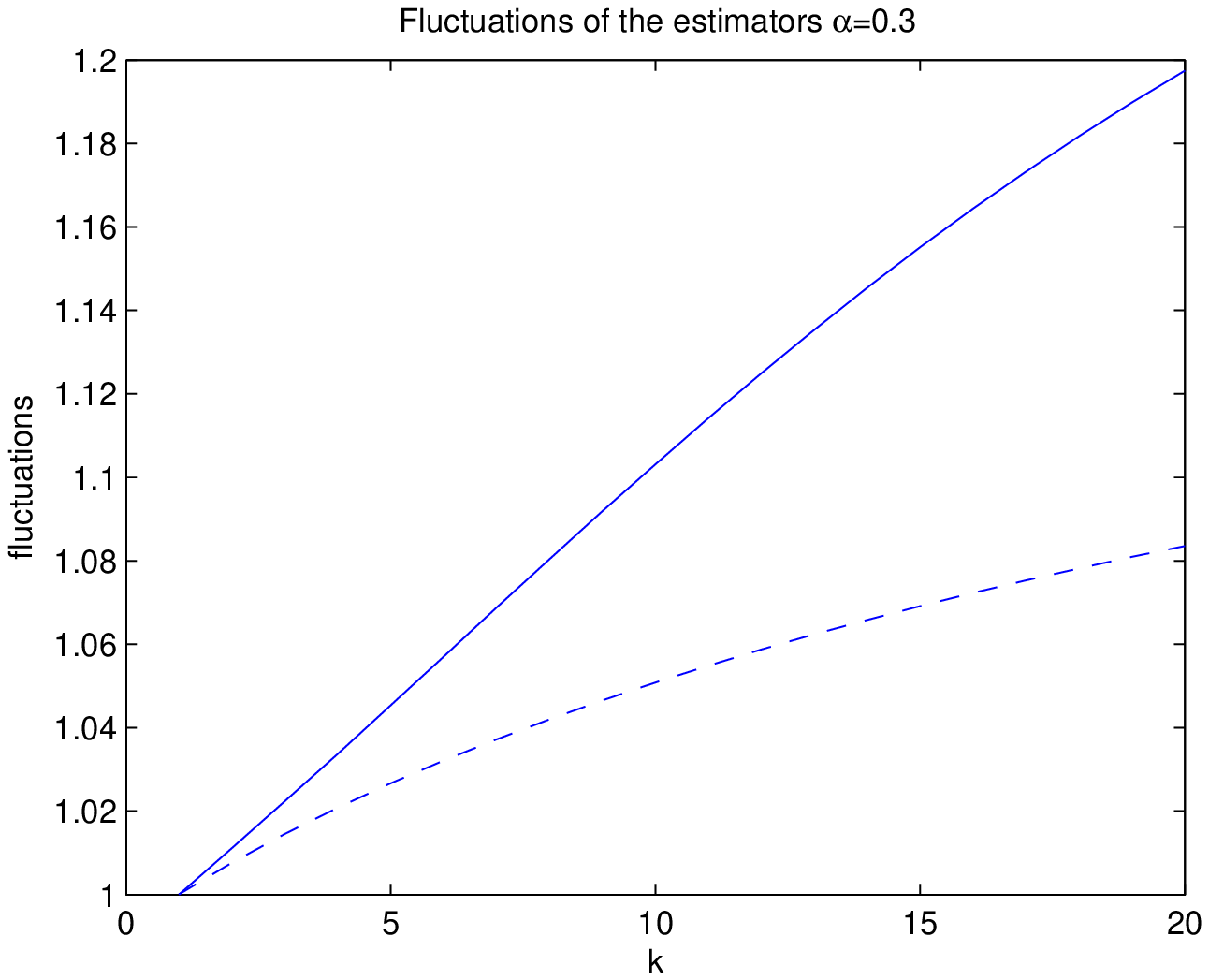}\\
\hspace{-0.5cm}\includegraphics[height=3cm,width=5cm]{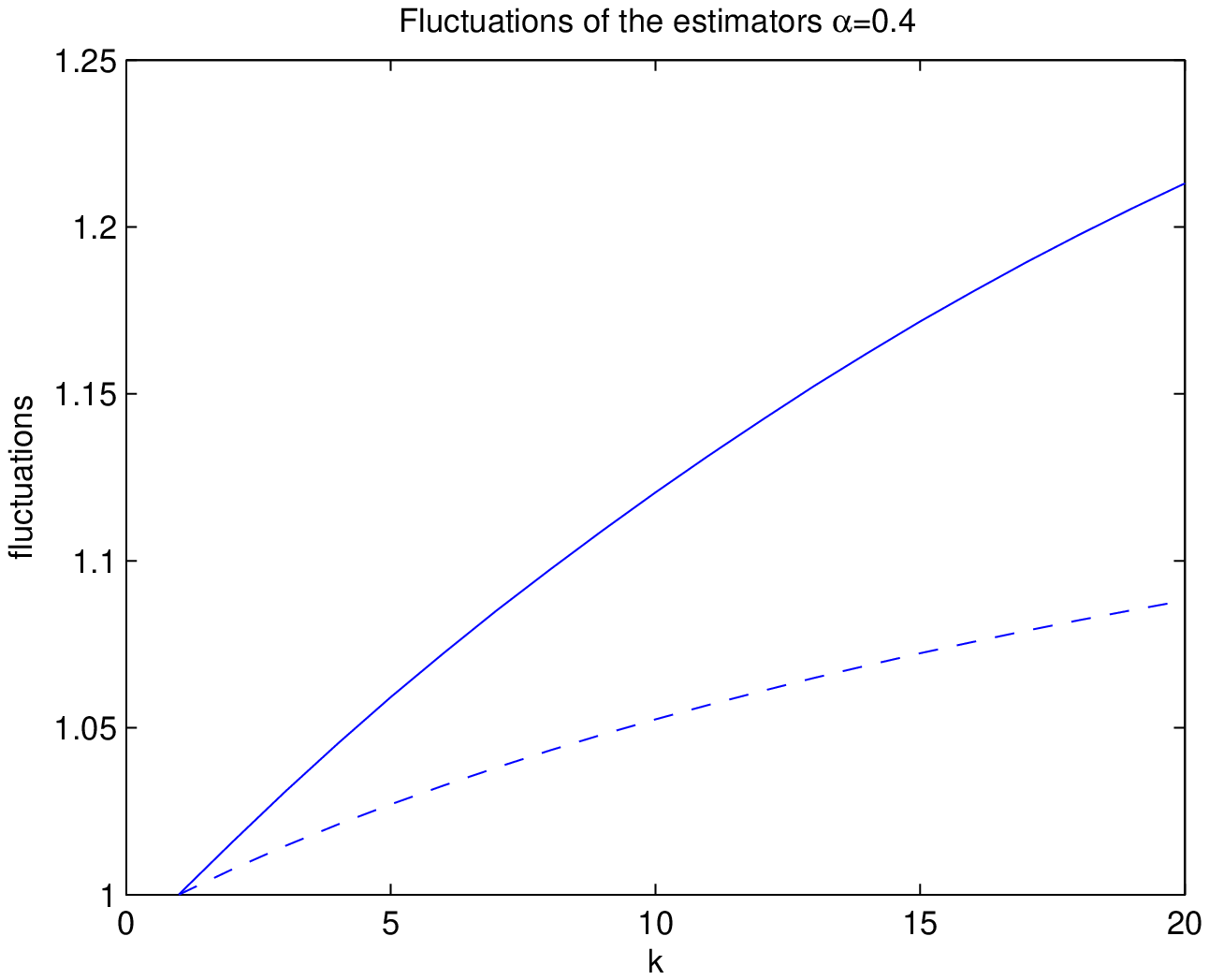} \hspace{-0.5cm}\includegraphics[height=3cm,width=5cm]{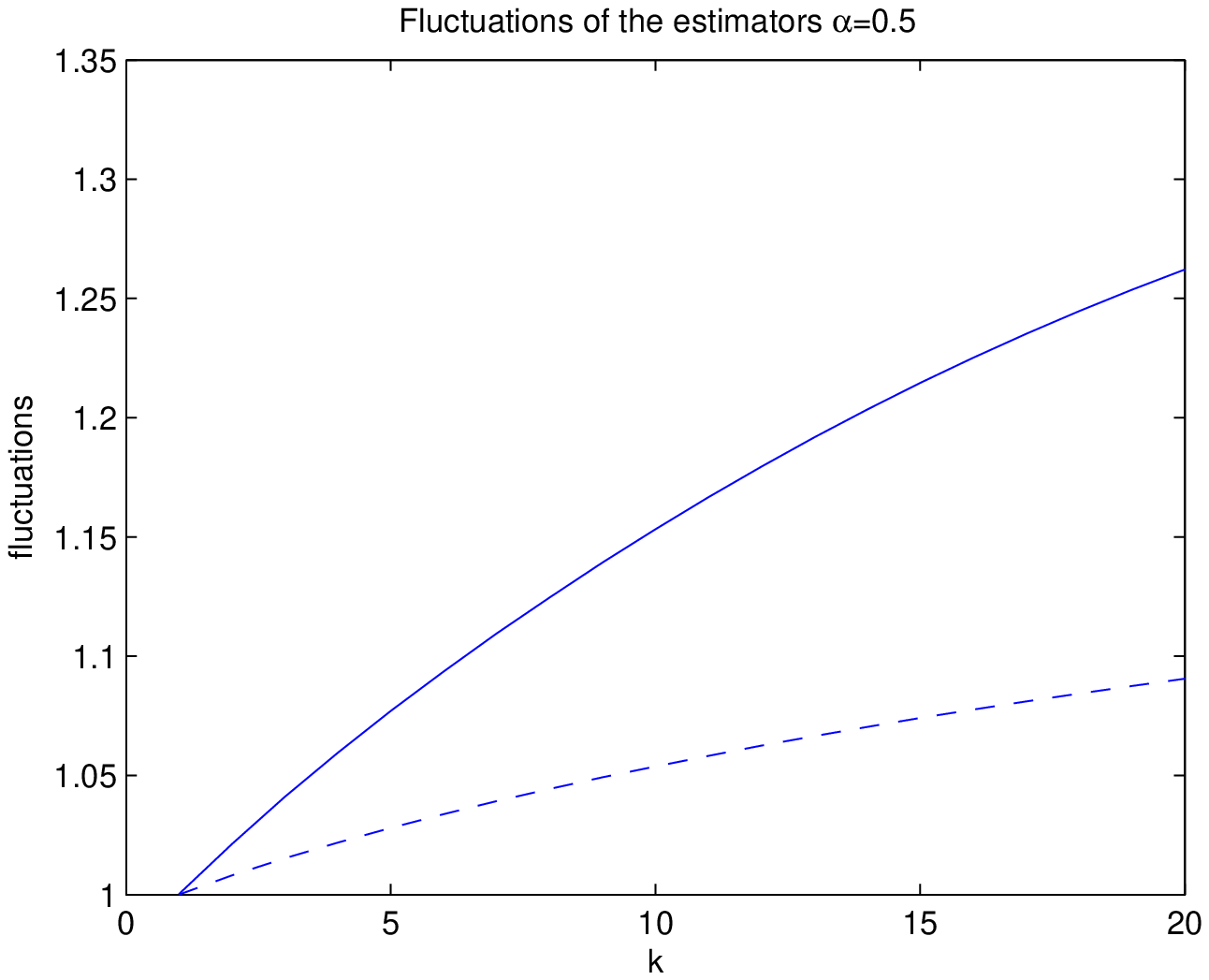} \hspace{-0.5cm}\includegraphics[height=3cm,width=5cm]{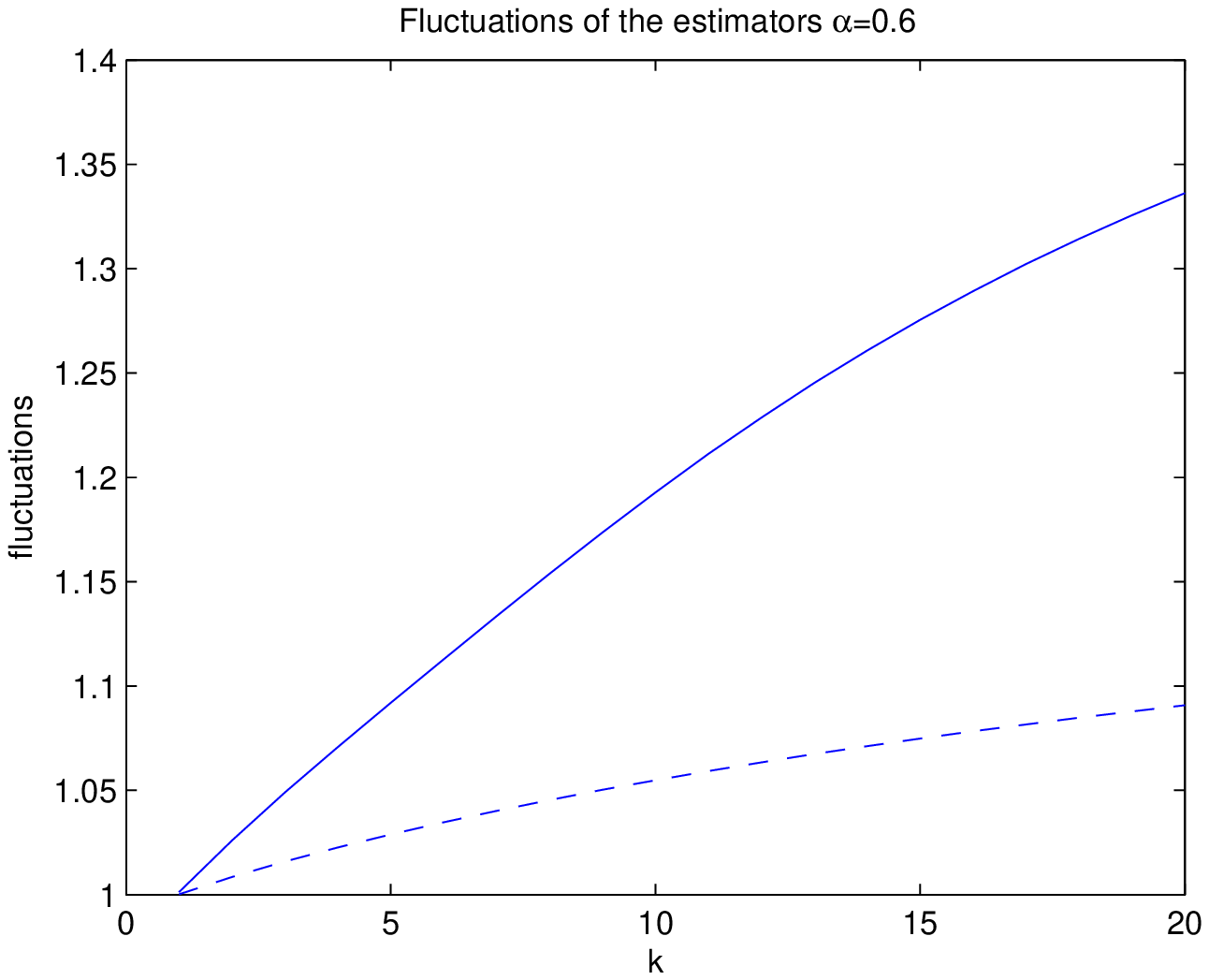}\\
\hspace{-0.5cm}\includegraphics[height=3cm,width=5cm]{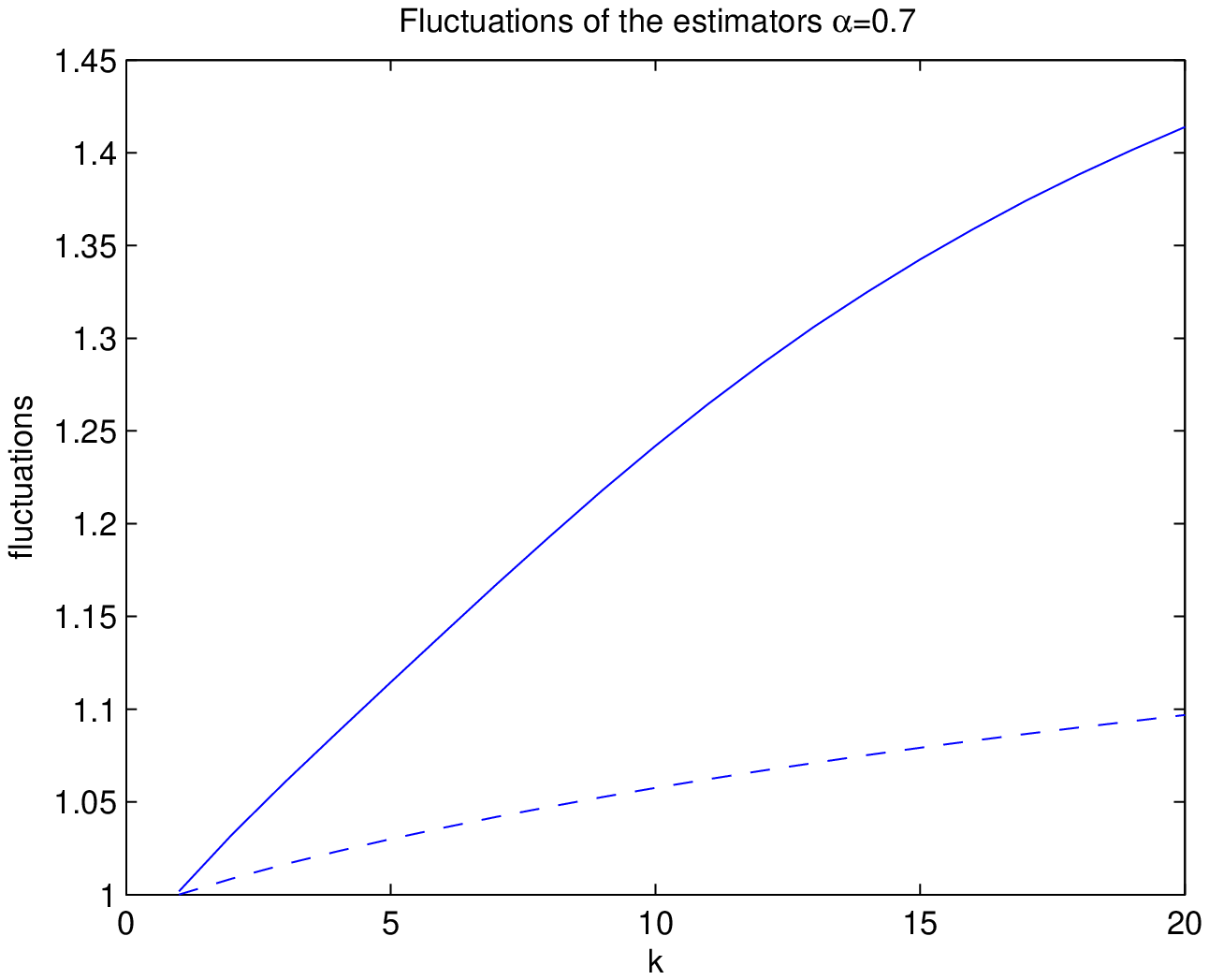} \hspace{-0.5cm}\includegraphics[height=3cm,width=5cm]{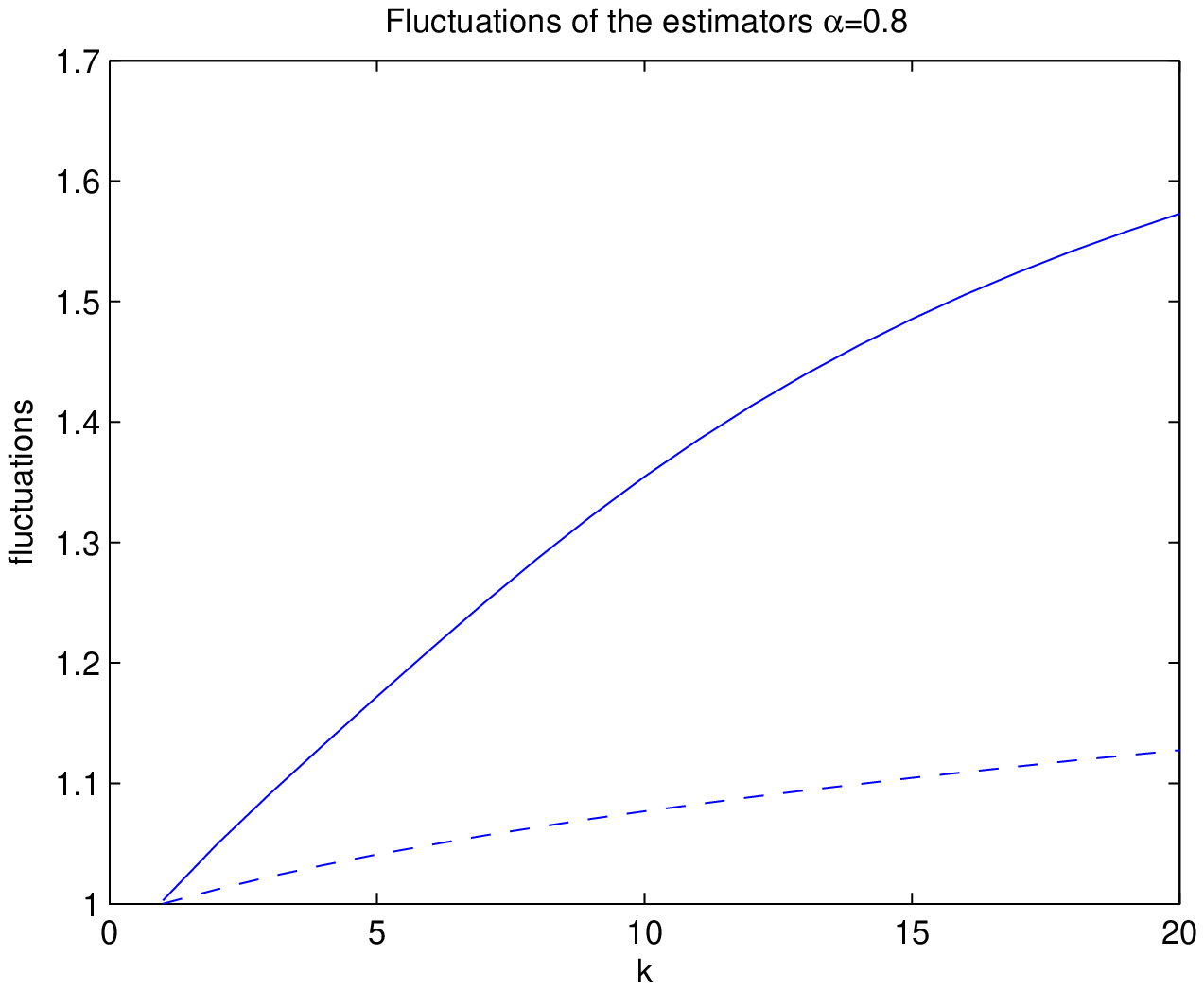} \hspace{-0.5cm}\includegraphics[height=3cm,width=5cm]{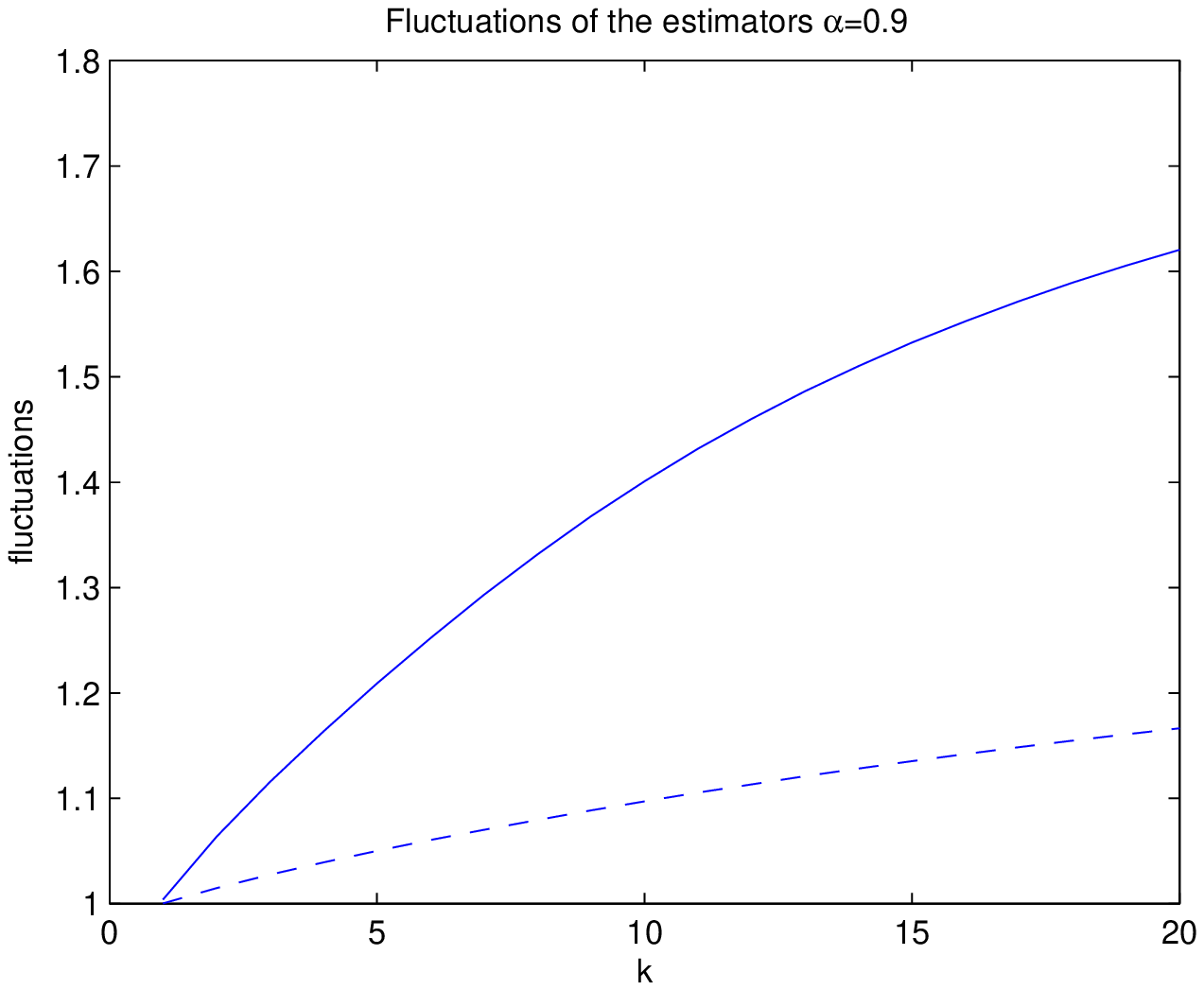}
\end{center}
\vspace{-0.5cm}\caption{{\small Moments of the estimators $\widehat f_n^{STCV}$ obtained on $2^{10}$ observations and $500$ simulations. In dashed lines are represented the moments of the kernel estimators.}}\label{fig:compdep}
\end{figure}

For small values $\alpha'\le 0.02$ and $k\le 4$, the moments of both estimators have similar values. But as $\alpha'$ growths, all the moments of $\widehat f_n^{STCV}$ explode more rapidly than the ones of kernel estimators as $k$ increases. Previous simulations studies show no behavior difference for $\widehat f_n^{STCV}$  between independent cases and dependent cases satisfying {\bf (D)}. In dependent cases that do not satisfy {\bf (D)}, the behavior of $\widehat f_n^{STCV}$ depends on the decrease rates of covariance terms. On the contrary, the behavior of the kernel estimators with rule of thumb width is more stable when {\bf (D)} is not satisfied.

\section{Proofs}\label{proof}
In this Section are collected all the proofs of this paper.

\subsection{Proofs of Lemma \ref{useful}}\label{prooflem}

Firstly we give the proof of the inequality \eqref{Rosenthal} of Lemma \ref{useful}. This proof is essentially based on Lemmas \ref{douklou} and \ref{L1L2}. The first one, Lemma \ref{douklou}, is simply the Theorem 2 of \cite{DoukhanLouhichi} that we recall here for completeness:
\begin{lem}{(Doukhan \& Louhichi (1999))}\label{douklou}
Let $(Z_i)_{1\le i\le n}$ be centeres variables.\\
Let $q$ be an even integer and $n\geq2$.\\
Suppose that for all $p = 2, \dots,q$ and for all $1\le s_1\le s_p\le n$
satisfying $\max s_{i+1}-s_i= s_{u+1}-s_u=r$, there exists $V_{p,n}$ such that:
$$
n \sum_{r=0}^{n-1}(r+1)^{p-2} \cov(Z_{s_1}\cdots
Z_{s_u},Z_{s_{u+1}}\cdots Z_{s_p})\le V_{p,n}.
$$
Then, we have \begin{equation}
\E|\sum_{i=1}^nZ_i|^q\le
\frac{(2q-2)!}{(q-1)!}\left\{V_{2,n}^{q/2}\vee V_{q,n}\right\}\;,
\end{equation}
\end{lem} We refer the reader to \cite{DoukhanLouhichi} for the proof of this result.
We apply this result on $Z_i=\psi_{j,k}(X_i)-\beta_{j,k}$, when {\bf
(D)} holds. We first determine the bounds $V_{p,n}$ as under Assumption {\bf
(D)} we have:
$$
\sum_{r=0}^{n-1}(r+1)^{p-2} \cov(Z_{s_1}\cdots
Z_{s_u},Z_{s_{u+1}}\cdots Z_{s_p})\leq \sum_{r=0}^{n-1}(r+1)^{p-2}\rho(r)p^2 (2^{j/2}M_2)^{p-2}.$$
Here we need the following analytic Lemma to bound the quantity $\sum_{r=0}^\infty(r+1)^p\rho(r)$:
\begin{lem}\label{L1L2}
If \eqref{condrho} is satisfied, {\it i.e.} if $\rho(r)\leq C_0 e^{-ar^b}$, then for
all integer $p$ we have
\begin{equation}\label{defmu}
\sum_{r=0}^\infty(r+1)^p \rho(r)\le C_1 C_2^p (p!)^{1/b},
\end{equation}
with some constants $C_1$ and $C_2$ that are depending on $a$ and $b$.
\end{lem}
Applying this result, whose proof is given at the end of this subsection, we directly obtain the new bound:
$$
\sum_{r=0}^{n-1}(r+1)^{p-2} \cov(Z_{s_1}\cdots \leq C_1 p^2(2^{j/2}M_2 C_2)^{p-2}((p-2)!)^{1/b}
$$
We set $
V_{p,n}=C_1 p^2(2^{j/2}M_2 C_2)^{p-2}((p-2)!)^{1/b} n,$
and, applying Lemma \ref{douklou}, we obtain that:
$$\E |\sum_{i=1}^n(Z_i-\E Z_0)|^q\le \frac{(2q-2)!}{(q-1)!}\left\{(C_1n)^{q/2}\right\}\vee\left\{ C_1(2^{j/2}M_2C_2)^{q-2}((q-2)!)^{1/b}\right\}.$$
Dividing by $n^{-q}$ and noticing that $2^{j/2}\le n$ for $0\le j\le \log n$, we derive that:
$$\E |\widehat\beta_{j,k}-\beta_{j,k})|^q\le n^{q/2}\frac{(2q-2)!}{(q-1)!}\left\{(C_1)^{q/2}\right\}\vee\left\{ C_1(M_2C_2)^{q-2}((q-2)!)^{1/b}\right\},$$ which corresponds to the inequality \eqref{Rosenthal}. In particular, for $p=2$ we have
$$
\E|\sum_{i=1}^n(Z_i-\E Z_0)|^2\le4C_1n.
$$

Now the inequality \eqref{bernstdep} of Lemma \ref{useful} is a direct application of Theorem 1 in \cite{DoukhanNeumann} with $S_n=\sum_{i=1}^n(Z_i-\E Z_0)$, $t=n\lambda$, $\nu=0$, $\mu=1/b$, $A_n=4 C_1n$ and
$B_n=2M_2 C_2 2^{(2+b)/b}2^{j/2}$. We refer the reader to \cite{DoukhanNeumann} for the
definition of the parameters $t,\nu,\mu,A_n$ and $B_n$.

\begin{proof}[Proof of Lemma \ref{L1L2}]
Define $g(x)=(1+x)e^{-ax^b}$ for all $x\le 0$. Studying its derivative, we can easily see that it exists $x_{a,b}$ such that the function $g$ decreases on $[x_{a,b},+\infty)$. If we denote $k\geq1$ the smallest integer greater than $x_{a,b}$, we can infer from \eqref{condrho} that
\begin{eqnarray*}
\sum_{r=0}^\infty(r+1)^p \rho(r) & \le & C_0\left(\sum_{r=0}^{k-1}(r+1)^p \rho(r)+\int_{k-1}^\infty(x+1)^p\exp(-ax^b)dx\right)\\
& \le & C_0\left(C_{a,b}+\int_0^\infty(x+1)^p\exp(-ax^b)dx\right).\end{eqnarray*}
With a convex inequality on $x\mapsto x^p$, we achieve the bound
$$
\sum_{r=0}^\infty(r+1)^p \rho(r) \le C_0\left(C_{a,b}+2^{p-1}\left(\int_0^\infty x^p\exp(-ax^b)dx+\int_0^\infty \exp(-ax^b)dx\right)\right).
$$
Then, writing $u=ax^b$,
$$
\sum_{r=0}^\infty(r+1)^p \rho(r) \le C_0\left(C_{a,b}+2^{p-1}b^{-1}\left[a^{-(p+1)/b}\Gamma\left(\frac{p+1}b\right)+a^{-1/b}\Gamma\left(\frac{1}{b}\right)\right]\right),
$$
with $\Gamma$ defined by $\Gamma(x)=\int_0^\infty u^{x-1}\exp(-u) du$ for all $x>0$.
Let $\lceil x\rceil$ denotes the largest integer smaller than $x>0$ and note $\overline \Gamma =\sup_{x\in]0,1]}\Gamma(x)$. Using the inequalities $\Gamma((p+1)/b)\le \overline \Gamma\lceil(p+1)/b\rceil!\le \overline \Gamma e^{p/b}(p\,!)^{1/b}$ and the fact that this last bound is also available for $\Gamma(1/b)$, we have:
\begin{eqnarray*}
\sum_{r=0}^\infty(r+1)^p \rho(r) & \le & C_0\left(C_{a,b}+2^{p-1}b^{-1}\left[a^{-(p+1)/b}+a^{-1/b}\right]\overline \Gamma e^{p/b}(p\,!)^{1/b}\right)\\
& \le & C_0\left(C_{a,b}+2^{-1}b^{-1}a^{-1/b}\overline \Gamma\right) \left(2e^{1/b}(a^{-1/b}\vee 1)\right)^p (p\,!)^{1/b}.
\end{eqnarray*}
The Lemma comes immediately by choosing the appropriated constants.
\end{proof}

\subsection{Proof of Theorem \ref{th}}\label{proofth}

The proof of Theorem \ref{th} is very similar to the one in the iid case given in \cite{DonJohnKerkPic}. In the sequel, $C$ denotes a positive real number that does not depend on $n$, $j$ nor $k$. Its value may vary from an equation to another. Let us fix $1\le p<\infty$ and consider only the cases where $1\le\pi\le p$. The cases where $\pi>p$ follow from the case $\pi=p$ applying Jensen's inequality on the error term $\E\|\widehat f_n-f\|_p^p$. Theorem \ref{th} provides the convergence rate for
$f$ belonging to the Besov Ball $\mathcal
B^s_{\pi,r}(M_1)$. In particular $f\in L^2$ and it can be written as:
$$f=\underbrace{\sum_{k=0}^{2^j_0-1}\alpha_{j_0,k}\phi_{j_0,k}}_{E_{j_0}f}+ \underbrace{\sum_{j=j_0}^{\infty}\sum_{k=0}^{2^j-1}\beta_{j,k}\psi_{j,k}}_{D_{j_0}f}.$$
We decompose the estimators $\widehat f_n$ of $f$ in the same way:
$$\widehat f_n= \underbrace{\sum_{k=0}^{2^j_0-1}\widehat\alpha_{j_0,k}
\phi_{j_0,k}}_{\widehat E_{j_0}f}+\underbrace{\sum_{j=j_0}^{j_1}
\sum_{k=0}^{2^j-1}\gamma_{\lambda_j}(\widehat\beta_{j,k})\psi_{j,k}}_{\widehat
D_{j_0}f},$$ where the $\gamma_{\lambda_j}$ denotes without distinction the soft and hard-threshold function.

Thanks to Minkowski's inequality, the risk of $\widehat f_n$ is divided in two terms:
$$\E[\|\widehat f_n -f\|_p^p]\leq
2^{p-1}\left(\underbrace{\E[\|\widehat E_{j0}f -E_{j0}f\|_p^p]}_{T_1}
+\underbrace{\E[\|\widehat D_{j_0,j_1}f -D_{j_0}f\|_p^p]}_{T_2}\right).$$\\
To study the convergence
rates of these terms the main tools are the following Lemmas given respectively in \cite{Meyer} and \cite{DonJohnKerkPic}. Here for any $p\ge1$ we denote $\|\cdot\|_{\ell_p}$ the $\ell_p$-norm defined by $\|a\|_{\ell_p}^p=\sum_i |a_i|^p$ for any sequence of real number $(a_i)_{i\ge 0}$.

\begin{lem}\label{Meyer}
Let $\delta$ denote with no distinction $\phi$ and $\psi$. For any $1\le p\le \infty$, there exists $c_1,c_2>0$ such that for all $j\ge0$ and all sequence $(a_{k})_{0\le k\le 2^j-1}$ we have
$$
c_1 2^{j(p/2-1)}\|a\|_{\ell_p}^p\le \left\|\sum_{k=0}^{2^j-1}a_{k}\delta_{j,k}\right\|_p^p \le c_2 2^{j(p/2-1)}\|a\|_{\ell_p}^p.
$$
\end{lem}

\begin{lem}\label{Donoho}
For any $1\le p< \infty$ there exists a constant $C>0$ such that for all $0\le j^-\le j^+$ and any triangular arrays $(a_{j,k})_{j^-\le j\le j^+,\,0\le k\le 2^j-1}$ we have
$$
\left\|\sum_{j=j^-}^{j^+}\sum_{k=0}^{2^j-1}a_{j,k}\delta_{j,k}\right\|_p^p \le C \begin{cases}\displaystyle \sum_{j=j^-}^{j^+} 2^{j(p/2-1)}\|a_j\|_{\ell_p}^p,&\mbox{ if }1\le p\le2,\\
\displaystyle \left(\sum_{j=j^-}^{j^+}2^{j\beta p/(p-2)}\right)^{(p-2)/2}\sum_{j=j^-}^{j^+} 2^{j(p/2-1-\beta p/2)}\|a_j\|_{\ell_p}^p,&\mbox{ if }p>2.
\end{cases}
$$
Here $\beta$ is an arbitrary real number. 
\end{lem}

\subsubsection[Bias of scale estimation]{Bias of scale estimation $T_1$}

Lemma \ref{Meyer} for $\delta=\phi$ and $a_{k}=\widehat\alpha_{j_0,k}-\alpha_{j_0,k}$ yields the existence of $C>0$ such that
$$
T_1=\E\left\|
\sum_{k=0}^{2^j_0-1}(\widehat\alpha_{j_0,k}-\alpha_{j_0,k})\phi_{j_0,k}\right\|_p^p\le C\,
2^{j_0(p/2-1)}\sum_{k=0}^{2^{j_0}-1}\E|\widehat\alpha_{j_0,k}-\alpha_{j_0,k}|^p.
$$
Thanks to Lemma \ref{useful}, the term
$T_1$ is bounded by
$$T_1\le C \,2^{j_0 p/2}n^{-p/2}.$$ Note that the
choice of $j_0$ in Theorem \ref{th} implies that the order of the bound is
$(2^{j_0}/n)^{p/2}\le Cn^{-pN/(2+2N)}$ negligible compare with $n^{-ps/(1+2s)}$ thanks to the hypothesis $N> 2s$. Note that
\begin{equation}\label{comparealpha}
\frac{\alpha_-}{\alpha_+}=1+\frac{\epsilon}{2sp\pi(1+2(s-1/\pi))},
\end{equation}
in order that if $\epsilon> 0$ then $\alpha_+<\alpha_-$, if $\epsilon< 0$ then $\alpha_+>\alpha_-$ and if $\epsilon= 0$ then $\alpha_+=\alpha_-$. We conclude that for all possible choices of $\epsilon$ the term $T_1$ is negligible.

\subsubsection[Details term]{Details term $T_2$}

The proof is based on multiple applications of Lemma \ref{Donoho} with $a_{j,k}=\gamma_{\lambda_j}(\widehat\beta_{j,k})-\beta_{j,k}$ and different $j^-$ and $j^+$. Studying the expectation of the loss, according to Lemma \ref{Donoho} and from the linearity of the expectation, the key terms of the bounds are 
$$
\E\|\gamma_{\lambda_j}(\widehat\beta_{j,k})-\beta_{j,k}\|_{\ell_p}^p=\sum_{k=0}^{2^j-1}\E|\gamma_{\lambda_j}(\widehat\beta_{j,k})-\beta_{j,k}|^p.
$$
The following Lemma gives upper bounds for the terms $\E|\gamma_{\lambda_j}(\widehat\beta_{j,k})-\beta_{j,k}|^p$ for any $j,k$:
\begin{lem}\label{boundbeta}
Under the assumptions of Theorem \ref{th}, there exists a constant $C>0$ such that for all $j_0\le j\le j_1$ and all $0\le k\le 2^j-1$ we have
$$
\E|\gamma_{\lambda_j}(\widehat\beta_{j,k})-\beta_{j,k}|^p\le C (\mbox{\bf B1}\wedge\mbox{\bf B2}\wedge\mbox{\bf B3}),
$$
where {\bf B1}=$|\beta_{j,k}|^p$, {\bf B2}=$|\beta_{j,k}|^\pi \lambda_j^{p-\pi}$ and {\bf B3}=$\lambda_{j}^p$.
\end{lem}
\begin{proof}
From the definition of the hard and soft-threshold functions $\gamma_{\lambda_j}$ we have :
\begin{equation}\label{decomp}
|\gamma_{\lambda_j}(\widehat\beta_{j,k})-\beta_{j,k}|^p\le 2^{p-1}\left(|\widehat\beta_{j,k}-\beta_{j,k}|^p+|\lambda_j|^p\right)\1_{\{|\widehat\beta_{j,k}|>\lambda_j\}}+|\beta_{j,k}|^p\1_{\{|\widehat\beta_{j,k}|\le \lambda_j\}}.
\end{equation}
The idea is to introduce the difference $|\widehat \beta_{j,k}-\beta_{j,k}|$ and bound it using Lemma \ref{useful}. More precisely, we have 
$$
\1_{\{|\widehat\beta_{j,k}|>\lambda_j\}}\le\1_{\{|\widehat\beta_{j,k}-\beta_{j,k}|\ge\lambda_j/2\}}+\1_{\{|\beta_{j,k}|\ge\lambda_j/2\}},
$$
and the expectation of the first term of the bound in \eqref{decomp} is bounded by:
\begin{multline*}
\left(\E|\widehat\beta_{j,k}-\beta_{j,k}|^{2p}\right)^{1/2}\left(\P(|\widehat\beta_{j,k}-\beta_{j,k}|\ge\lambda_j/2)\right)^{1/2}+|\lambda_j|^{p}\P(|\widehat\beta_{j,k}-\beta_{j,k}|\ge\lambda_j/2)\\+(\E|\widehat\beta_{j,k}-\beta_{j,k}|^{p}+|\lambda_j|^{p})\1_{\{|\beta_{j,k}|\ge \lambda_j/2\}}.
\end{multline*}
Inequality \eqref{Rosenthal} in Lemma \ref{useful} provides that $\E|\widehat\beta_{j,k}-\beta_{j,k}|^{p}\le C n^{-p/2}$ and consequently the expectation terms of the sum are smaller than $|\lambda_j|^{p}$. Applying \eqref{bernstdep} of Lemma \ref{useful} with $\lambda=\lambda_j$ satisfying $\lambda_j \sqrt{n}=K j$, we infer the existence of a constant $c>0$ such that 
$$
\left(\P(|\widehat\beta_{j,k}-\beta_{j,k}|\ge\lambda_j/2)\right)^{1/2}\le C2^{-cKj}
$$
from the inequality
$$
(2^{j/2}/\sqrt{n})^b\le (\log n)^{-2b-3}\le C j^{-2b-3}.
$$

Using $
\1_{\{|\widehat\beta_{j,k}|\le \lambda_j\}}\le\1_{\{|\widehat\beta_{j,k}-\beta_{j,k}|\ge\lambda_j\}}+\1_{\{|\beta_{j,k}|\le 2\lambda_j\}}$, the expectation of the second term of the bound in \eqref{decomp} is lower than
$$
|\beta_{j,k}|^p\left(\P(|\widehat\beta_{j,k}-\beta_{j,k}|\ge\lambda_j)+\1_{\{|\beta_{j,k}|\le 2\lambda_j\}}\right).
$$
As above the probability term is bounded by $2^{-cKj}$. 

Using all the previous bounds in inequality \eqref{decomp} leads to an upper bound for $\E|\gamma_{\lambda_j}(\widehat\beta_{j,k})-\beta_{j,k}|^p$ with the expression $B(p,j,k)+(\lambda_j^p+|\beta_{j,k}|^p)2^{-cKj}$ where
$$
B(p,j,k)=\lambda_j^p\1_{\{|\beta_{j,k}|\ge \lambda_j\}}+|\beta_{j,k}|^p\1_{\{|\beta_{j,k}|\le \lambda_{j}\}}.
$$
We investigate three ways for bounding $B(p,j,k)$, using the fact that $\1_{\{a\le b\}}\le (b/a)^\alpha$ for all $a,b>0$ and all $\alpha\ge 0$:
\begin{description}
\item[B1] $B(p,j,k)\le C |\beta_{j,k}|^p$ using the indicator in the first term of the sum,
\vspace{16pt}

\item[B2] $B(p,j,k)\le C |\beta_{j,k}|^\pi \lambda_j^{p-\pi}$ with $a=\lambda_j$, $b=|\beta_{j,k}|$ and $\alpha=\pi$ in the first term and $a=|\beta_{j,k}|$, $b=\lambda_j$ and $\alpha=p-\pi$ in the second term of the sum,
\vspace{16pt}

\item[B3] $B(p,j,k)\le C \lambda_{j}^p$ using the indicator in the second term of the sum.
\end{description}
Thanks to a suitable choice of $K$ it is obvious that $(\lambda_j^p+|\beta_{j,k}|^p)2^{-cKj}$ is smaller than $(\mbox{\bf B1}\wedge\mbox{\bf B2}\wedge\mbox{\bf B3})$. The result of Lemma \ref{boundbeta} follows. 
\end{proof}

The end of the proof of Theorem \ref{th} is based on successive uses of {\bf B1}, {\bf B2} and {\bf B3} depending on the resolution levels $j$.

Let us first consider the highest multi resolution levels $j>j_+$ where {\bf B1} is the most efficient bound given in Lemma \ref{useful}. The integer $j_+$ will be fixed later. From Minkowski's inequality we obtain the following bound for $T_2$, up to a constant:
$$ \underbrace{\E\Big\|\sum_{j=j_0}^{j_+}\sum_{k=0}^{2^{j}-1}(\gamma_{\lambda_j}(\widehat\beta_{j,k})-\beta_{j,k})\psi_{j,k}\Big\|_p^p}_{T_{21}}+\underbrace{\E\Big\|\sum_{ j_+<j\le j_1}\sum_{k=0}^{2^{j}-1}\gamma_{\lambda_j}(\widehat\beta_{j,k})\psi_{j,k}-\sum_{ j_+<j}\sum_{k=0}^{2^{j}-1}\beta_{j,k}\psi_{j,k}\Big\|_p^p}_{T_{22}}.
$$
The term $T_{22}$ is bounded applying Lemma \ref{Donoho} with $j^-=j_++1$, $j^+=\infty$ and $a_{j,k}=\gamma_{\lambda_j}(\widehat\beta_{j,k})-\beta_{j,k}$ for $j_+< j\le j_1$ and $0\le k\le 2^j-1$, $a_{j,k}=-\beta_{j,k}$ for $j>j_1$ and $0\le k\le 2^j-1$. Using Lemma \ref{boundbeta} with {\bf B1}, we achieve that $\E\|a_j\|_{\ell_p}^p\le C \|\beta_j\|_{\ell_p}^p$ for all $j>j_+$ and we have:
$$
T_{22}\le C \begin{cases}\displaystyle \sum_{j>j_+} 2^{j(p/2-1)}\|\beta_j\|_{\ell_p}^p,&\mbox{ if }1\le p\le2,\\
\displaystyle \left(\sum_{j>j_+}2^{j\beta p/(p-2)}\right)^{(p-2)/2}\sum_{j>j_+} 2^{j(p/2-1-\beta p/2)}\|\beta_j\|_{\ell_p}^p,&\mbox{ if }p>2.
\end{cases}
$$
The Sobolev inclusion $\mathcal B^s_{\pi,r}\subset
\mathcal B^{s'}_{p,r}\subset
\mathcal B^{s'}_{p,\infty}$ with
$s'=s-1/\pi+1/p$ leads us to choose $\beta=-2s'$. Noting that $\sum_{j>j_+}2^{-j2s'p/(p-2)}\le C 2^{-j_+s'p2/(p-2)}$ we obtain the inequalities:
$$
T_{22}\le C2^{-j_+s'p}\sum_{j=j_+}^{\infty} 2^{j(s'p+p/2-1)}\|\beta_j\|_{\ell_p}^p\le C\|f\|^p_{s',p,\infty}2^{-j_+s'p}.
$$
We can choose $j_+$ as the largest integer such as
$$
2^{j_+}\le \left(\frac{n}{\log n}\right)^\frac{\alpha}{s'}.
$$
We have to check that $j_+\le j_1$ for $n$ sufficiently large, {\it i.e.} $s'>\alpha$. When $\epsilon\leq 0$ we have $\alpha=s'/(1+2(s-1/\pi))$ and $\alpha< s'$ because $s> 1/\pi$ by hypothesis. When $\epsilon> 0$, equality \eqref{comparealpha} implies that $\alpha_+\le\alpha_-< s'$ and then obviously $s'>\alpha$ for all the possible values of $\epsilon$. With this choice of $j_+$, the rate of convergence of $T_{22}$ is then the one stated in the Theorem.

We study the convergence of $T_{21}$ using the bounds {\bf B2} and {\bf B3} given in Lemma \ref{boundbeta}. Let us investigate further these two cases in order to compare their efficiency. On the one hand, using {\bf B2} and noting that $\|f\|_{s,\pi,\infty}\le \|f\|_{s,\pi,r}<\infty$, we obtain the inequality:
$$
\E\|\gamma_{\lambda_j}(\widehat\beta_{j.})-\beta_{j.}\|_{\ell_p}^p\le C \lambda_j^{p-\pi}\|\beta_j\|_{\ell_\pi}^\pi\le C \lambda_j^{p-\pi}2^{-j(s\pi+\pi/2-1)}.
$$
On the other hand, using {\bf B3} we have directly
$$
\E\|\gamma_{\lambda_j}(\widehat\beta_{j.})-\beta_{j.}\|_{\ell_p}^p\le C 2^j\lambda_j^{p}.
$$
The rate of these two upper bounds are equivalent when $\lambda_j^2=2^{-j(2s+1)}$. Replacing $\lambda_j$ by its value,
it follows that {\bf B2} is more efficient for $j> j_-$ with
$j_-$ the largest integer such that 
$$ 
2^{j_-}\le\left(\frac{n}{\log n}\right)^\frac{1}{1+2s}.
$$
We check that $j_-\le j_+$ for all $n$ and all possible values of $\epsilon$ as $\alpha=s/(1+2s)>s'/(1+2s)$ if $\epsilon\ge 0$ and $s'/\alpha=1+2s-2/\pi<1+2s$ if $\epsilon\le 0$.

We decompose again $T_{21}$ using Minkowski's inequality. It gives the following upper bound, up to a constant:
$$ \underbrace{\E\Big\|\sum_{j=j_0}^{j_-}\sum_{k=0}^{2^{j}-1}(\gamma_{\lambda_j}(\widehat\beta_{j,k})-\beta_{j,k})\psi_{j,k}\Big\|_p^p}_{T_{211}}+\underbrace{\E\Big\|\sum_{j_-<j\le j_+}\sum_{k=0}^{2^{j}-1}(\gamma_{\lambda_j}(\widehat\beta_{j,k})-\beta_{j,k})\psi_{j,k}\Big\|_p^p}_{T_{212}}.
$$
We control the term $T_{211}$ using {\bf B3} according to the discussion above and applying Lemma \ref{Donoho} with $\lambda_j^2\le C\log n/n$:
$$
T_{211}\le C(\log n/n)^{p/2} \begin{cases}\displaystyle \sum_{j=j_0}^{j_1} 2^{jp/2},&\mbox{ if }1\le p\le2,\\
\displaystyle \left(\sum_{j_0}^{j_-}2^{j\beta p/(p-2)}\right)^{(p-2)/2}\sum_{j=j_0}^{j_1} 2^{jp(1-\beta)/2},&\mbox{ if }p>2.
\end{cases}
$$
Let us choose $\beta=1/2$ in order to obtain the inequalities
$$
T_{211}\le C(\log n/n)^{-p/2}2^{j_-p/2}\le C (\log n/n)^{p(1-1/(2s+1))/2}\le C (\log n/n)^{p\alpha}\qquad\mbox{ if }\epsilon \ge 0.
$$
This term is negligible if $\epsilon<0$ using \eqref{comparealpha}.

To conclude, it remains to bound $T_{212}$ using Lemma \ref{Donoho} and {\bf B2}. We use that $\lambda_j^2\le C\log n/n$ and we let appear the symbol $\epsilon=s\pi-(p-\pi)/2$
$$
T_{211}\le C(\log n/n)^{(p-\pi)/2} \begin{cases}\displaystyle \sum_{j_-<j\le j_+} 2^{-j\epsilon},&\mbox{ if }1\le p\le2,\\
\displaystyle \left(\sum_{j_-<j\le j_+}2^{j\beta p/(p-2)}\right)^{(p-2)/2}\sum_{j_-<j\le j_+} 2^{-j(\epsilon+p\beta/2)},&\mbox{ if }p>2.
\end{cases}
$$
From now we have to distinguish the cases where $\epsilon\neq0$ to these where $\epsilon=0$. 
\begin{itemize}
\item{If $\epsilon\neq 0$:}
Let us take $\beta=-\epsilon/p$. Then if $\epsilon<0$, we obtain the inequalities:
$$
T_{211}\le C(\log n/n)^{(p-\pi)/2}2^{-j_+\epsilon}\le C (\log n/n)^{\alpha\epsilon/s'+(p-\pi)/2},
$$
and we conclude using the equality $\alpha\epsilon/s'+(p-\pi)/2=p\alpha$. If $\epsilon>0$, then
$$
T_{211}\le C(\log n/n)^{(p-\pi)/2}2^{-j_-\epsilon}\le C (\log n/n)^{p(1-1/(1+2s))/2}\le C (\log n/n)^{p\alpha}.
$$

\item{If $\epsilon=0$:}
Then $p/2=\pi s+\pi/2$ and as $s\pi>1$ we deduce that $p>2$. Moreover we notice that $(p-\pi)/2=\alpha p$ and thus
$$
2^{j(p/2-1)}\E\|\gamma_{\lambda_j}(\widehat\beta_{j.})-\beta_{j}\|_{\ell_p}^p\le C \lambda_j^{(p-\pi)/2}\|\beta_j\|_{\ell_\pi}^\pi\le C \lambda_j^{2\alpha p}2^{j(\pi s+\pi/2-1)}\|\beta_{j}\|_{\ell_\pi}^\pi.
$$
Let us denote 
$
t_j=2^{j(\pi s+\pi/2-1)}\sum_{k=0}^{2^{j}-1}\|\beta_{j}\|_{\ell_\pi}^\pi
$. Notice that from the definition of the Besov norms we have $\sum_j t_j^{r/\pi}\le C \|f\|_{s,\pi,r}^r$. From Lemma \ref{Donoho} with $\beta=0$ we achieve
$$
T_{211}\le C(\log n/n)^{\alpha p} j_+^{p/2-1}\sum_{j_-<j\le j_+} t_j,\mbox{ as }p>2.
$$
If $r\le \pi$ then $\sum_{j_-<j\le j_+} t_j\le C$ and the result of Theorem \ref{th} follows. If $\pi< r$ we use Holder's inequality at the powers $r/\pi$ and $r/(r-\pi)$:
$$
\sum_{j=j_0}^{j_+}t_j\le C\|f\|_{s,\pi,r}^\pi \left(\sum_{j=j_0}^{j_+} t_j^{r/(r-\pi)}\right)^{1-\pi/r}\le C j_+^{1-\pi/r}\le C(\log n)^{1-\pi/r}.
$$ 
\end{itemize}

\subsection{Proofs of results given in Sections \ref{simus} and \ref{model}}\label{proofsimus}
In this Section are collected the proof of Proposition~\ref{phiD}, Lemma \ref{tildephi}, Proposition \ref{etaD} and  Propostion \ref{contrex}. Denoting with no distinction $\psi_{j,k}(x)-\beta_{j,k}$ and $\phi_{j,k}(x)-\alpha_{j,k}$ as $\widetilde\delta_{j,k}(x)$ for any $j,k$, we collect here some inequalities useful in this Section: $\E|\widetilde\delta_{j,k}(X_0)|\le2\|f\|_\infty\|\delta\|_1
2^{-j/2}$, $\E|\widetilde\delta_{j,k}(X_0)|^2\le \|f\|_\infty$, $\|\widetilde\delta_{j,k}\|_\infty\le
2\|\delta\|_\infty2^{j/2}$, $\Lip \widetilde\delta_{j,k}\le
\Lip\widetilde\delta 2^{3j/2}$ and $\|\widetilde\delta_{j,k}\|_{BV}\le(\|\widetilde\delta\|_\infty+A\Lip\widetilde\delta)2^{j/2+1}$ for all $j\ge1$. The last assertion comes from the fact that $\widetilde\delta_{j,k}$ is a bounded Lipschitz function supported by $[(-A+k)2^{-j},(A+k)2^{-j}]$. 

\begin{proof}[Proof of Proposition~\ref{phiD}]
As for any $j,k$ the function $\widetilde\delta_{j,k}$ has bounded variations we have
$$
\left| \cov \left(\widetilde\delta_{j,k}(X_{s_1})\cdots \widetilde\delta_{j,k}(X_{s_u}), \widetilde\delta_{j,k}(X_{s_{u +
1}})\cdots \widetilde\delta_{j,k}( X_{s_{u+v}}) \right) \right| \leq v
\mathbb{E}\left|\widetilde\delta_{j,k}(X_{s_1})\cdots \widetilde\delta_{j,k}(X_{s_u})\right| \|\widetilde\delta_{j,k}\|_{BV}^{v}
\widetilde \phi_v(r).
$$
Noticing that $\|\widetilde\delta_{j,k}\|_\infty\le \|\widetilde\delta_{j,k}\|_{BV}$, it follows
\begin{eqnarray*}
C_{u,v}^{j,k}(r)&\le& v \|\widetilde\delta_{j,k}\|_{BV}^{u+v-2}\|\widetilde\delta_{j,k}\|_{BV}\E|\widetilde\delta_{j,k}(X_0)| \widetilde \phi_v(r)\\
&\le& (u+v+uv)
(2^{j/2+1}c(\|\delta\|_\infty+A\Lip\delta))^{u+v-1}2\|f\|_\infty\|\delta\|_1 2^{-j/2}\exp(-ar^b).
\end{eqnarray*}
Then Proposition~\ref{phiD} is proved.
\end{proof}

\begin{proof}[Proof of Lemma \ref{tildephi}]
The proof is very close to the one given in \cite{Dedecker2007a}. First notice that \eqref{condC'} for $r=0$ implies that $c\ge 1=\sup_{g\in BV_1}\|g\|_{BV}$. From \eqref{condan} we infer that $\widetilde \phi(\sigma(X_0),X_r)\le \exp(ar^{-b})$ applying Lemma 4 of
\cite{DedeckerPrieur} on \eqref{condan}. From the Markov property we get $\widetilde \phi(\sigma(\{X_j,j\le0\}),X_r)=\widetilde \phi(\sigma(X_0),X_r)$. Now for any $\ell\le 1$, for all $r\le i_1\le \cdots\le i_\ell$ consider any
$g_{i_j}\in BV_1$ for $1\le j\le \ell$.

Let us prove that we can restrict ourselves to the $\ell$-uplets satisfying $i_1<\cdots<i_\ell$. On the one hand, if $i_j=i_{j'}=r$ then we have $\|\E(g_r(X_r)g_r(X_r)|X_{i_1}=\cdot)\|_{BV}\le c\|g_r^2\|_{BV}$ from \eqref{condC'}. As from assumption $\|g_r\|_\infty\le \|g_r\|_{BV}\le 1$ then $\|g_r^2\|_{BV}\le \|g_r\|_{BV}+\|g_r\|_{BV}\le 2$ and we achieve the bound $\|\E(g_r(X_r)g_r(X_r)|X_{i_1}=\cdot)\|_{BV}\le 2c$. On the other hand, if $i_j<i_{j'}$ then, from the Markov property, we get the equation
$$
\|\E(g_{i_j}(X_{i_j})g_{i_{j'}}(X_{i_{j'}})|X_{i_1}=\cdot)\|_{BV}=\|\E(g_{i_j}(X_{i_j})\E(g_{i_{j'}}(X_{i_{j'}})|X_{i_j})|X_{i_1}=\cdot)\|_{BV}.
$$
We proceed in two steps. Firstly from \eqref{condC'}, we infer thet $x\mapsto g_{i_j,i_{j'}}(x)=\E(g_{i_{j'}}(X_{i_{j'}})|X_{i_{j}}=x)$ has variations bounded by $c$. Notice also that $\|g_{i_j,i_{j'}}\|_\infty\le \|g_{i_{j'}}\|_\infty\le 1$ and therefore we deduce the bound $\|g_{i_j}g_{i_j,i_{j'}}\|_{BV}\le\|g_{i_j}\|_{BV}+ \|g_{i_j,i_{j'}}\|_{BV}\le 1+c$. Secondly, using \eqref{condC'} on $g_{i_j,i_{j'}}/\|g_{i_j,i_{j'}}\|_{BV}$, we infer that $\|\E(g_{i_j}(X_{i_j})g_{i_{j'}}(X_{i_{j'}})|X_{i_1}=\cdot)\|_{BV}\le c(1+c)\le c+c^2$.

As $c\ge 1$, this bound is larger than the one in the case $i_j=i_{j'}$. Now, by straightforward recurrences in the worst cases $i_1<\cdots<i_\ell$, we get that
$$
\|\E(g_{i_1}(X_{i_1})\cdots g_{i_\ell}(X_{i_\ell})|X_{i_1}=\cdot)\|_{BV}\le c+\cdots + c^\ell\le \ell \,c^\ell~~\mbox{ as }c\le 1.
$$
Then, denoting $g_{i_1,\ldots,c_\ell}(x)=\E(g_{i_1}(X_{i_1})\cdots g_{i_\ell}(X_{i_\ell})|X_{i_1}=x)$ for all $x$, we have almost surely the equation
$$
\E(g_{i_1}(X_{i_1})\cdots g_{i_\ell}(X_{i_\ell})|X_0)-\E(g_{i_1}(X_{i_1})\cdots g_{i_\ell}(X_{i_\ell}))=\E(g_{i_1,\ldots,c_\ell}(X_{i_1})|X_0)-\E(g_{i_1,\ldots,c_\ell}(X{i_1}))
$$
From the definition of coefficients $\widetilde \phi$ we get
$$
\|\E(g_{i_1,\ldots,c_\ell}(X_{i_1})|X_0)-\E(g_{i_1,\ldots,c_\ell}(X_{i_1}))\|_\infty\le \ell c^\ell\widetilde \phi(\sigma(\{X_j,j\le 0\}),X_r).
$$
For all $1\le \ell \le v$ this bound holds uniformly for all $g_{i_j}\in BV_1$, $1\le j\le \ell$, using the definition of $\widetilde \phi_v(r)$ we conclude that $\widetilde \phi_v(r)\le c^v\widetilde \phi (\sigma(\{X_j,j\le 0\}),X_r)\le c^v\exp(-ar^b)$.
\end{proof}

\begin{proof}[Proof of Proposition \ref{etaD}]
We use the direct bound given in the
proof of Lemma~1 of \cite{RagacheWintenberger} under {\bf (J)}:
$$
\left| \mbox{Cov} \left(\widetilde\delta_{j,k}(X_{i_1})\cdots \widetilde\delta_{j,k}(X_{i_u}), \widetilde\delta_{j,k}(X_{i_{u +
1}})\cdots \widetilde\delta_{j,k}( X_{i_{u+v}}) \right) \right| \leq
2^3(2^{j/2}2\|\delta\|_\infty)^{u+v-2}\gamma(r)
$$
with $\gamma(r)=\mathbb{E}|\widetilde\delta_{j,k}(X_0)\widetilde\delta_{j,k}(X_r)|\vee
(\mathbb{E}|\widetilde\delta_{j,k}(X_0)|)^2\le (\|f_r\|_\infty\vee 2\|f\|_\infty)\|\delta\|_1^2
2^{-j}$. Noticing that $\gamma\wedge \beta\le \gamma^{1/4} \beta^{3/4}$ for any positive numbers $\gamma$ and $\beta$, we combine the two bounds on the covariance terms and we infer that \eqref{cpr} is satisfied with $M_2=2\|\delta\|_\infty$ and
$$
\rho(r)= 2^{9/2} \|\delta\|_1^{3/2}(\Lip \delta)^{1/2}(\|f_r\|_\infty\vee \|f\|_\infty)^{3/4}\lambda(r)^{1/4}.
$$
Assumptions of Proposition~\ref{etaD} on respective decrease and increase rates of $\lambda(r)$ and $\|f_r\|_\infty$ yield the existence of $C_0>0$ such that $\rho(r)\le C_0\exp(-a'r^{b'}/4)$.
\end{proof}
\begin{proof}[Proof of Propostion \ref{contrex}]
We give the proof for $\widehat f_n$ but it also holds for $\widehat f_n^{STCV}$. We begin with recalling the result of Corollary 7.1 in \cite{Gouezel2004}:
\begin{lem}[Gou\"{e}zel, 2004]\label{lemgouezel}
For any Lipschitz function $\delta_1$, bounded measurable function $\delta_2$ such that $\delta_1,\delta_2=0$ in a neighborhood of $0$, then for any $0<\alpha'<1$, there exists some constant $C>0$ such that
\begin{equation}\label{equivcov}
\cov(\delta_1(X_0),\delta_2(X_r))\sim C \int\delta_1(x)dx\int\delta_2(x)dx\; r^{1-1/\alpha'} \mbox{ when }r\to\infty.
\end{equation}
\end{lem}
We use the decomposition of $f\in \L^2$ in the orthogonal basis of wavelets functions  and we obtain
\begin{eqnarray*}
\E(\|\widehat f_n-f\|_2^2)&\ge& \E\|\sum_{k\in S_{j_0}}(\widehat \alpha_{j_0-k}-\alpha_{j_0-k})\phi_{j_0,k}\|_2^2\\
&\ge&\sum_{k\in S_{j_0}}\E(\widehat \alpha_{j_0-k}-\alpha_{j_0-k})^2.
\end{eqnarray*}
If we develop $\E(\widehat \alpha_{j_0-k}-\alpha_{j_0-k})^2$ using the covariance terms and denoting $\widetilde\delta_{j,k}(x)=\phi_{j,k}(x)-\alpha_{j,k}$ for $x\in[0.01,1]$ and null elsewhere, it comes:
$$
\E(\widehat \alpha_{j_0-k}-\alpha_{j_0-k})^2=\frac{1}{n}\E(\widetilde\delta_{j,k}(X_i)^2)+2\sum_{r=1}^{n-1}\frac{n-r}{n^2}\cov(\widetilde\delta_{j_0,k}(X_0),\widetilde\delta_{j_0,k}(X_r)).
$$

We want to apply Lemma \ref{lemgouezel} with $\delta_1=\delta_2=\widetilde\delta_{j_0,k}$. We check easily the assumptions of this Lemma because of the definition of $\,\widetilde\delta_{j_0,k}$, resulting from the fact that we estimate the density on $[0.01,1]$. Moreover $\int \phi_{j_0,k}=2^{-j_0/2}\int\phi$ with $\int\phi>0$ from assumption and then the covariance terms $\cov(\widetilde\delta_{j_0,k}(X_0),\widetilde\delta_{j_0,k}(X_r))$ are equivalent to $C2^{-j_0}\,r^{1-1/\alpha'}$ for some $C>0$ as $r$ goes to infinity. 

Let $n_0$ be such that for all $n\geq n_0$, $u_{r,n}=\frac{n-r}{n^2}\cov(\widetilde\delta_{j_0,k}(X_0),\widetilde\delta_{j_0,k}(X_r))$ is nonnegative. For some $m>m'>2$, we decompose the sum of covariance terms, for $n$ sufficiently large, in four sums: 
$$\sum_{r=1}^{n-1}\frac{n-r}{n^2}\cov(\widetilde\delta_{j_0,k}(X_0),\widetilde\delta_{j_0,k}(X_r))= \sum_{r=1}^{n_0}u_{r,n} 
+\sum_{r=n_0}^{[n/m]}u_{r,n}+\sum_{r=[n/m]}^{[n/m']}u_{r,n}+\sum_{r=[n/m']}^{n-1}u_{r,n},$$
where $[a]$ denotes the integer part of $a$.
The first term goes to $0$ with rate $n$. Then, by definition of $n_0$, the second and the last terms are nonnegative. Concerning the sum $\sum_{r=[n/m]}^{[n/m']}u_{r,n}$, the summands $u_{r,n}$ are all larger than $(m'-2)/(m'n)\cov(\widetilde\delta_{j_0,k}(X_0),\widetilde\delta_{j_0,k}(X_r))$ that is equivalent to $C2^{-j_0}n^{-1}r^{1-1/\alpha'}$ as $n\to \infty$. The minimax rate $\alpha$ is such that $\alpha<1/2$ and then by hypothesis $\alpha'\geq 1/(2\alpha+1)>1/2$. Consequently, when $n$ goes to infinity, the sum $\sum_{r=[n/m]}^{[n/m']}u_{r,n}$ is larger than a partial sum equivalent to $C2^{-j_0}n^{1-1/\alpha'}$ for some $C>0$ as $n\to \infty$. As we assume $2\alpha\ge1/\alpha'-1$ we obtain the existence of some $C>0$ such that  
$$n^{2\alpha}2^{j_0}\sum_{r=1}^{n-1}\frac{n-r}{n^2}\cov(\widetilde\delta_{j_0,k}(X_0),\widetilde\delta_{j_0,k}(X_r))\ge C, \mbox{ for $n$ sufficiently large}.$$
Collecting these facts and using that $|S_{j_0}|$ is up to a constant equal to $2^{j_0}$, we obtain that $n^{2\alpha}\E(\|\widehat f_n-f\|_2^2)$  is larger than some positive constant and the result of Proposition \ref{contrex} is proved.
\end{proof}
\noindent {\bf Acknowledgements.} We thank Karine Tribouley and Vincent Rivoirard for their helpful advices  to one of the author. We also 
thank referees for having pointed out lapses in the first version of this paper.


\bibliographystyle{acm}
\bibliography{biblio}

\end{document}